\documentclass[12pt]{amsart}

\usepackage{xcolor,mathbbol,amsmath,amssymb,amsxtra,anysize,times,tikz}
\usetikzlibrary{matrix,arrows,decorations.pathmorphing}
\usepackage{tikz-cd}
\usepackage{comment}
\usepackage[all]{xy}
\usepackage{stackengine}
\usepackage{mathrsfs}

\numberwithin{equation}{subsection}

\stackMath
\newcommand\tsup[2][2]{%
 \def\useanchorwidth{T}%
  \ifnum#1>1%
    \stackon[-.5pt]{\tsup[\numexpr#1-1\relax]{#2}}{\scriptscriptstyle\sim}%
  \else%
    \stackon[.5pt]{#2}{\scriptscriptstyle\sim}%
  \fi%
}

\newtheorem{lemma}{Lemma}[section]
\newtheorem{theorem}[lemma]{Theorem}%[section]
\newtheorem{corollary}[lemma]{Corollary}%[section]
\theoremstyle{definition}   %\MV{seeing how it looks with this style}
\newtheorem{example}[lemma]{Example}%[section]
\newtheorem{remark}[lemma]{Remark}%[section]
\newtheorem{definition}[lemma]{Definition}%[section] 
 
\newcommand{\dnorm}[1]{\left\lVert#1\right\rVert_d}

\newcommand{\Unzf}{\mathfrak{U}_0(\mathfrak{gl}_n[t,t^{-1}])}
\newcommand{\Unz}{\mathfrak{U}_n}
\newcommand{\Snd}{\mathbb{S}^{\operatorname{aff}}_0(n,d)}
\newcommand{\Fnd}{\mathcal{F}_n^d}
\newcommand{\Fmu}{\mathcal{F}_{\mu}}

\DeclareMathOperator{\spn}{span}
\DeclareMathOperator{\supp}{supp}
\DeclareMathOperator{\id}{id}

\title{Equivariant K-theory of the variety of partial flags}
\author{Sergey Arkhipov}
\address{Matematisk Institut, Aarhus Universitet, Ny Munkegade, DK-8000, Aarhus C, Denmark}
\email{hippie@math.au.dk}

\author{Mikhail Mazin}
\address{Mathematics Department, Kansas State University, 138 Cardwell Hall, Manhattan, KS 66506}
\email{mmazin@math.ksu.edu}

\begin{document}

\begin{abstract}
We use Drinfeld style generators and relations to define an algebra $\mathfrak{U}_n$ which is a ``$q=0$'' version of the affine quantum group of $\mathfrak{gl}_n.$ We then use the convolution product on the equivariant $K$-theory of varieties of pairs of partial flags in a $d$-dimensional vector space $V$ to define affine $0$-Schur algebras ${\mathbb S}_0^{\operatorname{aff}}(n,d)$ and to prove that for every $d$ there exists a surjective homomorphism from $\mathfrak{U}_n$ to ${\mathbb S}_0^{\operatorname{aff}}(n,d).$
\end{abstract}

\maketitle

\section*{Introduction}
Our paper started as an attempt to study certain degenerations of quantum groups and of quantum Schur algebras. Different versions of the famous Schur-Weyl correspondence were the main inspiration for our study. Let us recall various levels of quantum Schur-Weyl correspondence from a geometric perspective. 
\subsection*{Finite sets: quantum Schur-Weyl correspondence.}
 Let us begin with the original combinatorial setting for the finite Hecke and Schur algebras. Consider a finite field $\mathbb{k}$ with $q=p^r$ elements. 
Denote the set of full flags in the vector space $V=\mathbb{k}^d$ by ${\mathcal B}={\mathcal B}(\mathbb{k})$. Notice that $G:=Gl(V)$ acts on the set ${\mathcal B}\times{\mathcal B}$ diagonally.  Consider the $\mathbb{C}$-vector space of $G$-invariant functions on ${\mathcal B}\times{\mathcal B}$ denoted by ${\mathcal H}_q^{\operatorname{fin}}$.
It is known and will be recalled below in various contexts that ${\mathcal H}_q^{\operatorname{fin}}$ is equipped with an associative algebra structure given by convolution. The obtained algebra is called the finite Hecke algebra for the parameter $q$.

Fix another integer number $n$.  Denote the set of partial $n$-step flags in the vector space $V=\mathbb{k}^d$ by ${\mathcal F}_n^d= {\mathcal F}_n^d(\mathbb{k})$. Notice that the group $G=Gl(V)$ acts on the set ${\mathcal F}_n^d\times{\mathcal F}_n^d$ diagonally.  Consider the $\mathbb{C}$-vector space of $G$-invariant functions on ${\mathcal F}_n^d\times{\mathcal F}_n^d$ denoted by ${\mathbb S}_q^{\operatorname{fin}}={\mathbb S}_q^{\operatorname{fin}}(n,d)$. 
Again, the standard construction of convolution makes ${\mathbb S}_q^{\operatorname{fin}}(n,d)$
into an associative algebra called the quantum finite Schur algebra.

The third geometric datum of the picture is given by the set ${\mathcal B}\times{\mathcal F}_n^d$, with the diagonal action of $G$. The vector space of $G$-invariant functions forms a 
${\mathbb H}_q^{\operatorname{fin}}$-${\mathbb S}_q^{\operatorname{fin}}(n,d)$-bimodule. It is the key ingredient to prove the double centralizer property for the Hecke and Schur algebras known as the quantum finite Schur-Weyl correspondence.

%\textcolor{red}{Misha: this is very classical, but what are good references for the previous five paragraphs? Perhaps a textbook...}

In their groundbreaking work, Beilinson, Lusztig and MacPherson constructed a surjective map  from the quantum group $\mathfrak{U}_q(\mathfrak{gl}_n)$ to the algebra ${\mathbb S}_q^{\operatorname{fin}}(n,d)$ (see \cite{BLM90}). The map is explicitly given in terms of quantum Chevalley generators, in particular, the elements $E_i$ (resp. $F_i$) are given by characteristic functions of certain closed orbits in ${\mathcal F}_n^d\times{\mathcal F}_n^d$ that can be viewed as elementary building blocks for the convolution algebra. 
Furthermore, the Beilinson, Lusztig, and MacPherson developed a stabilization procedure for the structure constants in the algebra ${\mathbb S}_q^{\operatorname{fin}}(n,d)$, so that in the limit one obtains the quantum group $\mathfrak{U}_q(\mathfrak{gl}_n)$. This first step to relate quantum groups to the geometry of partial flags started a variety of generalizations.

\subsection*{Finite sets: setting $q=0$.} The dependence of the structure constants in the Beilinson-Lusztig-MacPherson construction on the number $q$ of elements in the base field $\mathbb{k}$ is polynomial. The $0$-Schur algebra ${\mathbb S}_0^{\operatorname{fin}}(n,d)$ is then obtained by specializing at $q=0.$ Jensen and Su in \cite{JS15} showed that the algebra ${\mathbb S}_0^{\operatorname{fin}}(n,d)$ can be obtained geometrically by replacing the convolution product in the Beilinson-Lusztig-MacPherson construction with the \textbf{generic multiplication} of orbits (see Section $6$ in \cite{JS15}). Applying a similar approach to the complete flags Jensen and Su also obtained a geometric realization of the $0$-Hecke algebra ${\mathbb H}_0^{\operatorname{fin}}$ (see Section $10$ in \cite{JS15}). See also \cite{S10,JS15,JSY16,JSY20} for further results in this direction.
%considered the same finite sets ${\mathcal{B}}(\mathbb{k})$ and ${\mathcal{F}}_n^d(\mathbb{k})$ and the same $\mathbb{C}$-vector spaces of $G$-invariant functions but changed the definition of the product: namely they defined formally what it means to take only the ``principal part'' for the products in the $q$-expansion. 

%This way they obtained the algebras ${\mathbb H}_0^{\operatorname{fin}}$ and ${\mathbb S}_0^{\operatorname{fin}}(n,d)$ respectively. Notice that these algebras do not depend on $q$. It turns out that ${\mathbb H}_0^{\operatorname{fin}}$ is obtained from ${\mathbb H}_q^{\operatorname{fin}}$ by literally substituting $q=0$ in the formulas for the product (in a natural basis). The situation on the Schur side is more interesting.

Krob and Thibon (\cite{KT99}) defined a version of the algebra $\mathfrak{U}_q(\mathfrak{gl}_n)$ for $q=0$. Notice that the relations on the standard quantum Chevalley generators do not allow to substitute $0$ literally. The authors performed a certain renormalization to define the finite 0-quantum group $\mathfrak{U}_0(\mathfrak{gl}_n)$.
The relations on quantum Chevalley generators degenerate, in particular, quantum Serre relations stay cubic but become binomial, known under the name of {\em plactic relations} (see \ref{subsection: plactic}). 
Jensen, Su, and Yang obtained the same cubic binomial relations in their study of the $0$-Schur algebra ${\mathbb S}_0^{\operatorname{fin}}(n,d).$ %(\cite{JS15,JSY16,JSY20}).  

\subsection*{Introducing geometry: quantum affine Schur-Weyl correspondence.}
Next we replace the finite field $k$ by $\mathbb{C}$. We consider the vector space  $V:=\mathbb{C}^d$. The sets $\mathcal B$ and ${\mathcal F}_n^d$ become algebraic varieties over complex numbers. 

Kazhdan, Lusztig, Ginzburg, and Vasserot considered the Steinberg varieties on the Weyl and Schur sides $St$ and $St_n^d$ respectively (see \cite{GV93,KL80,KL87,V98}). Those are $G:=Gl(d,\mathbb{C})$-invariant Lagrangian subvarieties in 
$T^*({\mathcal B})^2$ and $T^*({\mathcal F_n^d})^2$ respectively. Notice that the multiplicative group $\mathbb{C}^*$ acts on them by dilations of cotangent directions.
The complexified $G\times\mathbb{C}^*$-equivariant K-group of $St$ (resp, of $St_n^d)$ is denoted by 
${\mathbb H}_q^{\operatorname{aff}}$ (resp. by ${\mathbb S}_q^{\operatorname{aff}}(n,d)$).
One modifies the convolution product construction for the K-group case. The obtained algebras are called the  affine Hecke algebra and the quantum affine Schur algebra respectively. 

Kazhdan, Lusztig, and Ginzburg proved that the former is indeed isomorphic to the Hecke algebra for the affine Weyl group, for the parameter $q$ coming from the $\mathbb{C}^*$ action. 
Ginzburg and Vasserot proved that the quantum loop algebra $\mathfrak{U}_q(\mathfrak{gl}_n[t,t^{-1}])$, in the Drinfeld realization, maps surjectively onto ${\mathbb S}_q^{\operatorname{aff}}(n,d)$.
Finally an immediate geometric modification of the bimodule construction provides a double centralizer property for the algebras ${\mathbb S}_q^{\operatorname{aff}}(n,d)$ and
${\mathbb H}_q^{\operatorname{aff}}$ known as affine quantum Schur-Weyl correspondence.

\subsection*{Introducing geometry: Weyl side for $q=0$.} A modification of the affine Weyl side is due to Kostant and Kumar in \cite{KK90}. Ignore the cotangent directions and consider the complexified  $G$-equivariant K-group of the variety $\mathcal B\times\mathcal B$, equipped with the convolution product. 
The authors proved that the obtained algebra ${\mathbb H}_0^{\operatorname{aff}}$ is isomorphic to the algebra obtained by substitution $q=0$ into the definition of the affine Hecke algebra. The algebra is called the affine $0$-Hecke algebra. It is the affinization of the algebra considered by Jensen, Su, and Yang.

\subsection*{The main objective: Schur side for ${q=0}$ and the affine quantum group.} Notice that, similarly to the case of $\mathfrak{U}_q(\mathfrak{gl}_n)$, the relations do not allow a direct specialization of the parameter $q$ to $0$.
One goal of the present paper is to define an algebra $\Unzf$ by Drinfeld style generators and relations. Then we interpret those generators and relations geometrically via equivariant K-groups of the varieties of partial flags.

Following the case studied by Kostant and Kumar, we consider the $G$-equivariant $K$-group of
the complex variety ${\mathcal{F}}_n^d\times {\mathcal{F}}_n^d$ and introduce the usual convolution product. We call the resulting algebra the \textit{affine $0$-Schur algebra} and denote it ${\mathbb S}_0^{\operatorname{aff}}(n,d)$. 

Our main construction provides a surjective map of the algebra $\Unzf$ onto ${\mathbb S}_0^{\operatorname{aff}}(n,d)$. The map is given explicitly by the images of the Drinfeld generators. Namely, we modify the Vasserot's construction from \cite{V98} to our case and consider explicit K-theory classes supported on the same closed orbits as in Vasserot's construction (also the same orbits as in the Beilinson-Lusztig-MacPherson construction from \cite{BLM90}).

One difference both from the Jensen-Su-Yang story and from the Vasserot algebra is that the relations we impose on Drinfeld generators are {\em quadratic}. However, the positive part of the Krob-Thibon $0$-quantum group can be mapped into our algebra, as the cubic plactic relations on the generators of the \textit{constant part} follow from our quadratic relations.

\subsection*{Structure of the paper.} The rest of the paper is structured as follows. Section \ref{section: Def and Prel} is devoted to reviewing the definitions and basic properties as well as formulating the main result. In Subsection \ref{subsection: U_n} we define the algebra $\Unz:=\Unzf$ using generators and relations. In Subsection \ref{subsection: K-theory} we recall the definitions and basic properties of the equivariant $K$-theory. In Subsection \ref{subsection: partial flags} we recall the definitions of the variety of partial flags $\mathcal{F}_n^d,$ the correspondence endomorphisms on $K^G(\mathcal{F}_n^d),$ and use the convolution product on $K^G(\mathcal{F}_n^d\times \mathcal{F}_n^d)$ to define the affine $0$-Schur algebra $\mathbb{S}_0^{\operatorname{aff}}(n,d).$ Finally, we formulate our main result: 

\smallskip

\noindent\textbf{Theorem A.} (Theorem \ref{theorem: main}). For every positive integer $d$ there exist a surjective morphism $$\phi_d:\Unz\to\mathbb{S}^{\operatorname{aff}}_0(n,d).$$

\smallskip

Section \ref{section: orbits} is devoted to the combinatorics of the orbit structure of the variety of pairs of partial flags $\mathcal{F}_n^d\times\mathcal{F}_n^d$ under the diagonal action of $G:=GL(V).$ In Subsection \ref{subsection: orbits matrices} we follow the work of Beilinson, Lusztig, MacPherson (\cite{BLM90}) and Vasserot (\cite{V98}) to establish that the orbits in $\mathcal{F}_n^d\times\mathcal{F}_n^d$ are enumerated by $n\times n$ matrices with non-negative integer entries, with the total sum of entries equal to $d.$ In Subsection \ref{subsection: bruhat order} we study the Bruhat order on these matrices given by the inclusion order on the closures of orbits. We establish explicit combinatorial conditions for the cover relations and for the order itself (see Theorem \ref{Theorem: partial order and cover relation}; an elementary combinatorial proof is provided in the Appendix). We further use these results to identify matrices corresponding to closed orbits in $\mathcal{F}_n^d\times\mathcal{F}_n^d$ and to give a geometric characterization of closed orbits (Corollary \ref{corollary: closed orbits}), and also to identify the matrices corresponding to open orbits (Corollary \ref{corollary: open orbits}). 

In Subsection \ref{subsection: convolution product for supports} we study an associative operation on the $n\times n$ matrices with non-negative integer entries arising from the generic multiplication of orbits in the variety of pairs of partial flags $\mathcal{F}_n^d\times\mathcal{F}_n^d.$ Computing the product of two arbitrary matrices is rather complicated and even mysterious, but the multiplication by an almost diagonal matrix (see Definition \ref{definition: almost diagonal}) can be described explicitly (see Theorem \ref{theorem: matrix mult by almost diag}).  The principal result of this subsection is the following:

\smallskip

\noindent\textbf{Theorem B} (Corollaries \ref{corollary:matrix factorization} and \ref{corollary:strictly monotone factorization}). Any matrix can factored in to a product of diagonal and almost diagonal matrices (with respect to the product introduced in Definition \ref{definition: almost diagonal}). 

\smallskip

For most of this subsection we follow \cite{BLM90} and \cite{V98}, but some of the approaches are new. In addition, we provide multiple examples.

In Section \ref{section: convolution algebra} we analyze the structure of the convolution algebra $\mathbb{S}^{\operatorname{aff}}_0(n,d).$ In Subsections \ref{subsection: K-theory of F} and \ref{subsection: K-theory of FxF} we review the equvariant $K$-theory of the variety of partial flags $\mathcal{F}_n^d$ and the orbits inside the variety of pairs of partial flags $\mathcal{F}_n^d\times\mathcal{F}_n^d.$ Those $K$-groups are isomorphic to partially symmetric Laurent polynomials in $d$ variables. In Subsection \ref{subsection: push-forward pull-back} we establish explicit formulas for the push-forward and pull-back morphisms along natural projections between orbits in the variety of pairs of partial flags (see Lemma \ref{Lemma: pull-back push-forward}). The formulas are in terms of the partially symmetric Laurent polynomials, introduced in Subsection \ref{subsection: K-theory of FxF}. We further derive explicit corollaries from these formulas for certain specific classes, used in the next Section (see Corollaries \ref{corollary: pi_*^-} and \ref{corollary: pi_*^+}, and Lemma \ref{lemma: push-forward for other steps}).

In Subsection \ref{subsection: convolution algebra} we use the results from Subsection \ref{subsection: convolution product for supports} and Section \ref{section: convolution algebra} to prove 

\smallskip

\noindent\textbf{Theorem C} (Theorem \ref{theorem: Snd generated by diag and almost diag}). The convolution algebra $\mathbb{S}^{\operatorname{aff}}_0(n,d)$ is generated by the classes supported on the orbits corresponding to diagonal and almost diagonal matrices.

\smallskip

In Subsection \ref{subsection: local generators} we introduce special classes $\mathfrak{e}_{\mu,k}(p),\mathfrak{f}_{\mu,k}(p),$ and $\mathfrak{h}_{\mu,n}(r)$ in $\mathbb{S}^{\operatorname{aff}}_0(n,d),$ supported on orbits corresponding to diagonal and almost diagonal matrices, and prove a series of key relations on these classes. We further show that 

\smallskip

\noindent\textbf{Theorem D} (Theorem \ref{theorem: Snd generated by E F H}).
The classes $\mathfrak{e}_{\mu,k}(p),\mathfrak{f}_{\mu,k}(p),$ and $\mathfrak{h}_{\mu,n}(r)$ generate the convolution algebra $\mathbb{S}^{\operatorname{aff}}_0(n,d).$

\smallskip

In Section \ref{section: the map phi_d} is devoted the map $\phi_d: \Unz\to \mathbb{S}^{\operatorname{aff}}_0(n,d)$ and the proof of the main theorem. In Subsection \ref{subsection: definition of phi_d} we introduce the map $\phi_d$ and use the relations from Subsection \ref{subsection: local generators} to prove that it is well defined. In Subsection \ref{subsection: surjectivity of phi_d} we use Theorem \ref{theorem: Snd generated by E F H} to prove surjectivity of the map $\phi_d.$ 

\section{Definitions and Preliminaries}\label{section: Def and Prel}

\subsection{The algebra $\Unz:=\Unzf$.}\label{subsection: U_n}

\begin{definition}
The algebra $\Unz$ is defined via generators and relations as follows. The set of generators consists of $E_i(p), F_i(q),$ and $H_n(r),$ for $0<i<n$ and $p,q,r\in\mathbb{Z}.$ The generators satisfy the following list of relations.

\bigskip

\textbf{$1.$ Relations on $E$ generators:}

\begin{equation}\tag{E1}\label{rel:2 term E}
E_i(p)E_i(q)=-E_i(q-1)E_i(p+1),
\end{equation}
for all $0<i<n$ and $p,q\in\mathbb{Z}.$

\begin{equation}\tag{E2}\label{rel:3 term E}
E_{i+1}(p)E_i(q)=E_i(q)E_{i+1}(p)-E_i(q+1)E_{i+1}(p-1),
\end{equation}
for all $0<i<n-1$ and $p,q\in\mathbb{Z}.$

\begin{equation}\tag{E3}\label{rel:comm E}
E_i(p)E_j(q)=E_j(q)E_i(p),
\end{equation}
whenever $|i-j|>1,\ 0<i,j<n$ and $p,q\in\mathbb{Z}.$

\bigskip

\textbf{$2.$ Relations on $F$ generators:}

\begin{equation}\tag{F1}\label{rel:2 term F}
F_i(p)F_i(q)=-F_i(q+1)F_i(p-1),
\end{equation}
for all $0<i<n$ and $p,q\in\mathbb{Z}.$

\begin{equation}\tag{F2}\label{rel:3 term F}
F_i(p)F_{i+1}(q)=F_{i+1}(q)F_i(p)-F_{i+1}(q-1)F_i(p+1),
\end{equation}
for all $0<i<n-1$ and $p,q\in\mathbb{Z}.$

\begin{equation}\tag{F3}\label{rel:comm F}
F_i(p)F_j(q)=F_j(q)F_i(p),
\end{equation}
whenever $|i-j|>1,\ 0<i,j<n$ and $p,q\in\mathbb{Z}.$

\bigskip

\textbf{$3.$ Relations between $E$ and $F$ generators:}

\begin{equation}\tag{EF1}\label{rel:comm E and F}
E_i(p)F_j(q)=F_j(q)E_i(p),
\end{equation}
For all $i\neq j,\ 0<i,j<n$ and $p,q\in\mathbb{Z}.$

\begin{equation}\tag{EF2}\label{rel: p+q}
E_i(p)F_i(q)-F_i(q)E_i(p)=E_i(p')F_i(q')-F_i(q')E_i(p'),
\end{equation}
whenever $p+q=p'+q'.$ For convenience, we introduce the notation:
\begin{equation}\label{def: H}
H_i(p+q):=E_i(p)F_i(q)-F_i(q)E_i(p).
\end{equation}

\begin{remark}
Note that $H_i(r)$ for $0<i<n,$ $r\in\mathbb{Z}$ are defined as the commutators of certain generators $E_i$ and $F_i$, as above, while $H_n(r),$ $r\in\mathbb{Z}$ are generators (see also Remark \ref{remark: U tilda} below).
\end{remark}

\bigskip

\textbf{$4.$ The $H$ elements commute:}

\begin{equation}\tag{H1}\label{rel:comm H}
H_i(p)H_j(q)=H_j(q)H_i(p),
\end{equation}
For all $0<i,j\le n$ and $p,q\in\mathbb{Z}.$

\bigskip

\textbf{$5.$ Relations between $H$ and $E$:}

\begin{equation}\tag{HE1}\label{rel: H_n vs E 3 term}
H_n(p)E_{n-1}(q)=E_{n-1}(q)H_n(p)-E_{n-1}(q+1)H_n(p-1),
\end{equation}
for all $p,q\in\mathbb{Z}.$

\begin{equation}\tag{HE2}\label{rel: H vs E 2 term}
H_i(p)E_i(q)=-E_i(q-1)H_i(p+1),
\end{equation}
for all $0<i<n$ and $p,q\in\mathbb{Z}.$

\begin{equation}\tag{HE3}\label{rel: H_n vs E commutes}
H_n(p)E_{i}(q)=E_{i}(q)H_n(p),
\end{equation}
for all $0<i<n-1$ and $p,q\in\mathbb{Z}.$

\bigskip

\textbf{$6.$ Relations between $H$ and $F$:}

\begin{equation}\tag{HF1}\label{rel: H_n vs F 3 term}
F_{n-1}(p)H_n(q)=H_n(q)F_{n-1}(p)-H_n(q-1)F_{n-1}(p+1),
\end{equation}
for all $p,q\in\mathbb{Z}.$

\begin{equation}\tag{HF2}\label{rel: H vs F 2 term}
H_i(p)F_i(q)=-F_i(q+1)H_i(p-1),
\end{equation}
for all $0<i<n$ and $p,q\in\mathbb{Z}.$

\begin{equation}\tag{HF3}\label{rel: H_n vs F commutes}
H_n(p)F_{i}(q)=F_{i}(q)H_n(p),
\end{equation}
for all $0<i<n-1$ and $p,q\in\mathbb{Z}.$
\end{definition}

\begin{remark}
Note that relation (\ref{rel:3 term E}) and the definition of $H_i(p)$ for $0<i<n$ (formula (\ref{def: H})) imply the analogs of relation (\ref{rel: H_n vs E 3 term}):
\begin{equation}\tag{HE1'}\label{rel: H vs E 3 term}
H_i(p)E_{i-1}(q)=E_{i-1}(q)H_i(p)-E_{i-1}(q+1)H_i(p-1),
\end{equation}
and
\begin{equation}\tag{HE1''}\label{rel: H vs E 3 term v2}
E_i(q)H_{i-1}(p)=H_{i-1}(p)E_i(q)-H_{i-1}(p+1)E_i(q-1),
\end{equation}
for all $0<i<n$ and $p,q\in\mathbb{Z}.$

Similarly, one also get the analogs of relation (\ref{rel: H_n vs F 3 term}): 
\begin{equation}\tag{HF1'}\label{rel: H vs F 3 term}
F_{i-1}(p)H_i(q)=H_i(q)F_{i-1}(p)-H_i(q-1)F_{i-1}(p+1),
\end{equation}
and
\begin{equation}\tag{HF1''}\label{rel: H vs F 3 term v2}
H_{i-1}(q)F_i(p)=F_i(p)H_{i-1}(q)-F_i(p-1)H_{i-1}(q+1),
\end{equation}
for all $0<i<n$ and $p,q\in\mathbb{Z}.$

Finally, it immediately follows from the definition of $H_i(p)$ for $0<i<n$ that $H_i(p)$ commutes with both $E_j(q)$ and $F_j(q)$ for $|i-j|>1:$

\begin{equation}\tag{HE3'}
H_i(p)E_j(q)=E_j(q)H_i(p),
\end{equation}
for all $|i-j|>1$ and $p,q\in\mathbb{Z}.$

\begin{equation}\tag{HF3'}
H_i(p)F_j(q)=F_j(q)H_i(p),
\end{equation}
for all $|i-j|>1$ and $p,q\in\mathbb{Z}.$
\end{remark}

\begin{remark}\label{remark: U tilda}
One can think about the extra generators $H_n(p)$ in the following manner. Let $\Unz'\subset \Unz$ be the subalgebra generated by $E_i(p)$ and $F_j(q),$ $0<i,j<n,$ $p,q\in\mathbb{Z}.$ One has a natural embedding $\Unz'\subset \mathfrak{U}_{n+1}',$ where $\mathfrak{U}_{n+1}'$ has the extra generators $E_n(p)$ and $F_n(q),$ $p,q\in\mathbb{Z}.$ Then one also gets
\begin{equation*}
\Unz'\subset \Unz\subset \mathfrak{U}_{n+1}',
\end{equation*}
where $\Unz$ is generated by $E_i(p), F_j(q),$ $0<i,j<n$ and $p,q\in\mathbb{Z},$ and $H_n(r),$ $r\in \mathbb{Z}.$ This approach allows to shorten the list of relations as relations \ref{rel: H_n vs E 3 term}, \ref{rel: H_n vs E commutes}, \ref{rel: H_n vs F 3 term}, and \ref{rel: H_n vs F commutes} become unnecessary.

\end{remark}

\subsection{Relations on the generators $E_i(0)$ and the Krob-Thibon $0$-quantum group}\label{subsection: plactic}

The positive part $\mathfrak{U}_0(\mathfrak{gl}_n)^+$ of the Krob-Thibon's $0$-quantum group is generated by elements $e_i,$ $0<i<n,$ subject to the \textit{plactic relations} (see \cite{KT99}):

\begin{equation}\label{rel:plactic}
e_ie_j=e_je_i\ \text{for}\ |i-j|>1,\ 
e_i^2e_{i+1}=e_ie_{i+1}e_i,\ 
e_{i+1}e_ie_{i+1}=e_ie_{i+1}^2.
\end{equation}

On the other hand, we have

\begin{lemma}
In the algebra $\mathfrak{U}_n$ relations \ref{rel:2 term E}, \ref{rel:3 term E} and \ref{rel:comm E} imply
\begin{equation*}
E_i(0)^2E_{i+1}(0)=E_i(0)E_{i+1}(0)E_i(0),
\end{equation*}

and

\begin{equation*}
E_{i+1}(0)E_i(0)E_{i+1}(0)=E_i(0)E_{i+1}(0)^2.
\end{equation*}
\end{lemma}

\begin{proof}
Indeed, using relation \ref{rel:2 term E} one gets:
\begin{equation*}
E_i(p)E_i(p+1)=-E_i(p)E_i(p+1),
\end{equation*}
therefore,
\begin{equation*}
E_i(p)E_i(p+1)=0.
\end{equation*}
Hence, using relation \ref{rel:3 term E} one gets
\begin{equation*}
E_i(0)E_{i+1}(0)E_i(0)=E_i(0)^2E_{i+1}(0)-E_i(0)E_i(1)E_{i+1}(-1)=E_i(0)^2E_{i+1}(0),
\end{equation*}
and
\begin{equation*}
E_{i+1}(0)E_i(0)E_{i+1}(0)=E_i(0)E_{i+1}(0)^2-E_i(1)E_{i+1}(-1)E_{i+1}(0)=E_i(0)E_{i+1}(0)^2.
\end{equation*}
\end{proof}

\begin{definition}
Let $\mathfrak{U}_n^+\subset \mathfrak{U}_n$ denote the subalgebra of $\mathfrak{U}_n$ generated by elements $E_i(p),$ $0<i<n,$ $p\in\mathbb{Z}.$
\end{definition}

\begin{corollary}
The map $i^+:\mathfrak{U}_0(\mathfrak{gl}_n)^+\to \mathfrak{U}_n^+$ defined by $i^+(e_i):=E_i(0)$ is well-defined.
\end{corollary}

We expect this map to be injective, but do not attempt to prove this here. The difficulty of proving injectivity can be underscored by the following considerations.

The relations \ref{rel:2 term E}, \ref{rel:3 term E}, and \ref{rel:comm E} imply that $\mathfrak{U}_n^+$ has a basis consisting of monomials

\begin{equation*}%\label{equation: monomial basis E}
E_{i_1}(p_{1,1})\cdot\ldots\cdot E_{i_1}(p_{1,k_1})\cdot E_{i_2}(p_{2,1})\cdot\ldots\cdot E_{i_2}(p_{2,k_2})\cdot\ldots\cdot E_{i_l}(p_{l,1}) \cdot\ldots\cdot E_{i_l}(p_{l,k_l}),
\end{equation*}  
where $i_1<i_2<\ldots<i_l$ and for each $1\le j\le l$ one has $p_{j,1}\ge p_{j,2}\ge\ldots\ge p_{j,k_j}.$ %In particular, one can say that $\mathfrak{U}_n^+$ has the size of the polynomial ring in variables $E_i(p).$
However, monomials

\begin{equation*}
    E_{i_1}(0)\cdot\ldots\cdot E_{i_l}(0),
\end{equation*}
where $i_1<i_2<\ldots<i_l$ do not form a basis of the subalgebra in $\mathfrak{U}_n^+$ generated by $E_i(0),$ $0<i<n.$

At the same time, one can construct a basis in $\mathfrak{U}_0(\mathfrak{gl}_n)^+$ as follows. For $i<j$ denote

\begin{equation*}
e_{ij}:=e_ie_{i+1}\ldots e_{j-1}.
\end{equation*}

Then there is a basis in $\mathfrak{U}_0(\mathfrak{gl}_n)^+$ consisting of the monomials in $e_{ij}$ of the form

\begin{equation*}
e_{i_1j_1}\ldots e_{i_lj_l},
\end{equation*}
where $i_1\ge i_2\ge \ldots \ge i_l,$ and $i_k=i_{k+1}$ implies $j_k\le j_{k+1}.$ Proving that the images of these elements are linearly independent in $\mathfrak{U}_n^+$ seems tricky.
 
\subsection{Equivariant $K$-theory.}\label{subsection: K-theory}

For a complex linear algebraic group $G$ and a $G$-variety $X,$ let $Coh^G(X)$ denote the category of $G$-equivariant coherent sheaves on $X$ and let $K^G(X)$ denote the complexified Grothendieck group of $Coh^G(X)$. Given $\mathcal{A}\in Coh^G(X)$ let $[\mathcal{A}]$ denote its class in $K^G(X).$ $K^G(X)$ then has a natural $R(G)$-module structure, where $R(G)=K^G(pt)$ is the complexified representation ring of $G.$ Here we summarize some basic properties of $K^G(X)$ that are going to be important for our constructions later in the paper. We refer to \cite[Sections 5.2 and 5.3]{CG97} for a systematic treatment of the subject. %Some of the properties can also be found in \cite{V98}.

\noindent\textbf{Induction. }%\cite[p. 249]{CG97}. 
For an algebraic subgroup $H\subset G$ and any $H$-variety $X$ one gets 

\begin{equation}\label{equation: induction}
K^H(X)\simeq K^G(G\times_H X),
\end{equation}
where $G\times_H X$ is the space of orbits of the free action of $H$ on $G\times X.$ In particular, for $X=pt$ one gets

\begin{equation}\label{equation: induction from pt}
K^G(G/H)\simeq K^H(pt)\simeq R(H).
\end{equation}

\noindent\textbf{Reduction. }%\cite[p. 249]{CG97}. 
Any algebraic group $G$ can be written as a semidirect product $R\ltimes U$ where $R$ is reductive and $U$ is the unipotent radical of $G.$ Then for any $G$-variety $X$ one gets

\begin{equation}\label{equation: reduction}
K^G(X)\simeq K^R(X).
\end{equation}

\noindent\textbf{Pushforward. }%\cite[p. 248]{CG97}. 
For any proper $G$-equivariant map $\phi:X\to Y$ between two $G$-varieties $X$ and $Y$ there is a derived direct image morphism 
\begin{equation}\label{equation: pushforward}
R\phi_*: K^G(X)\to K^G(Y),\ \ R\phi_*:=\sum_{i\ge 0} (-1)^i [R^i f_*].
\end{equation} 
Going forward we will use a simplified notation $\phi_*:=R\phi_*.$ Note that $\phi_*$ is a $\mathbb{C}$-linear map.

\noindent\textbf{Pullback. }%\cite[p. 244]{CG97}. 
For any flat $G$-equivariant map $\phi:X\to Y$ between two $G$-varieties $X$ and $Y$ there is an inverse image morphism 
\begin{equation}\label{equation: pullback}
\phi^*: K^G(Y)\to K^G(X).
\end{equation}
 In particular, one gets the pullback map for open embeddings and for smooth fibrations. 

\noindent\textbf{Restriction to closed smooth subvarieties. }%\cite[p. 244]{CG97}. 
Let $i: N\to M$ be a smooth closed embedding of a smooth $G$-subvariety $N$ into a smooth $G$-variety $M.$ Then there is a restriction morphism 
\begin{equation}\label{equation: restriction}
i^*: K^G(M)\to K^G(N).
\end{equation}

Furthermore, for any closed $G$-subvariety $A\subset M$ there is also a restriction with supports morphism
\begin{equation}\label{equation: restriction with supports}
    i^*: K^G(A)\to K^G(A\cap N).
\end{equation}
Note that this morphism depends on the ambient spaces $N$ and $M.$

\noindent\textbf{Tensor product. }%\cite[p. 246]{CG97}. 
If $M$ is a smooth $G$-variety then $K^G(M)$ has the structure of a commutative associative $R(G)$-algebra with the multiplication given by
\begin{equation}\label{equation: tensor}
f\otimes g:= \Delta^* (f\boxtimes g),
\end{equation}
where $\boxtimes: K^G(M)\otimes K^G(M)\to K^G(M\times M)$ is the exterior tensor product and $\Delta: M\to M\times M$ is the diagonal embedding. Note that this multiplication is only defined for a smooth $G$-variety $M.$ Otherwise the restriction $\Delta^*$ to the diagonal might not be well defined. 

\noindent\textbf{Multiplicativity}\cite[p. 270]{V98}. Let $X$ and $Y$ be smooth $G$-varieties and let $\phi:X\to Y$ be either flat or a closed embedding. Then the pullback map $\phi^*: K^G(Y)\to K^G(X)$ is multiplicative with respect to the tensor product, i.e. for any two classes $f,g\in K^G(Y)$ one has

\begin{equation}\label{equation: multiplicativity of pullback}
    \phi^*(f\otimes g)=\phi^*(f)\otimes \phi^*(g).
\end{equation}

\noindent\textbf{Tensor product with supports.} One can refine the above definition as follows. Let $M$ be a smooth $G$-variety and let $X,Y\subset M$ be two closed $G$-invariant subvarieties in $M.$ Then one gets the tensor product map
\begin{equation}\label{equation: tensor with support}
\otimes: K^G(X)\otimes K^G(Y)\to K^G(X\cap Y).
\end{equation} 
Note that this map depends on the smooth ambient variety $M.$

\noindent\textbf{Projection formula.} Let $M$ and $N$ be two smooth $G$-varieties and let $f:N\to M$ be a proper flat map. Then for any $\alpha\in K^G(M)$ and $\beta\in K^G(N)$ one has
\begin{equation}\label{equation: projection formula}
f_*(\beta)\otimes \alpha=f_*(\beta\otimes f^*(\alpha)).
\end{equation}

\noindent\textbf{Proper base change.} Consider a Cartesian square
\begin{center}
\begin{tikzpicture}

%\node (TT) at (3,3.5) {$L=\overline{U}_i/U_i$};
\node (X) at (0,0) {$X$};
\node (Y) at (3,2) {$Y$};
\node (Z) at (3,0) {$Z$};
\node (W) at (0,2) {$W=X\otimes_Z Y$};
\draw
(X) edge[->,>=angle 90] node[above] {$g$} (Z)
(Y) edge[->,>=angle 90] node[right] {$f$} (Z)
(W) edge[->,>=angle 90] node[left] {$\tilde{f}$} (X)
(W) edge[->,>=angle 90] node[above] {$\tilde{g}$} (Y)
;
\end{tikzpicture}
\end{center}
where $f$ is proper and $g$ is flat. It then follows that $\tilde{f}$ is also proper, and $\tilde{g}$ is also flat. One gets

\begin{equation}\label{equation: base change}
g^*\circ f_*=\tilde{f}_*\circ\tilde{g}^*.
\end{equation}

\noindent\textbf{Conormal bundle.} Let $M$ be a smooth $G$-variety and $N\subset M$ be a smooth closed $G$-subvariety. Let $i:N\to M$ be the embedding. Condsider the conormal vector bundle $T^*_N M.$ Let
\begin{equation*}
\Lambda(T^*_N M):=\sum\limits_{k=0}^\infty (-1)^k [\Lambda^i T^*_N M]\in K^G(N)
\end{equation*}
be the $K^G$ class of the exterior algebra of the conormal bundle $T^*_N M.$ Then for any $\alpha\in K^G(N)$ one gets
\begin{equation}\label{equation: conormal bundle}
i^* (i_* \alpha) = \alpha \otimes \Lambda(T^*_N M).
\end{equation}

\noindent\textbf{Localization and Lefschetz formula}\cite[p. 270]{V98}. Let $T\simeq (\mathbb{C}^*)^n$ be a complex torus and $t\in T.$ In this case the complexified representation ring $R(T)$ can be identified with the ring of regular functions on $T$ by considering the characters of the representations. Consider the multiplicative set $S_t\subset R(T)$ consisting of elements of $R(T)$ such that the trace does not vanish at $t.$ Let $R(T)_t$ denote the localization of $R(T)$ with respect to $S_t.$

Let $M$ be a smooth $T$-variety, $M_t\subset M$ be the subvariety of $t$-fixed points, and $i:M_t\to M$ be the (smooth, closed) embedding. Consider the localized equivariant $K$ groups
\begin{equation*}
K^T(M)_t:=R(T)_t\otimes_{R(T)} K^T(M), 
\end{equation*}
and
\begin{equation*} 
K^T(M_t)_t:=R(T)_t\otimes_{R(T)} K^T(M_t).
\end{equation*}
Then $\Lambda(T^*_{M_t} M)$ is invertible in $K^T(M_t)_t$ and $i_*$ induces an isomorphism between $K^T(M_t)_t$ and $K^T(M)_t.$  

Let $p: M\to pt$ and $q:=p\circ i: M_t\to pt$ be projections to a point, and $\alpha\in K^G(M)$. Then one also gets the \textit{Lefschetz formula}:
\begin{equation}\label{equation: Lefschetz formula}
p_*(\alpha)=q_*(i^*(\alpha)\otimes \Lambda(T^*_{M_t} M)^{-1}).
\end{equation}

\noindent\textbf{Long exact sequence.} Let $X$ be a $G$-variety and let $i:Y\to X$ be a closed embedding of a $G$-subvariety. Let also $U=Y\setminus X$ be the complement and $j:U\to X$ be the complementary open embedding. Then one gets a long exact sequence:

\begin{equation}\label{equation: long exact sequence}
\begin{tikzcd}
\ldots \arrow[r] &K^G(Y) \arrow[r,"i_*"] & K^G(X) \arrow[r, "j^*"] & K^G(U)\arrow[r] & 0.
\end{tikzcd}    
\end{equation}

On the left the sequence extends with the use of higher $K$-theory groups (see \cite[p. 248]{CG97} for details). However, for the purposes of this paper we only need the part of the sequence displayed above.

\subsection{The variety of partial flags and the convolution algebra.}\label{subsection: partial flags}

%Convolution algebras associated with the varieties of pairs of partial flags were used to provide geometric frameworks for various versions of quantum groups. In \cite{BLM90} Beilinson, Lusztig, and MacPherson used point counting over finite fields to provide geometric setting for the quantum group of $\mathfrak{gl}_n.$ In \cite{V98} Vasserot used equivariant $K$-theory of the Steinberg variety to provide a geometric setting for the quantum group of the affine Lie algebra $\hat{\mathfrak{gl}}_n.$ In this paper we are using the equivariant $K$-theory of the varieties of pairs of partial flags to study algebra $\Unz$ defined in the previous section. Our construction can be thought of as a ``$q=0$'' version of the Vasserot's construction, although one cannot directly plug $q=0$ into Vasserot's formulas. Below we review the definitions of the varieties of partial flags and the convolution product.

Let $V\simeq\mathbb{C}^d,$ and $G:=GL(V).$ We say that $\mu$ is a \textit{composition} of $d$ of length $n$ if $\mu=(\mu_1,\ldots,\mu_n)\in\mathbb{Z}_{\ge 0}^n$ and $\sum \mu_i=d.$ We write $\mu \vDash d$ and $l(\mu)=n$ in that case. Note that the last part $\mu_n,$ as well as any other part, is allowed to be zero. We will also use the partial sums $d_k:=\sum\limits_{j=1}^k \mu_j$.

\begin{definition}
The variety of $\mu$-partial flags $\Fmu$ is defined by the following  
$$
\Fmu:=\{U_1\subset\ldots\subset U_n\subset V: \dim U_k=d_k\}.
$$ 
We will also use the notation $\Fnd:=\bigsqcup \Fmu.$
\end{definition}

%\subsection{Correspondences, variety of pairs of flags, and the convolution product}

Consider the variety of pairs of partial flags $\Fnd\times \Fnd.$ The group $G=GL(V)$ acts on it diagonally. Let  $\alpha\in K^G(\Fnd\times \Fnd)$ be an equivariant K-theory class. Consider the diagram

\smallskip

\begin{center}
\begin{tikzpicture}

%\node (TT) at (3,3.5) {$L=\overline{U}_i/U_i$};
\node (T) at (3,2) {$\Fnd\times \Fnd$};
\node (L1) at (0,0) {$\Fnd$};
\node (L2) at (6,0) {$\Fnd$};
\draw
%(TT) edge[->,>=angle 90] node[left] {} (T)
(T) edge[->,>=angle 90] node[left] {$\pi_1$} (L1)
(T) edge[->,>=angle 90] node[right] {$\pi_2$} (L2)
;
\end{tikzpicture}
\end{center}

\begin{definition}
The \textit{correspondence morphism} $\phi_\alpha:K^G(\Fnd)\to K^G(\Fnd)$ is defined as follows:

\begin{equation*}
\phi_\alpha(f):={\pi_2}_*(\pi_1^*f\otimes \alpha).
\end{equation*}
\end{definition}

Let $\beta\in K^G(\Fnd\times \Fnd)$ be another equivariant $K$-theory class. Consider the composition $\phi_\beta\circ\phi_\alpha:K^G(\Fnd)\to K^G(\Fnd).$ We get the following diagram:

\smallskip

\begin{center}
\begin{tikzpicture}

\node (TT) at (6,3.5) {$\Fnd\times \Fnd\times \Fnd$};
\node (T1) at (3,2) {$\Fnd\times \Fnd$};
\node (T2) at (9,2) {$\Fnd\times \Fnd$};
\node (L1) at (0,0) {$\Fnd$};
\node (L2) at (6,0) {$\Fnd$};
\node (L3) at (12,0) {$\Fnd$};
\draw
(TT) edge[->,>=angle 90] node[left] {$\pi_{12}$} (T1)
(TT) edge[->,>=angle 90] node[right] {$\pi_{23}$} (T2)
(T1) edge[->,>=angle 90] node[left] {$\pi_1$} (L1)
(T1) edge[->,>=angle 90] node[right] {$\pi_2$} (L2)
(T2) edge[->,>=angle 90] node[left] {$\pi_1$} (L2)
(T2) edge[->,>=angle 90] node[right] {$\pi_2$} (L3)
;
\end{tikzpicture}
\end{center}

Here $\pi_{12}:\Fnd\times \Fnd\times \Fnd\to \Fnd\times \Fnd$ and $\pi_{23}:\Fnd\times \Fnd\times \Fnd\to \Fnd\times \Fnd$ are the projections onto the first two factors and the second two factors correspondingly. Note that the variety of triples of partial flags is the fibered product of the two copies of the variety of pairs of partial flags over the projections $\pi_2$ and $\pi_1.$ Therefore, using the base change and the projection formula we get:

\begin{equation*}
\begin{aligned}
\phi_\beta\circ\phi_\alpha(f)&={\pi_2}_*\left(\beta\otimes[\pi_1^*\circ{\pi_2}_*(\alpha\otimes \pi_1^*f)]\right)\\
&={\pi_2}_*\left(\beta\otimes[{\pi_{23}}_*\circ\pi_{12}^*(\alpha\otimes \pi_1^*f)]\right)\\
&={p_3}_*(\pi_{12}^*(\alpha)\otimes \pi_{23}^*(\beta)\otimes p_1^*f),
\end{aligned}
\end{equation*}
where $p_1=\pi_1\circ\pi_{12}:\Fnd\times \Fnd\times \Fnd\to \Fnd$ and $p_3=\pi_2\circ\pi_{23}:\Fnd\times \Fnd\times \Fnd\to \Fnd$ are the projections onto the first and the third factors. Consider the diagram

\smallskip

\begin{center}
\begin{tikzpicture}

\node (TT) at (6,3.5) {$\Fnd\times \Fnd\times \Fnd$};
\node (T) at (6,1.5) {$\Fnd\times \Fnd$};
\node (L1) at (2.5,0) {$\Fnd$};
\node (L2) at (9.5,0) {$\Fnd$};
\draw
(TT) edge[->,>=angle 90] node[left] {$\pi_{13}$} (T)
(TT) edge[->,>=angle 90] node[left] {$p_1$} (L1)
(TT) edge[->,>=angle 90] node[right] {$p_3$} (L2)
(T) edge[->,>=angle 90] node[left] {$\pi_1$} (L1)
(T) edge[->,>=angle 90] node[right] {$\pi_2$} (L2)
;
\end{tikzpicture}
\end{center}
where $\pi_{13}$ is the projection on the first and the third factors. The triangles in the diagram are commutative. Therefore, using the projection formula we get

\begin{equation*}
\begin{aligned}
\phi_\beta\circ\phi_\alpha(f)&={p_3}_*(\pi_{12}^*(\alpha)\otimes \pi_{23}^*(\beta)\otimes p_1^*f)\\
&={\pi_2}_*{\pi_{13}}_*(\pi_{12}^*(\alpha)\otimes \pi_{23}^*(\beta)\otimes \pi_{13}^*\pi_1^*f)\\
&={\pi_2}_*[{\pi_{13}}_*(\pi_{12}^*\alpha\otimes \pi_{23}^*\beta)\otimes \pi_1^*f].\\
\end{aligned}
\end{equation*}

This motivates the definition of the convolution product on $K^G(\Fnd\times \Fnd):$

\begin{definition}
The \textit{convolution product} of the classes $\alpha, \beta\in K^G(\Fnd\times \Fnd)$ is defined by
\begin{equation*}
\alpha\star\beta:={\pi_{13}}_*(\pi_{12}^*\alpha\otimes \pi_{23}^*\beta).
\end{equation*}
\end{definition}

We conclude that $\Snd:=K^G(\Fnd\times \Fnd)$ is an algebra under the convolution multiplication, and it acts on $K^G(\Fnd)$ via the correspondence morphisms. The main result of this paper is the following theorem:

\begin{theorem}\label{theorem: main}
For every positive integer $d$ there exists a surjective morphism 
\begin{equation}
\phi_d:\Unz\to \Snd.
\end{equation}
\end{theorem}

\subsection{Transversality Lemma} The following lemma is important for computations:

\begin{lemma}\label{lemma: transversality}
Let $M,N\subset \Fnd\times \Fnd$ be two smooth locally closed $G$-subvarieties. Then their preimages $\pi_{12}^{-1}(M)$ and $\pi_{23}^{-1}(N)$ intersect transversely in $\mathcal{F}_n^d\times\mathcal{F}_n^d\times\mathcal{F}_n^d.$ In other words, for any point 
\begin{equation*}
    x\in\pi_{12}^{-1}(M)\cap\pi_{23}^{-1}(N)\subset\mathcal{F}_n^d\times\mathcal{F}_n^d\times\mathcal{F}_n^d,
\end{equation*}
one has
\begin{equation*}
    T_x(\pi_{12}^{-1}(M))+T_x(\pi_{23}^{-1}(N))=T_x(\mathcal{F}_n^d\times\mathcal{F}_n^d\times\mathcal{F}_n^d).
\end{equation*}
\end{lemma}

\begin{proof}
Without loss of generality one can assume that $M$ and $N$ are connected. Let compositions $\mu_1,$ $\mu_2,$ and $\mu_3$ be such that $x=(x_1,x_2,x_3)\in\mathcal{F}_{\mu_1}\times\mathcal{F}_{\mu_2}\times\mathcal{F}_{\mu_3}.$ Note that $\mathcal{F}_{\mu_1}\times\mathcal{F}_{\mu_2}\times\mathcal{F}_{\mu_3}$ is a connected component of $\mathcal{F}_n^d\times\mathcal{F}_n^d\times\mathcal{F}_n^d.$ We also conclude that $M\subset\mathcal{F}_{\mu_1}\times\mathcal{F}_{\mu_2}$ and $N\subset\mathcal{F}_{\mu_2}\times\mathcal{F}_{\mu_3}.$

Since $G$ acts transitively on  $\mathcal{F}_{\mu_2},$ the restriction $\psi:=\pi_2|_M$ of the projection $\pi_2:\mathcal{F}_{\mu_1}\times\mathcal{F}_{\mu_2}\to \mathcal{F}_{\mu_2}$ is surjective. Furthermore, it has no critical points. Indeed, otherwise, due to the transitivity of the $G$-action, every point in $\mathcal{F}_{\mu_2}$ would be a critical value of $\psi:M\to\mathcal{F}_{\mu_2}$, which contradicts Sard's Lemma (\cite{S42}). In particular, in the tangent space $T_{(x_1,x_2)}(\mathcal{F}_{\mu_1}\times\mathcal{F}_{\mu_2})$
%$T_{(x_1,x_2)}(\mathcal{F}_{\mu_1}\times\mathcal{F}_{\mu_2})=T_{x_1}(\mathcal{F}_{\mu_1})\oplus T_{x_2}(\mathcal{F}_{\mu_2})$ 
one has 
\begin{equation*}
    T_{x_2}(\mathcal{F}_{\mu_2})\simeq T_{(x_1,x_2)}(M)\slash \ker d\psi|_{(x_1,x_2)} =T_{(x_1,x_2)}(M)\slash \left((T_{(x_1,x_2)}(\mathcal{F}_{\mu_1}\times \{x_2\})\cap T_{(x_1,x_2)}(M)\right),
\end{equation*}
which implies that
\begin{equation*}
T_{(x_1,x_2)}(M)+T_{(x_1,x_2)}(\mathcal{F}_{\mu_1}\times \{x_2\})=T_{(x_1,x_2)}(\mathcal{F}_{\mu_1}\times\mathcal{F}_{\mu_2}).%=T_y(\mathcal{F}_{\mu_1})\oplus T_y(\mathcal{F}_{\mu_2})
\end{equation*}

Using the results above and making computations inside $T_x(\mathcal{F}_{\mu_1}\times\mathcal{F}_{\mu_2}\times\mathcal{F}_{\mu_3}),$
%$$T_x(\mathcal{F}_{\mu_1}\times\mathcal{F}_{\mu_2}\times\mathcal{F}_{\mu_3})=T_{x_1}(\mathcal{F}_{\mu_1})\oplus T_{x_2}(\mathcal{F}_{\mu_2})\oplus T_{x_3}(\mathcal{F}_{\mu_3}),$$ 
one gets
\begin{align*}
    T_x&\bigl(\pi_{12}^{-1}(M)\bigr)+T_x\bigl(\pi_{23}^{-1}(N)\bigr)=T_x\bigl(M\times\mathcal{F}_{\mu_3}\bigr)+T_x\bigl(\mathcal{F}_{\mu_1}\times N\bigr)\\
    =&T_x\Bigl(M\times\{x_3\}\Bigr)+T_x\Bigl(\mathcal{F}_{\mu_1}\times\{x_2\}\times\{x_3\}\Bigr)+\\
    &\qquad \qquad \qquad \qquad \qquad \qquad +T_x\Bigl(\{x_1\}\times\{x_2\}\times\mathcal{F}_{\mu_3}\Bigr)+T_x\Bigl(\{x_1\}\times N\Bigr)\\
%    =&\Bigl(T_{(x_1,x_2)}(M),T_{x_3}(\mathcal{F}_{\mu_3})\Bigr)+\Bigl(T_{x_1}(\mathcal{F}_{\mu_1}),0,0\Bigr)+\Bigl(0,T_{(x_2,x_3)}(N)\Bigr)\\
%    =&\Bigl(T_{(x_1,x_2)}(M)+\bigl(T_{x_1}(\mathcal{F}_{\mu_1}),0\bigr),T_{x_3}(\mathcal{F}_{\mu_3})\Bigr)+\Bigl(0,T_{(x_2,x_3)}(N)\Bigr)\\
%    =&(T_{(x_1,x_2)}(M)+T_{(x_1,x_2)}(\mathcal{F}_{\mu_1},x_2),T_{x_3}(\mathcal{F}_{\mu_3}))+(T_{x_1}(\mathcal{F}_{\mu_1}),T_{x_2,x_3}(N))\\
    =&T_x\Bigl(\mathcal{F}_{\mu_1}\times\mathcal{F}_{\mu_2}\times \{x_3\}\Bigr)+T_x\Bigl(\{x_1\}\times\{x_2\}\times\mathcal{F}_{\mu_3}\Bigr)=T_x\Bigl(\mathcal{F}_n^d\times\mathcal{F}_n^d\times\mathcal{F}_n^d\Bigr).
\end{align*}
\end{proof}

This remark will be important for computations because of the following standard fact. Suppose that $X$ is a smooth $G$-variety, $i_A:A\hookrightarrow X$ and $i_B:B\hookrightarrow X$ be closed embeddings of irreducible $G$-subvarieties, and suppose that $A$ and $B$ only intersect at smooth points, and that the intersection is transversal, in particular $i:A\cap B\hookrightarrow X$ is a closed smooth embedding. Then for any classes $\alpha\in K^G(A)$ and $\beta\in K^G(B)$ one has

\begin{equation}\label{equation: tensor with transversality}
    i_{A*}\alpha\otimes i_{B*}\beta = i_* (\alpha\otimes\beta) = i_* \bigl(\alpha|_{A\cap B}\otimes \beta|_{A\cap B}\bigr).
\end{equation}

Note that the tensor product on the right can be computed in $K^G(A\cap B).$ In particular, this implies the following corollary for the convolution products. 

\begin{corollary}\label{corollary: convolution with restriction}
    Let $A,B\subset \Fnd\times\Fnd$ be closed irreducible $G$-invariant subvarieties, and $U\subset \Fnd\times\Fnd$ be an open $G$-invariant subset. Suppose that $A\times\Fnd\cap\pi_{13}^{-1}(U)$ and $\Fnd\times B\cap\pi_{13}^{-1}(U)$ intersect only in smooth points (according to Lemma \ref{lemma: transversality} the intersection is transversal). Let $\tilde{\pi_{12}},$ $\tilde{\pi_{23}},$ and $\tilde{\pi_{13}}$ be the restrictions of the projections: 
    $$\tilde{\pi}_{ij}:=\pi_{ij}|_{A\times\Fnd\cap\Fnd\times B\cap\pi^{-1}_{13}(U)}.$$
    Then for any $\alpha\in K^G(A)$ and $\beta\in K^G(B)$ one gets
\begin{equation}
    (\alpha\star\beta)|_U=\tilde{\pi}_{13*}(\tilde{\pi}_{12}^*\alpha\otimes\tilde{\pi}_{23}^*\beta).
\end{equation}
    In particular, if $A$ and $B$ are smooth, e.g. closed orbits, one can take $U=\Fnd\times\Fnd.$
\end{corollary}
%given $K^G(\Fnd\times\Fnd)$ classes $\alpha$ and $\beta$ supported on closed orbits $\mathbb{O}_1$ and $\mathbb{O}_2$ respectively, one can compute the convolution product using the restrictions of the projections $\tilde{\pi}_{ij}:=\pi_{ij}|_{(\mathbb{O}_1\times\Fnd)\cap (\Fnd\times\mathbb{O}_2)}:$

%\begin{equation*}
%\alpha\star\beta={\pi_{13}}_*(\pi_{12}^*\alpha\otimes \pi_{23}^*\beta)=\tilde{\pi}_{13*}(\tilde{\pi}_{12}^*\alpha\otimes \tilde{\pi}_{23}^*\beta).
%\end{equation*}

In order to describe the morphisms $\phi_d$ we need first to explicitly describe the orbit structure of the variety of pairs of partial flags $\Fnd\times \Fnd$ and the convolution product on $\Snd).$ 

\section{Orbits Structure of the Varieties of Pairs and Triples of Partial Flags.}\label{section: orbits}

\subsection{Orbits and Matrices}\label{subsection: orbits matrices}

The orbit structure of the variety of pairs of partial flags was studied in \cite{V98,BLM90} and many other sources. In this section we largely follow \cite{V98,BLM90}, although some of the results did not previously appear in literature, up to the authors' knowledge. %For the sake of completeness of the exposition we provide most of the proofs.  

Note that the orbits of the $G$ action on the variety of partial flags $\Fnd$ coincide with the connected components $\Fmu,$ i.e. all orbits are open and closed, and enumerated by compositions $\mu$ of $d$ of length $n.$ The orbit structure of the diagonal $G$ action on the variety of pairs of partial flags $\Fnd\times \Fnd$ is more subtle. The key observation here is the following standard fact from elementary linear algebra:

\begin{lemma}
Let $\textbf{U}:=\{\{0\}=U_0\subset U_1\subset\ldots\subset U_n=V\}\in \Fnd$ and $\textbf{W}:=\{\{0\}=W_0\subset W_1\subset\ldots\subset W_m=V\}\in \mathcal{F}_m^d$ be two partial flags. Then there exists a basis $\{v_1,\ldots,v_d\}$ such that every element of both flags is a span of basic vectors. 
\end{lemma}

%\begin{proof}
%One constructs such a basis recursively, starting from a basis of $U_1\cap W_1$ and then gradually extending it to other intersections $U_i\cap W_j.$

%Suppose that $\mathcal{A}$ is a basis of $U_{i-1}\cap W_{j-1},$ $\mathcal{B}$ is such that $\mathcal{A}\sqcup\mathcal{B}$ is a basis of $U_i\cap W_{j-1},$ and $\mathcal{C}$ is such that $\mathcal{A}\sqcup\mathcal{C}$ is a basis of $U_{i-1}\cap W_j,$ then $\mathcal{A}\sqcup\mathcal{B}\sqcup\mathcal{C}$ is linearly independent. Indeed, otherwise one would have
%\begin{equation*}
%A+B+C=0,
%\end{equation*}
%where $A,$ $B,$ and $C$ are linear combinations of elements of $\mathcal{A},$ $\mathcal{B},$ and $\mathcal{C}$ respectively, and at least one of the vectors $A+B$ and $A+C$ is not zero. Without loss of generality, assume that $A+B\neq 0.$ But then 
%\begin{equation*}
%0\neq-C=A+B\in (U_i\cap W_{j-1})\cap(U_{i-1}\cap W_j)=U_{i-1}\cap W_{j-1},
%\end{equation*}
%which implies that $B=0,$ because $\mathcal{A}$ is a basis of $U_{i-1}\cap W_{j-1}$ and $\mathcal{A}\sqcup\mathcal{B}$ is linearly independent, which contradicts the linear independence of $\mathcal{A}\sqcup\mathcal{C}.$ One can then proceed to extend $\mathcal{A}\sqcup\mathcal{B}\sqcup\mathcal{C}$ to a basis of $U_i\cap W_j.$

%Using the above procedure, one can keep extending the bases of intersections of the elements of the flags $\textbf{U}$ and $\textbf{W}$ until a desired basis of the total space is constructed. 
%\end{proof}

\begin{definition}[Section 1.1 in \cite{BLM90}, Section 1.2 in \cite{V98}]
Given a pair of flags $(\textbf{U},\textbf{W})\in \Fnd\times \Fnd$ define the corresponding matrix $M:=\{m_{ij}\}$ in such a way that it satisfies  

\begin{equation*}
\dim U_a\cap W_b=\sum\limits_{(i,j)\le(a,b)}m_{ij}.
\end{equation*}

where by $(i,j)\le(a,b)$ we will always mean that $i\le a$ and $j\le b.$
\end{definition}

\begin{corollary}[\cite{V98,BLM90}]
Two pairs of flags belong to the same orbit of the diagonal action of $GL(V)$ on $\Fnd\times \Fnd$ if and only if the corresponding matrices are the same.
\end{corollary}

Let $\mathcal{M}(n,d)$ be the set of $n\times n$ matrices with non-negative integer entries and the total sum of entries equal to $d.$ For a matrix $M\in \mathcal{M}(n,d)$ let $\mathbb{O}_M\subset \Fnd\times \Fnd$ denote the corresponding orbit. 

\begin{example}
Let $M$ be a diagonal matrix and $(\textbf{U},\textbf{W})\in \mathbb{O}_M$. Then for every $i,j$ we get
\begin{equation*}
U_i\cap W_j=U_{\min(i,j)}=W_{\min(i,j)},
\end{equation*}
i.e. $\textbf{U}=\textbf{W}.$
\end{example}

Given a matrix $M=\{m_{ij}\}\in\mathcal{M}(n,d)$ and a pair of flags $(\textbf{U},\textbf{W})\in \mathbb{O}_M,$ one gets $\textbf{U}\in \mathcal{F}_{\textbf{R}^M}$ and $\textbf{W}\in \mathcal{F}_{\textbf{C}^M},$ where 
\begin{equation*}
\textbf{R}^M:=(R_1^M,\ldots, R_n^M),\ R_k^M:=\sum\limits_{i=1}^n m_{ki}, 1\le k\le n,
\end{equation*}
is the row-sum composition, and
\begin{equation*}
\textbf{C}^M:=(C_1^M,\ldots, C_n^M),\ C_k^M:=\sum\limits_{i=1}^n m_{ik}, 1\le k\le n,
\end{equation*}
is the column-sum composition. Hence one gets $\mathcal{M}(n,d)=\bigsqcup\limits_{\mu,\nu}\mathcal{M}_{\mu,\nu},$ where
\begin{equation*}
\mathcal{M}_{\mu,\nu}:=\{M\in\mathcal{M}(n,d)| \textbf{R}^M=\mu,\textbf{C}^M=\nu\},
\end{equation*}
and $\mathbb{O}_M$ and $\mathbb{O}_N$ belong to the same connected component of $\Fnd\times \Fnd$ if and only if $\textbf{R}^M=\textbf{R}^N$ and $\textbf{C}^M=\textbf{C}^N.$

\subsection{Bruhat order}\label{subsection: bruhat order}

The inclusion (partial) order on the closures of orbits defines a partial order on $\mathcal{M}(n,d).$ This order is usually called the \textit{Bruhat order,} and it was studied in various sources, see \cite{BLM90,MWZ99,S07,S18,V98} among others. Some of the material in this Section follow these sources, although an explicit description of the cover relations and the characterization of the minimal and maximal elements  in terms of the matrices, as well as an explicit characterization of the open orbits appear to be new.

\begin{definition}
We say that $N\prec M$ if $\mathbb{O}_N\subset\partial \mathbb{O}_M.$ Let also $N\lessdot M$ denote the cover relation, i.e. $N\lessdot N$ if and only if $N\prec M,$ and there is no matrix $K$ such that $N\prec K\prec M.$
\end{definition}

Let $M\in\mathcal{M}(n,d)$ be a matrix and $(\textbf{U},\textbf{W})\in \mathbb{O}_M$ be a pair of flags in the corresponding orbit. Let also $\mathcal{B}=\bigsqcup\limits_{1\le i,j\le n} \mathcal{B}_{ij}$ be a basis of $V$ such that for any $1\le k,l\le n$ one has
\begin{equation*}
U_k\cap W_l=\spn(\bigsqcup\limits_{i=1}^k\bigsqcup\limits_{j=1}^l\mathcal{B}_{ij}).
\end{equation*} 
Then the dimension of the stabilizer of $(\textbf{U},\textbf{W})$ is given by
\begin{equation*}
\dim\text{Stab}(\textbf{U},\textbf{W})=\sharp\{(u,v)|u\in\mathcal{B}_{ij},\ v\in\mathcal{B}_{qp},\ i\le q,\ j\le p\}=\sum\limits_{1\le i,j\le n} sm_{ij}m_{ij},
\end{equation*}
where for all $1\le k,l\le n$
\begin{equation*}
sm_{kl}:=\sum\limits_{i=1}^k\sum\limits_{j=1}^l m_{ij}.
\end{equation*}
Therefore, one gets
\begin{align*}
\dim (\mathbb{O}_M)&=d^2-\dim\text{Stab}(\mathbb{O}_M)=\sharp\{(u,v)|u\in\mathcal{B}_{ij},\ v\in\mathcal{B}_{qp},\ i>q\ \text{or}\ j>p\}\\
&=\sum\limits_{1\le i,j\le n} (d-sm_{ij})m_{ij}.
\end{align*}

\begin{definition}
We will use the notation $E_{ij}$ for the matrix with all entries equal to zero, except for the entry on $i$th row and $j$th column, which is equal to $1.$ %We will also use the notation
%\begin{equation*}
%S_{ij}^{op}:=E_{io}+E_{jp}-E_{ip}-E_{jo}.
%\end{equation*}
\end{definition}

The following result is standard although we couldn't locate this particular formulation in the literature (see the proof below for references).

\begin{theorem}\label{Theorem: partial order and cover relation}
Let $N=\{n_{ij}\}$ and $M=\{m_{ij}\}$ be two non-negative integer matrices. Then $N\preceq M$ if and only if the following two conditions are satisfied:
\begin{enumerate}
\item One has $\textbf{R}^M=\textbf{R}^N$ and $\textbf{C}^M=\textbf{C}^N.$ 
%For every $k$ one has $\sum\limits_i m_{ki}=\sum\limits_i n_{ki}$ and $\sum\limits_i m_{ik}=\sum\limits_i n_{ik},$
\item For every $k,l$ one has $\sum\limits_{i=1}^k\sum\limits_{j=1}^l n_{ij}\ge \sum\limits_{i=1}^k\sum\limits_{j=1}^l m_{ij}.$ 
\end{enumerate}
Moreover, one can also characterize the cover relations: one has $N\lessdot M$ if and only if there exist four integers $1\le a_1<a_2\le n$ and $1\le b_1<b_2\le n$ such that the following two conditions are satisfied:
\begin{enumerate}
\item[(i)] $M-N=E_{a_1b_2}+E_{a_2b_1}-E_{a_1b_1}-E_{a_2b_2}.$ %where $E_{ij}$ is the matrix with all entries equal to zero, except for the entry on the $i$th row and the $j$th column, which is equal to $1.$
\item[(ii)] $m_{ij}=n_{ij}=0\ $ for $\ a_1\le i\le a_2$ and $b_1\le j\le b_2$ except for 
\begin{equation*}
(i,j)\in \{(a_1,b_1),(a_1,b_2),(a_2,b_1),(a_2,b_2)\}.
\end{equation*}
\end{enumerate}
\end{theorem}

See Figure \ref{fig: cover relation} for an illustration of a cover relation.

\begin{figure}
    \centering
    $\begin{pmatrix}
           &\vdots& & & &\vdots&       \\
    \ldots & a+1  &0&0&0& b    &\ldots \\
           & 0  &0&0&0& 0      &   \\
  %         & 0  &0&0&0& 0      &   \\
           & 0  &0&0&0& 0      &   \\
    \ldots & c  &0&0&0& d+1    &\ldots \\
           &\vdots& & & &\vdots&       
    \end{pmatrix}\lessdot
    \begin{pmatrix}
           &\vdots& & & &\vdots&       \\
    \ldots & a  &0&0&0& b+1    &\ldots \\
           & 0  &0&0&0& 0      &   \\
 %          & 0  &0&0&0& 0      &   \\
           & 0  &0&0&0& 0      &   \\
    \ldots & c+1&0&0&0& d      &\ldots \\
           &\vdots& & & &\vdots&       
    \end{pmatrix}$
    \caption{Example of a cover relation. Here $a,b,c$ and $d$ are non-negative integers. The rectangular area of the matrix can be of any size and located anywhere in the matrix.}
    \label{fig: cover relation}
\end{figure}

\begin{proof}
Suppose that $N\preceq M.$ Then $\textbf{R}^M=\textbf{R}^N$ and $\textbf{C}^M=\textbf{C}^N,$ because $\mathbb{O}_M$ and $\mathbb{O}_N$ belong to the same connected component. Also, for any $1\le k,l\le n$ and any $(\textbf{U},\textbf{W})\in \mathbb{O}_M$ one has $\dim(U_k\cap W_l)=\sum\limits_{i=1}^k\sum\limits_{j=1}^l m_{ij}.$ As a family of pairs of flags approaches the boundary of the orbit $\mathbb{O}_M,$ the dimension of the intersection can only increase. Therefore, one gets 
\begin{equation*}
\sum\limits_{i=1}^k\sum\limits_{j=1}^l n_{ij}\ge \sum\limits_{i=1}^k\sum\limits_{j=1}^l m_{ij}.
\end{equation*}
The fact that conditions (1) and (2) imply that $N\prec M$ follows from work of Bongartz on representations of quivers (see \cite{B95,B96}). See also \cite{MWZ99} and \cite{S07} for a more elementary approach.

It is easy to check that condition (i) implies conditions (1) and (2). Characterizations of cover relations in terms of quiver representations, equivalent to conditions (i) and (ii), can be found in \cite{S07}.

%We include a complete proof of the Theorem using only combinatorics and linear algebra in the Appendix.  
\end{proof}

The following immediate corollary is also standard:

\begin{corollary}
Suppose that the orbit $\mathbb{O}_M\subset \Fnd\times \Fnd$ is closed. Then there exist sequences $1\le a_1\le\ldots\le a_k\le n$ and $1\le b_1\le\ldots\le b_k\le n$ such that for all $1\le i<n$ $(a_i,b_i)\neq(a_{i+1},b_{i+1}),$ and $m_{ij}\neq 0$ if and only if $(i,j)=(a_l,b_l)$ for some $1\le l\le k.$ 
\end{corollary}

Furthermore, one immediately gets a geometric description of closed orbits in $\Fnd\times\Fnd,$ which appears to be new:

\begin{corollary}\label{corollary: closed orbits}
In the notations as above, the closed orbit $\mathbb{O}_M$ is then the image of the injective map $\phi_M:\Fmu \to \Fnd\times \Fnd,$ where $\mu=(m_{a_1b_1},\ldots,m_{a_kb_k})$ and the map $\phi$ is given by $\phi(\textbf{V})=(\textbf{U},\textbf{W})$ where for all $i$
\begin{equation*}
U_i:=V_{\max \{j|a_j\le i\}}
\end{equation*}
and
\begin{equation*}
W_i:=V_{\max \{j|b_j\le i\}}.
\end{equation*}
\end{corollary}

Similarly, for open orbits one gets:

\begin{corollary}\label{corollary: open orbits}
Suppose that the orbit $\mathbb{O}_M\subset \Fnd\times \Fnd$ is open. Then there exist sequences $1\le a_1\le\ldots\le a_k\le n$ and $n\ge b_1\ge\ldots\ge b_k\ge 1$ such that for all $1\le i<n$ $(a_i,b_i)\neq(a_{i+1},b_{i+1}),$ and $m_{ij}\neq 0$ if and only if $(i,j)=(a_l,b_l)$ for some $1\le l\le k.$ %Moreover, the orbit $\mathbb{O}_M$ is then the image of the injective map $\phi_M:\textrm{Split}_\mu \to \Fnd\times \Fnd,$ where $\mu=(m_{a_1b_1},\ldots,m_{a_kb_k}),$ $\textrm{Split}_\mu$ is the space of direct sum decompositions 
%\begin{equation*}
%\textrm{Split}_\mu:=\{V=V_1\oplus\ldots\oplus V_k|\dim V_i=m_{a_ib_i}\},
%\end{equation*} 
%and the map $\phi_M$ is given by $\phi_M(V_1\oplus V_2\oplus\ldots\oplus V_k)=(\textbf{U},\textbf{W}),$ where for all $i$
%\begin{equation*}
%U_i:=\bigoplus\limits_{a_j\le i} V_j
%\end{equation*}
%and
%\begin{equation*}
%W_i:=\bigoplus\limits_{b_j\le i} V_j.
%\end{equation*}
\end{corollary}

\begin{remark} 
Note that for $\mathcal{M}_{1^d,1^d}$ the partial order considered above coincides with the Bruhat order on the permutation matrices. Moreover, in the general case the order can be obtained from the Bruhat order in the following manner. Let $\mu=(\mu_1,\ldots,\mu_n)$ and $\nu=(\nu_1,\ldots,\nu_m)$ be two compositions of $d.$ Consider the map
\begin{equation*}
p:\mathcal{M}_{1^d,1^d}\to \mathcal{M}_{\mu,\nu}
\end{equation*}
given by 
\begin{equation*}
p(M)_{kl}=\sum\limits_{\substack{
\mu_{k-1}<i\le \mu_k\\
\nu_{l-1}<j\le\nu_l}}
m_{ij}.
\end{equation*}

The partial order on $\mathcal{M}_{\mu,\nu}$ is uniquely defined by requiring that the map $p$ is weakly monotone %(note that the image of incomparable matrices might be comparable).
Note that $p$ is the combinatorial shadow of the forgetful map 
\begin{equation*}
\pi: \mathcal{F}_{1^d}\times \mathcal{F}_{1^d} \to \Fmu\times \mathcal{F}_{\nu},
\end{equation*}
i.e. $\pi(\mathbb{O}_M)=\mathbb{O}_{p(M)}$ for all $M\in\mathcal{M}_{1^d,1^d}.$

\end{remark}

\subsection{The support of a convolution product}\label{subsection: convolution product for supports} In order to study the convolution algebra 
$$\Snd=K^G(\Fnd\times \Fnd)$$
we need to first study the behavior of supports of the K-theory classes under the convolution product. Let $\mathbb{O}_M,\mathbb{O}_N\subset \Fnd\times \Fnd$ be two orbits and $M,N$ be the corresponding matrices. Let also $\alpha,\beta \in K^G(\Fnd)$ be two classes such that $\supp(\alpha)=\overline{\mathbb{O}_M}$ and $\supp(\beta)=\overline{\mathbb{O}_N}.$ Then the convolution product $\alpha\star\beta$ is supported on the closure of $\pi_{13}(\pi_{12}^{-1}(\mathbb{O}_M)\cap\pi_{23}^{-1}(\mathbb{O}_N)).$ Note that $\pi_{13}(\pi_{12}^{-1}(\mathbb{O}_M)\cap\pi_{23}^{-1}(\mathbb{O}_N))$ is non-empty if and only if the second flag in $\mathbb{O}_M$ is of the same type as the first flag in $\mathbb{O}_N,$ which is equivalent to $\textbf{C}^M=\textbf{R}^N.$

\begin{lemma}\cite[Lemma 6.1]{JS15}\label{lemma:irreducibility}
Suppose that $\textbf{C}^M=\textbf{R}^N.$ Then the subvarieties $\pi_{12}^{-1}(\mathbb{O}_M)\cap\pi_{23}^{-1}(\mathbb{O}_N)\subset \Fnd\times \Fnd\times \Fnd$ and, therefore, $\pi_{13}(\pi_{12}^{-1}(\mathbb{O}_M)\cap\pi_{23}^{-1}(\mathbb{O}_N))\subset \Fnd\times \Fnd$ are irreducible.
\end{lemma} 

Lemma \ref{lemma:irreducibility} motivates the following definition:

\begin{definition}[\cite{BLM90,V98,JS15}]
If $\textbf{C}^M=\textbf{R}^N,$ then we define $M\circ N$ to be the matrix such that $\mathbb{O}_{M\circ N}\subset \pi_{13}(\pi_{12}^{-1}(\mathbb{O}_M)\cap\pi_{23}^{-1}(\mathbb{O}_N))$ is the unique orbit that is open in $\pi_{13}(\pi_{12}^{-1}(\mathbb{O}_M)\cap\pi_{23}^{-1}(\mathbb{O}_N)).$ If $\textbf{C}^M\neq\textbf{R}^N$ then one sets $M\circ N:=0.$
\end{definition}

The product $M\circ N$ was studied in \cite{BLM90} and \cite{V98}.  In \cite{JS15} Jensen and Su studied this product using the quiver representation techniques. Here we establish basic properties of this operation. We mostly follow \cite{BLM90} and \cite{V98}, although we provide more details and some of the formulations are new.

In general, computing the product $M\circ N$ is rather complicated. However, if one of the matrices is of a specially simple form, then one can describe $M\circ N$ explicitly. 

\begin{definition}
An orbit $\mathbb{O}\subset \Fnd\times \Fnd\times \Fnd$ of the diagonal $GL(V)$ action is called \textit{combinatorial} if for every triple of flags $(\textbf{U},\textbf{V},\textbf{W})\in \mathbb{O}$ there exists a basis of $V$ such that all elements of flags $\textbf{U},\textbf{V},$ and $\textbf{W}$ are coordinate  subspaces.
\end{definition}

\begin{remark}
Not all orbits in the variety of triples of partial flags are combinatorial: for example, consider three distinct lines in a two dimensional space. Moreover, in most of the cases the variety of triples of partial flags contains infinitely many orbits of the diagonal $GL(V)$--action. See \cite{MWZ99} for a classification of multiple flag varieties with finitely many orbits.
\end{remark}

\begin{lemma}\label{lemma: orbits in preimage of closed orbits}
Let $M=\{m_{ij}\}\in\mathcal{M}(n,d)$ be a matrix such that $\mathbb{O}_M\subset \Fnd\times \Fnd$ is a closed orbit. Then for any $(l,m)\in\{(1,2),(1,3),(2,3)\}$ the preimage $\pi_{lm}^{-1}(\mathbb{O}_M)\subset \Fnd\times \Fnd\times \Fnd$ is a union of combinatorial orbits. 
\end{lemma}

\begin{proof}
 According to Corollary \ref{corollary: closed orbits}, there exist two sequences $1\le a_1\le\ldots\le a_k\le n$ and $1\le b_1\le\ldots\le b_k\le n,$ such that $m_{a_1b_1},\ldots,m_{a_kb_k}$ are the only non-zero entries of $M,$ and the map 
\begin{equation*}
\phi:\Fmu \to \mathbb{O}_M,
\end{equation*}
where $\mu=(m_{a_1b_1},\ldots,m_{a_kb_k}),$ given by $\phi(\textbf{V})=(\textbf{U},\textbf{W})$ where
\begin{equation*}
U_i:=\textbf{V}_{\max \{j|a_j\le i\}}
\end{equation*}
and
\begin{equation*}
W_i:=\textbf{V}_{\max \{j|b_j\le i\}}
\end{equation*}
is an isomorphism.

Let $(\textbf{U},\textbf{W})\in \mathbb{O}_M.$ Let $\textbf{V}=\phi^{-1}(\textbf{U},\textbf{W})\in \Fmu.$ %Then 
%\begin{equation*}
%\textbf{U}_i:=\textbf{V}_{\max \{j|a_j\le i\}}
%\end{equation*}
%and
%\begin{equation*}
%\textbf{W}_i:=\textbf{V}_{\max \{j|b_j\le i\}}.
%\end{equation*}
Consider any triple of flags $(\textbf{U},\textbf{W},\textbf{T})\in \pi_{12}^{-1}(\mathbb{O}_M).$ Let $\mathcal{B}$ be a basis of $V,$ such that all elements of flags $\textbf{V}$ and $\textbf{T}$ are spans of basic vectors. But then the elements of the flags $\textbf{U}$ and $\textbf{W}$ are also spans of basic vectors, as they coincide with elements of the flag $\textbf{V}.$ The preimages under $\pi_{13}$ and $\pi_{23}$ are treated in the same way.
\end{proof}

Similar to the case of the variety of pairs of partial flags, combinatorial orbits in the variety of triples of partial flags can be enumerated by three dimensional arrays of non-negative integers. However, the proof of Lemma \ref{lemma: orbits in preimage of closed orbits} suggests that the orbits in $\pi_{ij}^{-1}(\mathbb{O}_M)$ for a closed orbit $\mathbb{O}_M$ can be enumerated in an easier way. Let $\phi: \Fmu\xrightarrow{\sim} \mathbb{O}_M\subset \Fnd\times \Fnd$ be the equivariant isomorphism from Corollary \ref{corollary: closed orbits}. Let us extend it to the equivariant isomorphism
\begin{equation*}
\Phi:=(\phi,id): \Fmu\times \Fnd\xrightarrow{\sim} \mathbb{O}_M\times \Fnd\simeq \pi_{12}^{-1}(\mathbb{O}_M)\subset \Fnd\times \Fnd\times \Fnd.
\end{equation*}

The orbits in $\Fmu\times \Fnd$ are enumerated by $k\times n$ non-negative integer matrices $N=\{n_{ij}\},$ such that $\textbf{R}^N=\mu.$ Furthermore, the following result follows immediately from Corollary \ref{corollary: closed orbits}:

\begin{lemma}%\label{lemma: orbits in preimage of closed orbits}
Let $L=\{l_{ij}\}\in\mathcal{M}(n,d)$ be the matrix corresponding to the orbit $\pi_{13}(\Phi(\mathbb{O}_N))\subset \Fnd\times \Fnd.$ Then one gets
\begin{equation*}
l_{kl}=\sum\limits_{\{i|a_i=k\}} n_{il}.
\end{equation*}
In other words, the $k$th row of the matrix $L$ is obtained from the matrix $N$ by adding up rows corresponding to the non-zero entries of the $k$th row of the matrix $M.$

Similarly, if $K=\{k_{ij}\}\in\mathcal{M}(n,d)$ is the matrix corresponding to the orbit $\pi_{23}(\Phi(\mathbb{O}_N))\subset \Fnd\times \Fnd.$ Then one gets
\begin{equation*}
k_{kl}=\sum\limits_{\{i|b_i=k\}} n_{il},
\end{equation*}
i.e. the $k$th row of the matrix $K$ is obtained from the matrix $N$ by adding up rows corresponding to the non-zero entries of the $k$th column of the matrix $M.$
\end{lemma}

\begin{example}\label{example: comb orbits}
Consider the matrix
\begin{equation*}
M=\left(
\begin{array}{lll}
2 &1 &0 \\
0 &2 &0 \\
0 &1 &3
\end{array}
\right)\in \mathcal{M}(3,9).
\end{equation*}
Clearly, $M$ satisfies the conditions of Corollary \ref{corollary: closed orbits} with $\mu=(m_{11},m_{12},m_{22},m_{32},m_{33})=(2,1,2,1,3).$ We have the isomorphism $\Phi:\Fmu\times \mathcal{F}^9_3\xrightarrow{\sim}\pi_{12}^{-1}(\mathbb{O}_M)$ with orbits in $\Fmu\times \mathcal{F}^9_3$ enumerated by $5\times 3$ matrices with row-sums given by $\mu.$ For example, the matrix
\begin{equation*}
N=\left(
\begin{array}{lll}
0 &1 &1 \\
0 &1 &0 \\
1 &1 &0 \\
1 &0 &0 \\
0 &1 &2
\end{array}
\right)
\end{equation*}
corresponds to one of the orbits in $\Fmu\times \mathcal{F}_n.$ Furthermore, the orbit $\pi_{13}(\Phi(\mathbb{O}_N))\subset \Fnd\times \Fnd$ corresponds to the matrix  
\begin{equation*}
L=\left(
\begin{array}{lll}
0 &2 &1 \\
1 &1 &0 \\
1 &1 &2
\end{array}
\right),
\end{equation*}
obtained from $N$ by summing up the first two rows, and also the last two rows, since $M$ has two non-zero entries in the first row, and two non-zero entries in the last row. Similarly, the orbit $\pi_{23}(\Phi(\mathbb{O}_N))\subset \Fnd\times \Fnd$ corresponds to the matrix  
\begin{equation*}
K=\left(
\begin{array}{lll}
0 &1 &1 \\
2 &2 &0 \\
0 &1 &2
\end{array}
\right),
\end{equation*}
obtained from $N$ by summing up rows two, three, and four, since $M$ has three non-zero entries in the second column.
\end{example}

\begin{definition}\label{definition: almost diagonal}
Let $\mu=(\mu_1,\ldots,\mu_n)$ be a composition of $d-1,$ i.e. $\mu\in\mathbb{Z}_{\ge 0}^n$ and $\sum \mu_i=d-1.$ For $1\le k<n$ we will use the following notations
\begin{equation*}
\textbf{E}(\mu,k):=Diag(\mu)+E_{k,k+1},
\end{equation*} 
and
\begin{equation*}
\textbf{F}(\mu,k):=Diag(\mu)+E_{k+1,k},
\end{equation*} 
where $Diag(\mu)$ is the diagonal $n\times n$ matrix with numbers $\mu_1,\ldots, \mu_n$ along the diagonal. We are going to call such matrices \textit{almost diagonal}.
\end{definition}

\begin{definition}
We will use the notation $\textbf{e}_k=(0,\ldots,0,1,0,\ldots,0)\in \mathbb{Z}^n$ for the $k$th basic vector.
\end{definition}

Let $\textbf{E}(\mu,k)$ be an almost diagonal matrix, and let $M=\{m_{ij}\}\in\mathcal{M}_{d,n},$ such that $\textbf{R}^M=\textbf{C}^{E(\mu,k)}=\mu+\textbf{e}_{k+1}.$ Let also $\nu=(\mu_1,\ldots, \mu_k,1,\mu_{k+1},\ldots,\mu_n).$ Then the orbits in the preimage $\pi_{12}^{-1}(\mathbb{O}_{\textbf{E}(\mu,k)})\simeq \mathcal{F}_\nu\times \Fnd$ are enumerated by $(n+1)\times n$ non-negative integer matrices $N$ such that $\textbf{R}^N=\nu.$ Intersecting with the preimage $\pi_{23}^{-1}(\mathbb{O}_M)$ we get the extra condition that if one replaces $(k+1)$th and $(k+2)$th rows of $N$ by their sum, one obtains $M.$ 

The orbits in $\pi_{13}(\pi_{12}^{-1}(\mathbb{O}_{\textbf{E}(\mu,k)})\cap\pi_{23}^{-1}(\mathbb{O}_M))$ are then enumerated by matrices $L$ obtained from $N$ by replacing $k$th and $(k+1)$th rows by their sum. In other words, these orbits are enumerated by
\begin{equation*}
\mathcal{L}:=\{M-E_{k+1,j}+E_{k,j}|m_{k+1,j}>0\}.
\end{equation*}
By Theorem \ref{Theorem: partial order and cover relation}, the set $\mathcal{L}$ is totally ordered with the maximal element given by $M-E_{k+1,m}+E_{k,m},$ where $m=\max\{j|m_{k+1,j}>0\}.$ In particular, one gets the following standard result:

\begin{theorem}[\cite{BLM90,V98}]\label{theorem: matrix mult by almost diag}
With the notations as above, one has
\begin{equation*}
E(\mu,k)\circ M=M-E_{k+1,m}+E_{k,m},
\end{equation*}
where $m=\max\{j|m_{k+1,j}>0\}.$

Similarly, one gets $\textbf{F}(\mu,k)\circ M=0$ unless $\textbf{R}^M=\textbf{C}^{\textbf{F}(\mu,k)}=\mu+\textbf{e}_k,$ in which case one has
\begin{equation*}
F(\mu,k)\circ M=M-E_{k,m}+E_{k+1,m},
\end{equation*}
where $m=\min\{j|m_{k,j}>0\}.$
\end{theorem}

%\begin{proof}
%The Theorem follows immediately from Lemma \ref{lemma: orbits in preimage of closed orbits} and Theorem \ref{Theorem: partial order and cover relation}. Indeed, let $\nu=(\mu_1,\ldots,\mu_k,1,\mu_{k+1},\ldots,\mu_n).$ Then the orbits in the preimage $\pi_{12}^{-1}(\mathbb{O}_{\textbf{E}(\mu,k)})\simeq \mathcal{F}_\nu\times \Fnd$ are enumerated by $(n+1)\times n$ non-negative integer matrices $N$ such that $\textbf{R}^N=\nu.$ Intersecting with the preimage $\pi_{23}^{-1}(\mathbb{O}_M)$ we get the extra condition that if one replaces $(k+1)$th and $(k+2)$th rows of $N$ by their sum, one obtains $M.$ 

%The orbits in $\pi_{13}(\pi_{12}^{-1}(\mathbb{O}_{\textbf{E}(\mu,k)})\cap\pi_{23}^{-1}(\mathbb{O}_M))$ are then enumerated by matrices $L$ obtained from $N$ by replacing $k$th and $(k+1)$th rows by their sum. In other words, these orbits are enumerated by
%\begin{equation*}
%\mathcal{L}:=\{M-E_{k+1,j}+E_{k,j}|m_{k+1,j}>0\}.
%\end{equation*}
%By Theorem \ref{Theorem: partial order and cover relation}, the set $\mathcal{L}$ is totally ordered with the maximal element given by $M-E_{k+1,m}+E_{k,m},$ where $m=\max\{j|m_{k+1,j}>0\}.$ Therefore,
%\begin{equation*}
%\textbf{E}(\mu,k)\circ M=M-E_{k+1,m}+E_{k,m}.
%\end{equation*}
%The second part of the Theorem is verified in a similar manner.
%\end{proof}

\begin{corollary}[\cite{BLM90,V98}]\label{corollary:matrix factorization}
Every non-diagonal matrix $M=\{m_{ij}\}\in\mathcal{M}(n,d)$ can be factored into a $\circ$-product of several almost diagonal matrices.
\end{corollary}

This result is well known, but we present a new proof for completeness of the exposition. Our proof uses induction with respect to the \textit{diagonal norm} of a matrix defined as follows.

\begin{definition}\label{definition: dnorm}
For a matrix $M\in\mathcal{M}(n,d)$ define
\begin{equation*}
    \dnorm{M}:=\sum\limits_{1\le i,j\le n} |i-j|m_{ij}.
\end{equation*}
\end{definition}

\begin{proof}
%Consider the \textit{diagonal norm} function $DN:\mathcal{M}(n,d)\to \mathbb{Z}_{\ge 0}$ given by
%\begin{equation*}
%DN(M):=\sum\limits_{1\le i,j\le n} |i-j|m_{ij}.
%\end{equation*}
The proof proceeds by induction with respect to $\dnorm{M}.$ One has $\dnorm{M}=0$ if and only if $M$ is diagonal, and $\dnorm{M}=1$ if and only if $M$ is almost diagonal, in which cases the Corollary is trivial. Suppose that $\dnorm{M}>1$ and the Corollary is proven for all $N\in\mathcal{M}(n,d)$ such that $\dnorm{N}<\dnorm{M}.$ Then there exist $i\neq j$ such that $m_{ij}>0$. 

Suppose that $j>i.$ Then there exist $b=\max\{j|\ \exists i<j, m_{ij}>0\}$ and $a<b$ such that $m_{ab}>0.$ Applying Theorem \ref{theorem: matrix mult by almost diag} one gets
\begin{equation*}
M=\textbf{E}(\textbf{R}^M-\textbf{e}_a,a)\circ (M+E_{a+1,b}-E_{a,b})
\end{equation*} 
and $\dnorm{M+E_{a-1,b}-E_{a,b}}=\dnorm{M}-1.$

Similarly, if $j<i,$ then there exists $b=\min\{j|\exists i>j, m_{ij}>0\}$ and $a>b$ such that $m_{ab}>0.$ Applying Theorem \ref{theorem: matrix mult by almost diag} one gets
\begin{equation*}
M=\textbf{F}(\textbf{R}^M-\textbf{e}_a,a-1)\circ (M+E_{a-1,b}-E_{a,b})
\end{equation*} 
and $\dnorm{M+E_{a,b}-E_{a-1,b}}=\dnorm{M}-1.$
\end{proof}

\begin{corollary}\label{corollary: irreducibility with closure}
    Suppose that $\textbf{E}(\mu,k)\circ M\neq 0,$ then the intersection of the preimages $$\pi_{12}^{-1}(\mathbb{O}_{\textbf{E}(\mu,k)})\cap\pi_{23}^{-1}(\overline{\mathbb{0}_M})$$ is irreducible.
\end{corollary}

\begin{proof}
    Indeed, from the description of the Bruhat order and cover relations (Theorem \ref{Theorem: partial order and cover relation}), orbits in $\pi_{12}^{-1}(\mathbb{O}_{\textbf{E}(\mu,k)})\simeq \mathcal{F}_\nu\times \Fnd,$ and their images under $\pi_{23}$ it is clear that for any two orbits $\mathbb{O}_K,\mathbb{O}_L\subset\overline{\mathbb{O}_M}$ such that $K\lessdot L$ and an orbit $\mathbb{O}_{K'}\subset\pi_{12}^{-1}(\mathbb{O}_{\textbf{E}(\mu,k)})$ such that $\pi_{23}(\mathbb{O}_{K'})=\mathbb{O}_K,$ there exists an orbit $\mathbb{O}_{L'}\subset\pi_{12}^{-1}(\mathbb{O}_{\textbf{E}(\mu,k)})$ such that $\pi_{23}(\mathbb{O}_{L'})=\mathbb{O}_L$ and $K'\prec L'.$ Furthermore, $\pi_{12}^{-1}(\mathbb{O}_{\textbf{E}(\mu,k)})\cap\pi_{23}^{-1}(\mathbb{0}_M)$ is irreducible by Lemma \ref{lemma:irreducibility}.
\end{proof}

Let $M=\{m_{ij}\}, N=\{n_{ij}\}\in\mathcal{M}(n,d)$ and $M=\textbf{E}(\mu,k)\circ N.$ We know then that 

\begin{equation*}
N=M+E_{k+1,b}-E_{k,b},
\end{equation*}
where $b$ is such that $m_{kb}>0$ and $m_{k+1,j}=0$ for all $j>b.$ However, given a matrix $M,$ the set
\begin{equation*}
\{L\in\mathcal{M}(n,d)|M=\textbf{E}(\mu,k)\circ L\}
\end{equation*}
might consist of more than one matrix (see Example \ref{example: non unique quotient} below).

\begin{example}\label{example: non unique quotient}
One has

\begin{align*}
\left(
\begin{array}{llll}
1 &1 &0 &2 \\
0 &1 &0 &0 \\
0 &0 &0 &0 \\
0 &0 &0 &0 
\end{array}
\right)&=
\left(
\begin{array}{llll}
3 &1 &0 &0 \\
0 &1 &0 &0 \\
0 &0 &0 &0 \\
0 &0 &0 &0 
\end{array}
\right)\circ
\left(
\begin{array}{llll}
1 &0 &0 &2 \\
0 &2 &0 &0 \\
0 &0 &0 &0 \\
0 &0 &0 &0 
\end{array}
\right)\\
&=
\left(
\begin{array}{llll}
3 &1 &0 &0 \\
0 &1 &0 &0 \\
0 &0 &0 &0 \\
0 &0 &0 &0 
\end{array}
\right)\circ
\left(
\begin{array}{llll}
1 &1 &0 &1 \\
0 &1 &0 &1 \\
0 &0 &0 &0 \\
0 &0 &0 &0 
\end{array}
\right)
\end{align*}
\end{example}

The following Lemma follows immediately from Theorems \ref{theorem: matrix mult by almost diag} and \ref{Theorem: partial order and cover relation}:

\begin{lemma}\label{lemma: non-unique division}
Let $M=\{m_{ij}\}\in\mathcal{M}(n,d)$ be a matrix and $\textbf{E}(\mu,k)$ be such that the set
\begin{equation*}
\mathcal{S}:=\{L\in\mathcal{M}(n,d)|M=\textbf{E}(\mu,k)\circ L\}
\end{equation*}
is not empty. Then $\mathcal{S}$ is totally ordered, and the minimum of $\mathcal{S}$ is given by
\begin{equation*}
\min\mathcal{S}=M+E_{k+1,b}-E_{k,b},
\end{equation*}
where $b$ is the maximal integer satisfying $m_{kb}>0$ and $m_{k+1,j}=0$ for all $j>b.$

Similarly, if $\textbf{F}(\mu,k)$ is such that the set
\begin{equation*}
\mathcal{T}:=\{L\in\mathcal{M}(n,d)|M=\textbf{F}(\mu,k)\circ L\}
\end{equation*}
is not empty, then $\mathcal{T}$ is totally ordered, and the minimum of $\mathcal{T}$ is given by
\begin{equation*}
\min\mathcal{T}=M-E_{k+1,b}+E_{k,b},
\end{equation*}
where $b$ is the minimal integer satisfying $m_{k+1,b}>0$ and $m_{k,j}=0$ for all $j<b.$
\end{lemma}

\begin{corollary}\label{corollary:strictly monotone factorization}
For any matrix $M\in\mathcal{M}(n,d)$ such that $\dnorm{M}>0$ (i.e. $M$ is not diagonal, see Definition \ref{definition: dnorm}) there exists a factorization $M=D\circ N$ such that $D$ is almost diagonal, $\dnorm{N}<\dnorm{M},$ and for any $L\prec N$ one has $D\circ L\prec M.$
\end{corollary}

The following Lemma will be important in the next section:

\begin{lemma}\label{lemma: triple intersection}
Let $\mathbb{O}_1,\mathbb{O}_2,\mathbb{O}_3\subset \Fnd\times \Fnd$ be three orbits. Suppose that at least one of them is closed. Then the triple intersection $\pi_{12}^{-1}(\mathbb{O}_1)\cap\pi_{23}^{-1}(\mathbb{O}_2)\cap\pi_{13}^{-1}(\mathbb{O}_3)\subset \Fnd\times \Fnd\times \Fnd$ consists of a single orbit.
\end{lemma}

\begin{remark}
Note that the condition that one of the orbits is closed is essential: see Examples \ref{example: triple intersection 1} and \ref{example: triple intersection 2}.
\end{remark}

\begin{proof}
Let $\mathbb{O}_1$ be the closed orbit, then by Corollary \ref{corollary: closed orbits} we have $\mathbb{O}_1=\mathbb{O}_D$ where $D=\{d_{ij}\}\in\mathcal{M}(n,d)$ is such that there exist $a_1\le\ldots\le a_k$ and $b_1\le\ldots\le b_k$ and $\mu:=(\mu_1,\ldots,\mu_k)$ such that $d_{a_lb_l}=\mu_l>0$ for all $l$ and all other entries of $D$ are zero. Here for every $1\le l<k$ one has $(a_l,b_l)\neq (a_{l+1},b_{l+1}),$ i.e. either $a_l<a_{l+1}$ or $b_l<b_{l+1}.$ Let also $\mathbb{O}_2=\mathbb{O}_M$ and $\mathbb{O}_3=\mathbb{O}_N$ where $M,N\in\mathcal{M}(n,d).$

Using Lemma \ref{lemma: orbits in preimage of closed orbits} we conclude that all orbits in the triple intersection are combinatorial. Moreover, for every such orbit $\mathbb{O}$ there exists a $k\times n$ matrix $L$ such that $\mathbb{O}$ is isomorphic to $\mathbb{O}_L\subset \mathcal{F}_k^d\times \Fnd$ under the embedding
\begin{equation*}
\phi_D\times \id: \mathcal{F}_k^d\times \Fnd\to \Fnd\times \Fnd\times \Fnd,
\end{equation*}  
where $\phi_D: \mathcal{F}_k^d \to \Fnd\times \Fnd$ is given by
\begin{equation*}
\phi_D(\textbf{V}):=(\textbf{U},\textbf{W}),
\end{equation*}
where
\begin{equation*}
U_l:=V_{\max\{m|a_m\le l\}},
\end{equation*}
and  
\begin{equation*}
W_l:=V_{\max\{m|b_m\le l\}}.
\end{equation*}

In particular, matrix $L$ uniquely determines the orbit in the triple intersection, and matrices $M$ and $N$ can be obtained from $L$ as follows. For every $l$ one gets
\begin{equation*}
R^M_l=\sum_{a_i=l} R^L_i,
\end{equation*}
and
\begin{equation*}
R^N_l=\sum_{b_i=l} R^L_i.
\end{equation*}

Let us show that these conditions determine matrix $L$ uniquely. Suppose that there are two matrices $L$ and $L'$ satisfying these conditions. Let $i$ be the smallest number, such that $i$th rows of $L$ and $L'$ are not the same $R^L_i\neq R^{L'}_i$. We have $(a_i,b_i)\neq (a_{i+1},b_{i+1})$ by construction. Without loss of generality, let us assume that $a_i<a_{i+1}.$ Let $m=\min\{j|a_j=a_i\}.$ Then one gets
\begin{align*}
R^M_{a_i}&=R^L_m+R^L_{m+1}+\ldots+R^L_i\\
&=R^{L'}_m+R^{L'}_{m+1}+\ldots+R^{L'}_i.
\end{align*} 
But $R^L_j=R^{L'}_j$ for $j<i.$ Therefore, $R^L_i=R^{L'}_i.$ Contradiction.
\end{proof}

We finish the section with a few examples involving $\circ$-multiplication of matrices.

\begin{example}
Consider the matrix
\begin{equation*}
M=\{m_{ij}\}=
\left(
\begin{array}{ll}
0 &1 \\
1 &0
\end{array}
\right)\in\mathcal{M}(2,2).
\end{equation*}
Then the orbit $\mathbb{O}_M\subset \mathcal{F}_2^2\times \mathcal{F}_2^2$ consists of pairs of distinct lines in a two dimensional space $V$. 

The intersection of the preimages $\pi_{12}^{-1}(\mathbb{O}_M)\cap\pi_{23}^{-1}(\mathbb{O}_M)\subset \mathcal{F}_2^2\times \mathcal{F}_2^2\times \mathcal{F}_2^2$ consists of triples of lines $(l_1,l_2,l_3),$ such that $l_1\neq l_2$ and $l_2\neq l_3.$ There are two orbits in this intersection: $\mathbb{O}_1$ consisting of triples of distinct lines, and $\mathbb{O}_2$ consisting of triples of lines where the first and the third coincide, while the second is different from them. Note that orbit $\mathbb{O}_1$ is not combinatorial. The image of the intersection $\pi_{13}(\pi_{12}^{-1}(\mathbb{O}_M)\cap\pi_{23}^{-1}(\mathbb{O}_M))\subset \mathcal{F}_2^2\times \mathcal{F}_2^2$ also consists of two orbits: $\mathbb{O}_M$ and $\mathbb{O}_L,$ where
\begin{equation*}
L=
\left(
\begin{array}{ll}
1 &0 \\
0 &1
\end{array}
\right).
\end{equation*}
Clearly, $L\prec M,$ therefore $M\circ M=M.$ Obtaining this result in a direct combinatorial way is tricky because $\mathbb{O}_1$ is not a combinatorial orbit. In particular, there does not exist a $2\times 2\times 2$ array of non-negative integers $N=\{n_{ijk}\},$ such that for all $1\le i,j\le 2$
\begin{equation*}
n_{1ij}+n_{2ij}=n_{i1j}+n_{i2j}=n_{ij1}+n_{ij2}=m_{ij}.
\end{equation*}

On the other hand, using Theorem \ref{theorem: matrix mult by almost diag} and the factorization
\begin{equation*}
\left(
\begin{array}{ll}
0 &1 \\
1 &0
\end{array}
\right)=
\left(
\begin{array}{ll}
0 &1 \\
0 &1
\end{array}
\right)
\circ
\left(
\begin{array}{ll}
0 &0 \\
1 &1
\end{array}
\right),
\end{equation*}
one computes
\begin{align*}
&\left(
\begin{array}{ll}
0 &1 \\
1 &0
\end{array}
\right)
\circ
\left(
\begin{array}{ll}
0 &1 \\
1 &0
\end{array}
\right)=
\left[\left(
\begin{array}{ll}
0 &1 \\
0 &1
\end{array}
\right)
\circ
\left(
\begin{array}{ll}
0 &0 \\
1 &1
\end{array}
\right)\right]
\circ
\left(
\begin{array}{ll}
0 &1 \\
1 &0
\end{array}
\right)\\
&=
\left(
\begin{array}{ll}
0 &1 \\
0 &1
\end{array}
\right)
\circ
\left[\left(
\begin{array}{ll}
0 &0 \\
1 &1
\end{array}
\right)
\circ
\left(
\begin{array}{ll}
0 &1 \\
1 &0
\end{array}
\right)\right]
=\left(
\begin{array}{ll}
0 &1 \\
0 &1
\end{array}
\right)
\circ
\left(
\begin{array}{ll}
0 &0 \\
1 &1
\end{array}
\right)=
\left(
\begin{array}{ll}
0 &1 \\
1 &0
\end{array}
\right).
\end{align*}
\end{example}

\begin{example}\label{example: triple intersection 1}
Consider the orbit $\mathbb{O}_M$ where 
\begin{equation*}
M:=\left(
\begin{array}{cc}
0 & 1 \\
1 & 1 
\end{array}
\right),
\end{equation*}
i.e. the space of pairs of distinct lines in a three dimensional space $V$. Then the triple intersection $\pi_{12}^{-1}(\mathbb{O}_M)\cap \pi_{23}^{-1}(\mathbb{O}_M)\cap \pi_{13}^{-1}(\mathbb{O}_M)$ consist of triples of distinct lines in $V,$ which consists of two orbits: triples of lines that form a direct sum, and triples of distinct lines that lie in the same plane. 
\end{example}

\begin{example}\label{example: triple intersection 2}
In the previous example one of the orbits in the triple intersection was not combinatorial (triples of lines that lies in the same plane). However, it is not hard to come up with an example of a triple intersection containing more than one combinatorial orbit. This is equivalent to finding two distinct three dimensional arrays of non-negative integers $M=\{m_{ijk}\}$ and $N=\{n_{ijk}\}$ such that for every $a$ and $b$ one has
\begin{equation*}
\sum_i m_{iab}=\sum_i n_{iab},
\end{equation*}
\begin{equation*}
\sum_j m_{ajb}=\sum_j n_{ajb},
\end{equation*}
and
\begin{equation*}
\sum_k m_{abk}=\sum_k n_{abk}.
\end{equation*}
For example, one can take $2\times 2\times 2$ arrays given by $m_{ijk}=1$ for all $i,j,k,$ and 
\begin{equation*}
n_{ijk}=\left\{
\begin{array}{l}
2, \ i+j+k\ \textrm{is even},\\
0, \ i+j+k\ \textrm{is odd}.
\end{array}
\right.
\end{equation*}

\end{example}

\section{Equivariant K-theory and the convolution algebras $\Snd.$}\label{section: convolution algebra}

\subsection{Equivariant K-theory of the variety of partial flags}\label{subsection: K-theory of F}

Let $\mu$ be a composition of $d$ of length $n.$ Consider the variety of partial flags $\Fmu\simeq G/P_{\mu}$ where the parabolic subgroup $P_{\mu}\subset G$ is the stabilizer of a partial flag $p\in \Fmu.$ Note that $P_{\mu}\simeq L_{\mu}\ltimes U_{\mu}$ where $L_{\mu}\simeq GL_{\mu_1}\times\ldots\times GL_{\mu_n}\subset P_{\mu}$ is the Levi subgroup and $U_{\mu}$ is the unipotent radical of $P_{\mu}.$ One gets 

\begin{equation*}
K^G(\Fmu)\simeq K^{P_{\mu}}(pt)\simeq K^{L_{\mu}}(pt)\simeq R(L_\mu)\simeq C[x_1^{\pm 1},\ldots,x_d^{\pm 1}]^{S_{\bf\mu}}\simeq \bigotimes_i \Lambda^{\pm}_{\mu_i},
\end{equation*}
\noindent where $S_{\mu}=S_{\mu_1}\times\ldots\times S_{\mu_n}\subset S_d$ is the corresponding Young subgroup, and $\Lambda^{\pm}_t$ is the ring of symmetric Laurent polynomials in $t$ variables.

%\begin{lemma}
%One has the following isomorphism for the $G$-equivariant $K$-theory of $\Fmu:$
%$$
%K^G(\Fmu)\simeq \mathbb C[x_1^{\pm 1},\ldots,x_d^{\pm 1}]^{S_{\bf\mu}},
%$$
%\noindent where $S_{\mu}=S_{\mu_1}\times\ldots\times S_{\mu_n}\subset S_d.$ 
%\end{lemma}

%\begin{proof}
This isomorphism can be understood in the following way. Let $\pi:E\to \Fmu$ be a $G$-equivariant vector bundle. %Since $\Fmu$ is a homogeneous space, the whole bundle is determined by the fiber over just one point, and the action of the stabilizer on that fiber. 
Let us fix a basis $\mathcal{B}=\{e_1,\ldots,e_d\}\subset V,$ let $T\subset G$ be the maximal torus acting by scaling the coordinates in this basis, and let $\{x_1,\ldots,x_d\}$ be the corresponding characters of $T,$ forming a basis in the character lattice of $T.$ The $T$-fixed points in $\Fmu$ are the flags consisting of the coordinate subspaces.  

Let us introduce the following notations:

\begin{equation*}
\mathcal{B}=\mathcal{B}_1\sqcup\ldots\sqcup \mathcal{B}_n,
\end{equation*}
where
\begin{equation*}
\mathcal{B}_k:=\{e_{d_{k-1}+1},\ldots, e_{d_k}\}.
\end{equation*}
(Here, as before $d_l=\sum_{a\le l} \mu_a$ are the partial sums.)

We also set
\begin{equation*}
V_k:=\spn(\mathcal{B}_k),
\end{equation*}
so that $V=\bigoplus V_k.$

Let $p_\mu\in \Fmu$ be the $T$-fixed point given by

\begin{equation*}
p_{\mu}:=\{V_1\subset V_1\oplus V_2\subset\ldots\subset V_1\oplus\ldots\oplus V_n=V\}.
\end{equation*}

%\begin{definition}
%Let us use the following notation: 
%$$
%V_k:=\left<e_1,\ldots,e_k\right>
%$$
%for any $1\le k\le d.$ We will call the coordinate flag
%$$
%x:=\{V_{d_1} \subset V_{d_2} \subset\ldots\subset V\}\in \Fmu
%$$ 
%the \textit{standard $T$-fixed point in $\Fmu$}.
%\end{definition} 
The stabilizer of the fixed point $p_\mu$ is the parabolic subgroup of $G$ consisting of block upper triangular matrices with blocks of sizes $\mu_1,\ldots,\mu_n,$ so the $T$-character of the fiber $E_x\subset E$ belongs to $\mathbb C[x_1^{\pm 1},\ldots,x_d^{\pm 1}]^{S_{\bf\mu}}.$ 
%\end{proof}

%Note that the $T$-character determines the representation of the parabolic subgroup only up to its class in the Grothendieck group. One also gets

%\begin{equation*}
%K^G(\Fmu)\simeq \mathbb C[x_1^{\pm 1},\ldots,x_d^{\pm 1}]^{S_{\bf\mu}}\simeq\bigotimes_i \Lambda^{\pm}_{\mu_i},
%\end{equation*}
%where $\Lambda^{\pm}_t$ is the ring of symmetric Laurent polynomials in $t$ variables. 

Note that one could have chosen a different torus fixed point. The $T$-characters of the fibers of equivariant vector bundles are then related by the following Lemma:

\begin{lemma}\label{lemma: characters at fixed points}
Let $p:E\to \Fmu$ be a $G$-equivariant vector bundle. Let $x,x'\in \Fmu^T$ be two $T$-fixed points. It follows that there exists a permutation $\omega\in S_d$ such that $\omega(x)=x'$ (here $S_d$ acts on $V$ by permuting the basic vectors $\{e_1,\ldots,e_d\}$). Let $f$ and $f'$ be the $T$-characters of the fibers $E_x$ and $E_{x'}$ correspondingly. Then one has
$$
\omega(f)=f',
$$
where $\omega$ acts by permuting the variables in the Laurent polynomial $f.$
\end{lemma}

\begin{proof}
The operator $\omega\in S_d\subset G$ acts on the vector bundle $E\to \Fmu$ and sends the fiber $E_x$ to the fiber $E_{x'}.$ It follows that the $T$-actions on $E_x$ and $E_{x'}$ are related by the automorphism $t\mapsto \omega^{-1}t\omega,$ which implies the formula.
\end{proof}

\subsection{Equivariant K-theory of the variety of pairs of partial flags.}\label{subsection: K-theory of FxF}

Let $M=\{m_{ij}\}$ be an $n\times k$ matrix with non-negative integer entries, and let $\mathbb{O}_M\subset \Fnd\times \mathcal{F}^d_k$ be the corresponding orbit. Let $\mathcal{B}=\{e_1,\ldots,e_d\}=\bigsqcup \mathcal{B}_{ij}$ be a decomposition of the basis $\mathcal{B}$ of $V,$ such that $\sharp\mathcal{B}_{ij}=m_{ij},$ and let $V_{ij}:=\spn\mathcal{B}_{ij}.$ 
Let also $\{1,\ldots,d\}=\bigsqcup \mathcal{I}_{ij}$ be the corresponding decomposition of the index set, so that for each $i$ and $j$ one has $\mathcal{B}_{ij}=\{e_a|a\in I_{ij}\}.$ Let $p\in \mathbb{O}_M$ be the pair of flags $p:=(\textbf{U},\textbf{W}),$ where 

\begin{equation*}
U_i:=\spn \bigoplus_{\substack{1\le a\le i\\ 1\le b\le k}} V_{ab}\quad \text{and} \quad W_i:=\spn \bigoplus_{\substack{1\le a\le n\\ 1\le b\le i}} V_{ab}.
\end{equation*}

The stabilizer $St_p\subset G$ is then the intersection of two parabolic subgroups $St_p=St_{\textbf{U}}\cap St_{\textbf{W}}.$ Also, every element of the stabilizer $St_p$ can be uniquely factored as $gh$ where

\begin{equation*}
g\in GL(p):=\prod_{\substack{1\le a\le n,\\ 1\le b\le k}} GL(V_{ab})\quad \text{and} \quad h\in N(p):=Id + \bigoplus_{\substack{1\le a\le n,\\ 1\le b\le k}} \text{Hom}(V_{ab},\bigoplus_{\substack{c\le a,\ d\le b \\ (c,d)\neq(a,b)}}V_{cd}).
\end{equation*}

Therefore, we get a semidirect product decomposition $St_p=GL(p)\ltimes N(p),$ where $GL(p)$ is semisimple and $N(p)$ is the unipotent radical.

Let $S_p\simeq \prod_{ij} S_{m_{ij}}\subset S_d$  (recall that $m_{ij}=\sharp\mathcal{B}_{ij}=\sharp\mathcal{I}_{ij}$) be the subgroup of permutations preserving the subdivision $\{1,\ldots,d\}=\bigsqcup I_{ij.}$ Similar to before, one gets

%\begin{lemma}\label{lemma: K^G(\mathbb{O}_M)}
%There is an isomorphism
\begin{equation*}\label{equation: K^G(O_M)}
K^G(\mathbb{O}_M)\simeq K^{GL(p)}(pt)\simeq\mathbb{C}[x_1^{\pm 1},\ldots,x_d^{\pm 1}]^{S_p}\simeq\bigotimes_{i,j} \Lambda^{\pm}_{m_{ij}}.
\end{equation*}

The isomorphism can be exhibited in the following way. Given an equivariant vector bundle $W$ over $\mathbb{O}_M$ the partially symmetric polynomial corresponding to the class $[W]\in K^G(\mathbb{O}_M)$ equals to the $T$-character of the fiber $W_p$ over the point $p\in \mathbb{O}_M.$

%where $S_M=S_{m_{11}}\times S_{m_{12}}\times\ldots\times S_{m_{nn}}\subset S_d$ is the subgroup of the symmetric group preserving the subdivision $\{1,\ldots,d\}=\{1,m_{11}\},\{m_{11}+1,\ldots, m_{11}+m_{12}\},\ldots, \{d-m_{nn}+1,\ldots, d\}.$ 
%\end{lemma}

Note that different choices of the subdivision $\mathcal{B}=\bigsqcup \mathcal{B}_{ij}$ satisfying the condition $\sharp\mathcal{B}_{ij}=m_{ij}$ lead to different fixed points in $\mathbb{O}_M$ and different (but conjugated) subgroups in $S_d.$ The isomorphisms from equation  \ref{equation: K^G(O_M)} for different fixed points are related in a way similar to Lemma \ref{lemma: characters at fixed points}:

\begin{lemma}\label{lemma: characters at fixed points (double)}
Let $p:E\to \mathbb{O}_M$ be a $G$-equivariant vector bundle over an orbit $\mathbb{O}_M\subset \Fnd\times \Fnd$. Let $x,x'\in \mathbb{O}_M^T$ be two $T$-fixed points. It follows that there exists a permutation $\omega\in S_d$ such that $\omega(x)=x'$ (here $S_d$ acts on $V$ by permuting the basic vectors $\{e_1,\ldots,e_d\}$). Let $f$ and $f'$ be the $T$-characters of the fibers $E_x$ and $E_{x'}$ correspondingly. Then one has
$$
\omega(f)=f',
$$
where $\omega$ acts by permuting the variable in the Laurent polynomial $f.$
\end{lemma}

This lemma is proved in the same way as Lemma \ref{lemma: characters at fixed points}.

\subsection{Push-forward and pull-back for orbits}\label{subsection: push-forward pull-back}

Let $\mu=(\mu_1,\ldots,\mu_{n+1})$ and $\nu=(\nu_1,\ldots,\nu_k)$ be two composition of $d,$ and let $\mu'=(\mu_1,\ldots,\mu_l+\mu_{l+1},\ldots,\mu_{n+1}).$ Let $M$ be an $(n+1)\times k$ matrix with non negative integer entries such that $\textbf{R}^M=\mu$ and $\textbf{C}^M=\nu,$ and let the $n\times k$ matrix $M'$ be obtained from $M$ by summing up the $l$th and $(l+1)$th rows, so that $\textbf{R}^{M'}=\mu'$ and $\textbf{C}^{M'}=\nu.$ Let $\mathbb{O}_M\subset \Fmu\times \mathcal{F}_\nu$ and $\mathbb{O}_{M'}\subset \mathcal{F}_{\mu'}\times \mathcal{F}_\nu$ be the corresponding orbits. Then one has a natural equivariant projection $\pi:\mathbb{O}_M\to \mathbb{O}_{M'}$ given by forgetting the $k$th step of the first flag. This gives rise to the pull-back and push-forward operators:

$$
\pi^*:K^G(\mathbb{O}_{M'})\to K^G(\mathbb{O}_M),
$$
and
$$
\pi_*:K^G(\mathbb{O}_M)\to K^G(\mathbb{O}_{M'}).
$$

Let $p\in \mathbb{O}_M$ and $p'\in \mathbb{O}_{M'}$ be torus fixed points such that $\pi(p)=p',$ and let $S_p\subset S_d$ and $S_{p'}\subset S_d$ be the corresponding subgroups. Note that $S_p\subset S_{p'}.$ One can write explicit formulas for the operators $\pi^*$ and $\pi_*$ in terms of the isomorphisms $K^G(\mathbb{O}_M)\simeq\mathbb{C}[x_1^{\pm 1},\ldots,x_d^{\pm 1}]^{S_p}$ and $K^G(\mathbb{O}_{M'})\simeq\mathbb{C}[x_1^{\pm 1},\ldots,x_d^{\pm 1}]^{S_{p'}}:$

\begin{lemma}\label{Lemma: pull-back push-forward}
Let $f({\bf x})\in \mathbb{C}[x_1^{\pm 1},\ldots,x_d^{\pm 1}]^{S_{p'}}\subset \mathbb{C}[x_1^{\pm 1},\ldots,x_d^{\pm 1}]^{S_p},$ and $g({\bf x})\in \mathbb{C}[x_1^{\pm 1},\ldots,x_d^{\pm 1}]^{S_p}.$ One has the following formulas:
$$
\pi^* f({\bf x})=f({\bf x}),
$$
and
$$
\pi_* g({\bf x})=Sym_{S_p}^{S_{p'}}\frac{g({\bf x})}{\prod\limits_{1\le a\le k}\prod\limits_{\substack{i\in I_{l,a},\\ j\in I_{l+1,a}}} (1-\frac{x_i}{x_j})},
$$
where $Sym_{S_p}^{S_{p'}}(*)$ is the summation of $\sigma(*)$ over a choice of representatives $\sigma$ of the left cosets of $S_p$ in $S_{p'}.$ 
\end{lemma}

\begin{remark}
Note that the second formula above does not depend on the choice of representatives, because both $g({\bf x})$ and $\prod\limits_{1\le a\le k}\prod\limits_{\substack{i\in I_{l,a},\\ j\in I_{l+1,a}}} (1-\frac{x_i}{x_j})$ belong to $\mathbb{C}[x_1^{\pm 1},\ldots,x_d^{\pm 1}]^{S_p}.$
\end{remark}

\begin{proof}
%Let $X\in \Fmu$ and $Y\in \mathcal{F}_{\mu'}$ be the standard flags, i.e.
%$$
%X=\left\{\{0\}=V_0\subset V_{d_1}\subset\ldots\subset V_{d_{k-1}}\subset V_{d_k}\subset V_{d_{k+1}}\subset\ldots\subset V_n=V\right\}
%$$
%and 
%$$
%Y=\left\{\{0\}=V_0\subset V_{d_1}\subset\ldots\subset V_{d_{k-1}}\subset V_{d_{k+1}}\subset\ldots\subset V_n=V\right\},
%$$
%where $V_m=\left<e_1,\ldots,e_m\right>$ as before. 

Since $\pi(p)=p',$ for any vector bundle $E\to \mathbb{O}_{M'}$ the fiber of $\pi^* E$ at $p$ is by definition isomorphic to the fiber of $E$ at $p',$ so the first part of the Lemma follows immediately.

We will use the  Lefschetz formula (see \ref{equation: Lefschetz formula}) to compute the pushforward $\pi_*.$ There are finitely many torus fixed points in both $\mathbb{O}_M$ and $\mathbb{O}_{M'}:$ those are precisely the flags consisting of coordinate subspaces. It follows that for a generic point $t\in T\subset G$ the $t$-fixed points coincide with the full torus fixed points. We need to compute the $T$-character of the fiber of $\pi_* f$ at $p'.$ Let $\widetilde{\pi}:=\pi|_{\pi^{-1}(p')}:\pi^{-1}(p')\to {p'}.$ We will use the Lefschetz formula to compute $\widetilde{\pi}_* (f|_{\pi^{-1}(p')}),$ which is exactly what we need.

The symmetric group $S_d$ naturally acts both on the set of torus fixed points in $\mathbb{O}_M$ and the set of torus fixed points in $\mathbb{O}_{M'}.$ Moreover, both actions are transitive, they commute with the projection $\pi,$ and the stabilizers of $p\in \mathbb{O}_M$ and $p'\in \mathbb{O}_{M'}$ are $S_p$ and $S_{p'}$ correspondingly. Therefore, $S_{p'}$ acts transitively on the torus fixed points in the fiber $\pi^{-1}(p'),$ and the stabilizer of $p\in\pi^{-1}(p')$ is $S_p\subset S_{p'}.$ We conclude that for any collection of representatives $\{\omega_1,\ldots,\omega_N\}$ of the left cosets of $S_p$ in $S_{p'}$ the torus fixed points in the fiber $\pi^{-1}(p')$ are exactly $\{\omega_1(p),\ldots,\omega_N(p)\},$ and every point is listed exactly once.

The last ingredient we need is the class of the exterior algebra of the conormal bundle to the fixed point $p$ in $K^T(\pi^{-1}(p')).$ Since the fixed points are isolated, this is simply the $T$-character of the exterior algebra of the cotangent space to the fiber $\pi^{-1}(p')$ at $p.$ Fixing the maximal torus $T$ also fixes the $T$-invariant direct sum decomposition $V=\bigoplus V_{ij},$ so that the pair of flags $p=(\textbf{U},\textbf{W})$ is given by 

\begin{equation*}
U_i:=\spn \bigoplus_{\substack{1\le a\le i\\ 1\le b\le k}} V_{ab}\quad \text{and} \quad W_i:=\spn \bigoplus_{\substack{1\le a\le n+1\\ 1\le b\le i}} V_{ab},
\end{equation*}
as in Section \ref{subsection: K-theory of FxF}. Let also $M=\{m_{ij}\}.$ Then $\pi^{-1}(p')$ is isomorphic to the product of Grassmanians:

\begin{equation*}
\pi^{-1}(p')\simeq \prod_{a=1}^k Gr(m_{la},V_{l,a}\oplus V_{l+1,a}),
\end{equation*}
where $Gr(m_{la},V_{l,a}\oplus V_{l+1,a})$ is the Grassmanian of $m_{la}$-dimensional subspaces in $V_{l,a}\oplus V_{l+1,a}.$ The fixed point $p$ corresponds to $(V_{l,1},\ldots, V_{l,k}).$ For each $a$ the eigenvalues of the cotangent space of the Grassmanian  $Gr(m_{la},V_{l,a}\oplus V_{l+1,a})$ at the point $V_{l,a}$ are $\frac{x_i}{x_j}$ for $i\in I_{l,a}$ and $j\in I_{l+1,a},$ each with multiplicity $1.$ Therefore, for the exterior algebra one gets 
$$
\prod\limits_{1\le a\le k}\prod\limits_{\substack{i\in I_{l,a},\\ j\in I_{l+1,a}}} (1-\frac{x_i}{x_j})
$$
Using Lemma \ref{lemma: characters at fixed points} and the Lefschetz formula we get the desired result.
\end{proof}

\begin{remark}
Note that if $\nu=(d)$ one gets $\mathbb{O}_{M'}=\mathcal{F}_{\mu'}$ and $\mathbb{O}_M=\Fmu.$ Therefore, Lemma \ref{Lemma: pull-back push-forward} is also applicable to the natural projections of the varieties of partial flags. For $\pi:\mathcal{F}_{\mu'}\to \Fmu$ one gets
$$
\pi^* f({\bf x})=f({\bf x}),
$$
and
\begin{equation}\label{equation: pi_* for partial flags}
\pi_* g({\bf x})=Sym_{S_\mu}^{S_{\mu'}}\frac{g({\bf x})}{\prod\limits_{i\in I_l,\ j\in I_{l+1}} (1-\frac{x_i}{x_j})},
\end{equation}
where $\{1,\ldots,d\}=\bigsqcup\limits_{a=1}^{n+1} I_a$ is the corresponding subdivision. We illustrate this formula in Example \ref{example: pi_*}.
\end{remark}

Note that since $\pi:\mathbb{O}_M\to \mathbb{O}_{M'}$ is a smooth fibration with connected fibers, one has $\pi_* 1=1.$ One can obtain this from the formula in Lemma \ref{Lemma: pull-back push-forward} although it requires some manipulation with determinants. One also gets the projection formula:

\begin{equation*}
\pi_*(\pi^*(f)g)=f\pi_*(g).
\end{equation*}

This can also be obtained from the formula in Lemma \ref{Lemma: pull-back push-forward} directly, since if $f\in\mathbb{C}[x_1^{\pm 1},\ldots,x_d^{\pm 1}]^{S_{p'}}$ then one has

\begin{equation*}
Sym_{S_p}^{S_{p'}} fg=f Sym_{S_p}^{S_{p'}} g.
\end{equation*}

In particular, one obtains the following Lemma:

\begin{lemma}\label{lemma: push-forward surjectivity}
Let $\pi: \mathbb{O}_M\to \mathbb{O}_{M'}$ be a projection map as above. Then the push-forward operator $\pi_*:K^G(\mathbb{O}_M)\to K^G(\mathbb{O}_{M'})$ is surjective.
\end{lemma}

\begin{proof}
Indeed, for any $f\in K^G(\mathbb{O}_{M'})$ one has
\begin{equation*}
f=f\cdot 1=f\cdot \pi_*(1)=\pi_*(\pi^*(f)\cdot 1).
\end{equation*}
\end{proof}

%Clearly, the closure $\overline{\mathbb{O}_M}$ is a union of orbits. This determines a partial order on non-negative integer martices: 

%\begin{corollary}
%A $G$ orbit in $\Fnd\times \Fnd$ is closed if and only if it can be obtained from a connected component $\Fmu\subset \mathcal{F}_{2n-1}$ in the way described above.    
%\end{corollary}

%{\color{red} There is a combinatorial description of this partial order, but the proofs get a bit technical and it is not strictly required for our purposes.}

%We get the following combinatorial condition for the cover relations in this partial order.

%\begin{theorem}
%Let $M=\{m_{ij}\}$ and $N=\{n_{ij}\}$ be two non-negative integer matrices. Then $N\lessdot M$ if and only if the following condition is satisfied. There exist two distinct integers $1\le k<l\le n$ such that $n_{ij}$ 
%\begin{equation*}
%\begin{aligned}
%n_{kk}&=m_{kk}+1\\
%n_{ll}&=m_{ll}+1\\
%n_{kl}&=m_{kl}-1\\
%n_{lk}&=m_{lk}-1
%\end{aligned}
%\end{equation*}
%\end{theorem}

We will need to compute the push-forwards in some special cases. Let $\mu=(\mu_1,\ldots, \mu_n)$ be a composition of $d$. Let $\mu'=(\mu_1,\ldots,\mu_{k-1},1,\ldots,1,\mu_{k+1},\ldots,\mu_n)$ be its refinement, where the part $\mu_k$ is split into $\mu_k$ ones. Let $\mathcal{L}_1,\ldots,\mathcal{L}_{\mu_k}$ be the tautological line bundles over $\mathcal{F}_{\mu'}$ corresponding to the one-dimensional steps of the flag, left to right. Let $\pi:\mathcal{F}_{\mu'}\to \Fmu$ be the natural projection, and let $\lambda=(\lambda_1,\ldots,\lambda_{\mu_k})$ be a partition of length $\le\mu_k$ (i.e. $\lambda_1\ge\ldots\ge \lambda_{\mu_k}\ge 0$). According to Borel-Weil-Bott Theorem, one gets
\begin{equation}\label{formula: push-forward Schur functor}
\pi_*([\mathcal{L}_1]^{\lambda_{\mu_k}}[\mathcal{L}_2]^{\lambda_{\mu_k-1}}\ldots[\mathcal{L}_{\mu_k}]^{\lambda_1})=[\mathbb{S}_\lambda(\mathcal{T}_k)],
\end{equation}
where $\mathcal{T}_k$ is the tautological vector bundle over $\Fmu$ corresponding to the $k$th step of the flag, and $\mathbb{S}_\lambda$ is the Schur functor. In terms of the partially symmetric functions one gets
\begin{equation}\label{formula: push-forward Schur}
\pi_*(y_1^{\lambda_{\mu_k}}y_2^{\lambda_{\mu_k-1}}\ldots y_{\mu_k}^{\lambda_1})=s_\lambda(y_1,\ldots,y_{\mu_k}),
\end{equation}
where $y_1,\ldots,y_{\mu_k}$ are the torus characters corresponding to the $k$th step of the flag, and $s_\lambda$ is the Schur function. One can also obtain this formula directly from Lemma \ref{Lemma: pull-back push-forward}. 

Formula (\ref{formula: push-forward Schur}) also holds in larger generality. Let $\lambda=(\lambda_1,\ldots,\lambda_{\mu_k})\in \mathbb{Z}^{\mu_k}$ be an arbitrary integer vector. Let $s_\lambda$ be the \textit{generalized Schur function} defined as follows:

\begin{definition}
Let $\alpha=(\alpha_1,\ldots,\alpha_m)\in\mathbb{Z}^m$ be an integer vector. Let $\det_\alpha(y_1,\ldots,y_m)$ denote the determinant of the matrix $M_\alpha=\{m_{ij}\}=\{y_i^{\alpha_j}\}.$ Let also $\rho^m=(m-1,m-2,\ldots,0).$ In particular, $\det_{\rho^m}(y_1,\ldots,y_n)$ is the Vandermonde determinant.
\end{definition}

\begin{definition}
The Schur function $s_\lambda(y_1,\ldots,y_m)$ is defined as follows:
$$
s_\lambda=\frac{\det_{\lambda+\rho^m}(y_1,\ldots,y_m)}{\det_{\rho^m}(y_1,\ldots,y_m)}.
$$
\end{definition}

The following relation follows immediately from the above definition:

\begin{equation*}
s_\lambda(y_1,\ldots,y_m)=y_1\ldots y_ms_{\lambda -(1,\ldots,1)},
\end{equation*}
and
\begin{equation}\label{formula: general schur}
s_{\lambda_1,\ldots,\lambda_i,\lambda_{i+1},\ldots,\lambda_m}=-s_{\lambda_1,\ldots,\lambda_{i+1}-1,\lambda_i+1,\ldots,\lambda_m}.
\end{equation}
In particular, one gets that for every $\lambda\in (\mathbb{Z})^m$ the Schur function $s_\lambda$ is either zero, or equals to $\pm \frac{s_{\nu}(y_1,\ldots, y_m)}{(y_1\ldots y_m)^N}$ for a partition $\nu$ (i.e. $\nu_1\ge\ldots\ge\nu_m\ge 0$) and a positive integer $N$. One can use these formulas to define the \textit{generalized Schur functors} $\mathbb{S}_\lambda$ for any $\lambda\in \mathbb{Z}^m.$ One immediately gets 

\begin{lemma}
Formulas (\ref{formula: push-forward Schur functor}) and (\ref{formula: push-forward Schur}) hold for any $\lambda\in \mathbb{Z}^m.$
\end{lemma}

Formula (\ref{formula: general schur}) also implies that if for some $i$ and $k$ one has $\lambda_{i+k}=\lambda_i+k$ then $s_\lambda=0.$ In particular, for $0<k<m$ one gets

\begin{equation*}
s_{(-k)}(y_1,\ldots,y_m)=0.
\end{equation*}

In some cases it is also convenient to reduce generalized Schur functions to the form 
\begin{equation*}
\pm s_{\nu}(\frac{1}{y_1},\ldots, \frac{1}{y_m})(y_1\ldots y_m)^N,
\end{equation*}
\noindent where $\nu$ is a partition. Note that

\begin{align*}
\det\nolimits_{\rho^m}(y_1,\ldots,y_m)&=\det\nolimits_{(0,1,\ldots,m-1)}(\frac{1}{y_1},\ldots,\frac{1}{y_m})(y_1\ldots y_m)^{m-1}\\
&=(-1)^{\frac{m(m-1)}{2}}\det\nolimits_{\rho^m}(\frac{1}{y_1},\ldots,\frac{1}{y_m})(y_1\ldots y_m)^{m-1},
\end{align*}
%\begin{equation*}
%\det\nolimits_{\rho^m}(\frac{1}{y_1},\ldots,\frac{1}{y_m})=\frac{\det\nolimits_{(0,1,\ldots,m-1)}(y_1,\ldots,y_m)}{(y_1\ldots y_m)^{m-1}}=\frac{(-1)^{\frac{m(m-1)}{2}}\det\nolimits_{\rho^m}(y_1,\ldots,y_m)}{(y_1\ldots y_m)^{m-1}},
%\end{equation*}
and for $k\ge m$
\begin{align*}
\det\nolimits_{(-k,0,\ldots,0)+\rho^m}(y_1,\ldots,y_m)&=\det\nolimits_{(k-1,0,\ldots,m-2)}(\frac{1}{y_1},\ldots,\frac{1}{y_m})(y_1\ldots y_m)^{m-2}\\
&=(-1)^{\frac{(m-1)(m-2)}{2}}\det\nolimits_{(k-m,0,\ldots,0)+\rho^m}(\frac{1}{y_1},\ldots,\frac{1}{y_m})(y_1\ldots y_m)^{m-2}.
\end{align*}
Therefore, in this case one gets
\begin{align*}
s_{(-k)}(y_1,\ldots,y_m)&=\frac{(-1)^{m-1}s_{(k-m)}(\frac{1}{y_1},\ldots,\frac{1}{y_m})}{y_1\ldots y_m}\\
&=\frac{(-1)^{m-1}h_{k-m}(\frac{1}{y_1},\ldots,\frac{1}{y_m})}{y_1\ldots y_m},
\end{align*}
where $h_{k-m}$ is the complete homogeneous symmetric polynomial.

One can use these formulas to define generalized classes of Schur functors $\mathbb{S}_\lambda$ for a general $\lambda\in \mathbb{Z}^m.$ In particular, for $k\in\mathbb{Z}_{\ge 0}$ one gets

\begin{equation*}
\mathbb{S}_{(-k)}(V)=\left\{
\begin{array}{rl}
0,& 0<k<m,\\
(-1)^{m-1}[\Lambda_m(V^*)][Sym_{k-m}(V^*)],& k\ge m.
\end{array}
\right.
\end{equation*}

\begin{corollary}\label{corollary: pi_*^-}
Let $\lambda:=(\lambda_1,\ldots,\lambda_n)$ and $\lambda^-:=(\lambda_1,\ldots,\lambda_{i-1},1,\lambda_i-1,\lambda_{i+1},\ldots,\lambda_n)$ be compositions, and let $\pi^-:\mathcal{F}_{\lambda^-}\to \mathcal{F}_\lambda$ be the projection. Let $\mathcal{L}^-$ be the tautological line bundle over $\mathcal{F}_{\lambda^-}$ corresponding to the one dimensional step in the flag, and $\mathcal{T}_i$ be the tautological vector bundle over $\mathcal{F}_\lambda$ corresponding to the $i$th step. Then
\begin{equation*}
\pi_*^-([\mathcal{L}^-]^p)=
\left\{
\begin{array}{rl}
(-1)^{\lambda_i-1}[\Lambda_{\lambda_i}(\mathcal{T}_i)][Sym_{p-\lambda_i}(\mathcal{T}_i)],& p\ge \lambda_i,\\
0,& 0<p<\lambda_i,\\
\lbrack Sym_{-p}(\mathcal{T}_i^*)],& p\le 0.\\
\end{array}
\right.
\end{equation*}
In terms of the partially symmetric functions one gets 
\begin{equation*}
\pi_*^-(y_1^p)=
\left\{
\begin{array}{rl}
(-1)^{\lambda_i-1}y_1\ldots y_{\lambda_i} h_{p-\lambda_i}(y_1,\ldots,y_{\lambda_i}),& p\ge \lambda_i,\\
0, &0<p<\lambda_i,\\
h_{-p}(\frac{1}{y_1},\ldots,\frac{1}{y_{\lambda_i}}),& p\le 0.
\end{array}
\right.
\end{equation*}
\end{corollary}

\begin{proof}
Consider the composition $\lambda'=(\lambda_1,\ldots,\lambda_{i-1},1,\ldots,1,\lambda_{i+1},\ldots,\lambda_n)$ and the corresponding variety of partial flags $\mathcal{F}_{\lambda'}$ with the projection map $pr:\mathcal{F}_{\lambda'}\to \mathcal{F}_{\lambda^-}.$ Note that $pr^* (\mathcal{L}^-)=\mathcal{L}_1$ and $pr_*(\mathcal{L}^1)=\mathcal{L}^-$ where $\mathcal{L}_1$ is the line bundle on $\mathcal{F}_{\lambda'}$ corresponding to the first one-dimensional step. Therefore, using projection formula one gets
\begin{equation*}
\pi_*^-([\mathcal{L}^-]^p)=(\pi^-\circ pr)_*([\mathcal{L}_1]^p)=[\mathbb{S}_{(0,\ldots,0,p)}(\mathcal{T}_i)] 
\end{equation*}
(there are $\lambda_i-1$ zeros). Further, we get
\begin{equation*}
[\mathbb{S}_{(0,\ldots,0,p)}(\mathcal{T}_i)]=(-1)^{\lambda_i-1}[\mathbb{S}_{(p-\lambda_i+1,1,\ldots,1)}(\mathcal{T}_i)]=(-1)^{\lambda_i-1}[\Lambda_{\lambda_i}(\mathcal{T}_i)][\mathbb{S}_{(p-\lambda_i)}(\mathcal{T}_i)],
\end{equation*}
after which we use the calculations from the above.
\end{proof}

\begin{corollary}\label{corollary: pi_*^+}
Let also $\lambda^+:=(\lambda_1,\ldots,\lambda_{i-1},\lambda_i-1,1,\lambda_{i+1},\ldots,\lambda_n),$ let $\pi^+:\mathcal{F}_{\lambda^+}\to \mathcal{F}_\lambda$ be the projection, and let $\mathcal{L}^+$ be the tautological line bundle over $\mathcal{F}_{\lambda^+}$ corresponding to the one dimensional step in the flag. Then
\begin{equation*}
\pi_*^+([\mathcal{L}^+]^p)=
\left\{
\begin{array}{rl}
[Sym_p(\mathcal{T}_i)],& p\ge 0,\\
0,& -\lambda_i<p<0,\\
(-1)^{\lambda_i-1}[\Lambda_{\lambda_i}(\mathcal{T}_i^*)][Sym_{-p-\lambda_i}(\mathcal{T}_i^*)],& p\le -\lambda_i.
\end{array}
\right.
\end{equation*}
In terms of the partially symmetric functions one gets 
\begin{equation*}
\pi_*^+(y_{\lambda_i}^p)=
\left\{
\begin{array}{rl}
h_p(y_1,\ldots,y_{\lambda_i}),& p\ge 0,\\
0, &-\lambda_i<p<0,\\
\frac{(-1)^{\lambda_i-1}h_{-p-\lambda_i}(\frac{1}{y_1},\ldots,\frac{1}{y_{\lambda_i}})}{y_1y_2\ldots y_{\lambda_i}},& p\le -\lambda_i.
\end{array}
\right.
\end{equation*}
\end{corollary}

\begin{proof}
Similar to the above we get 
\begin{equation*}
\pi_*^+([\mathcal{L}^+]^p)=[\mathbb{S}_{(p)}(\mathcal{T}_i)].
\end{equation*}
Now the corollary follows immediately from the computations above.
\end{proof}

Finally, the classes derived from the tautological bundles not affected by the projection behave in the natural way:

\begin{lemma}\label{lemma: push-forward for other steps}
As before, let $\mu=(\mu_1,\ldots,\mu_{n+1}),$  $\mu'=(\mu_1,\ldots,\mu'_l,\mu''_l,\ldots,\mu_{n+1}),$ where $\mu'_l+\mu''_l=\mu_l,$ and let $\pi: \Fmu\to \mathcal{F}_{\mu'}$ be the natural projection. Let $\mathcal{T}_k$ and $\mathcal{T}'_k$ be the tautological bundles corresponding to the step $\mu_k,\ k\neq l,$ over $\Fmu$ and $\mathcal{F}_{\mu'}$ respectively. Then for any $\lambda\in\mathbb{Z}^{\mu_k}$ and any $f\in K^G(\Fmu)$ one has
\begin{equation*}
\pi_*(\mathbb{S}_\lambda(\mathcal{T}_k)\cdot f)=\mathbb{S}_\lambda(\mathcal{T}'_k)\cdot \pi_*(f).
\end{equation*}
\end{lemma}

\begin{proof}
One gets $\pi^*(\mathbb{S}_\lambda(\mathcal{T}'_k))=\mathbb{S}_\lambda(\mathcal{T}_k)$ immediately from definitions. Then, applying the projection formula, one gets
\begin{equation*}
\pi_*(\mathbb{S}_\lambda(\mathcal{T}_k)\cdot f)=\pi_*\left(\pi^*(\mathbb{S}_\lambda(\mathcal{T}'_k))\cdot f\right)=\mathbb{S}_\lambda(\mathcal{T}'_k)\cdot \pi_*(f).
\end{equation*}
\end{proof}

\begin{remark}
Lemma \ref{lemma: push-forward for other steps} also easily follows from the formulas of Lemma \ref{Lemma: pull-back push-forward}. Furthermore, Corollaries \ref{corollary: pi_*^-} and \ref{corollary: pi_*^+} can also be derived directly from Lemma \ref{Lemma: pull-back push-forward}. We illustrate this in the Example \ref{example: pi_*} below.
\end{remark}

\begin{example}\label{example: pi_*}
    Consider compositions $\mu'=(2,1)$ and $\mu=(3).$ Then $\mathcal{F}_{\mu'}\simeq\mathbb{P}^2$ and $\mathcal{F}_{\mu}$ is a point. Consider the natural projection $\pi:\mathcal{F}_{\mu'}\to\mathcal{F}_\mu$ which is simply the projection to the point in this case. We get that $K^G(\mathcal{F}_{\mu'})$ is isomorphic to the ring of Laurent polynomials in three variables $(x_1, x_2, x_3),$ symmetric in $x_1$ and $x_2,$ and $K^G(\mathcal{F}_{\mu})$ is isomorphic to the ring of Laurent polynomial symmetric in all three variables.

    Consider $g(x_1,x_2,x_3)=x_3^k\in K^G(\mathcal{F}_{\mu'}),$ $k>0.$ Let us compute $\pi_* g$ using Formula \ref{equation: pi_* for partial flags}:
\begin{align*}
    \pi_* x_3^k =& Sym_{S_{(2,1)}}^{S_{(3)}}\frac{x_3^k}{(1-\frac{x_1}{x_3})(1-\frac{x_2}{x_3})}\\
    =& (1+\sigma_{1,3}+\sigma_{2,3})\frac{x_3^{k+2}}{(x_3-x_1)(x_3-x_2)},
\end{align*}
where $\sigma_{1,3}$ and $\sigma_{2,3}$ are permutations switching variables $x_1$ and $x_3,$ and $x_2$ and $x_3,$ respectively. Continuing the computation, we get:
\begin{align*}
    \pi_* x_3^k 
    =& \frac{x_3^{k+2}}{(x_3-x_1)(x_3-x_2)}+\frac{x_1^{k+2}}{(x_1-x_3)(x_1-x_2)}+\frac{x_2^{k+2}}{(x_2-x_1)(x_2-x_3)}\\
    =&\frac{x_3^{k+2}(x_2-x_1)+x_1^{k+2}(x_3-x_2)-x_2^{k+2}(x_3-x_1)}{(x_3-x_1)(x_3-x_2)(x_2-x_1)}\\
    =&\frac{\det_{(k+2,1,0)}(x_1,x_2,x_3)}{\det_{(2,1,0)}(x_1,x_2,x_3)}=s_k(x_1,x_2,x_3)=h_k(x_1,x_2,x_3).
\end{align*}
\end{example}

\subsection{Convolution algebra}\label{subsection: convolution algebra}

\begin{lemma}\label{lemma: short exact sequence}
Let $N\subset M\subset \Fnd\times \Fnd$ be two closed invariant subvarieties in the variety of pairs of partial flags. Then one gets the short exact sequence
\begin{equation*}
\begin{tikzcd}
0 \arrow[r] &K^G(N) \arrow[r,"i_*"] & K^G(M) \arrow[r, "j^*"] & K^G(M\setminus N)\arrow[r] & 0,
\end{tikzcd}
\end{equation*}
\noindent where $i:N \hookrightarrow M$ is the closed embedding, and $j:M\setminus N \hookrightarrow M$ is the open embedding.
\end{lemma}
%\begin{lemma}\label{lemma: short exact sequence}
%For any orbit $\mathbb{O}_A\subset \Fnd\times \Fnd$ one gets a short exact sequence
%\begin{equation*}
%\begin{tikzcd}
%0 \arrow[r] &K^G(\partial \mathbb{O}_A) \arrow[r,"i_*"] & K^G(\overline{\mathbb{O}_A}) \arrow[r, "j_*"] & K^G(\mathbb{O}_A)\arrow[r] & 0,
%\end{tikzcd}
%\end{equation*}
%\noindent where $\overline{\mathbb{O}_A}$ is the closure of the orbit $\mathbb{O}_A,$ $\partial \mathbb{O}_A=\overline{\mathbb{O}_A}\setminus \mathbb{O}_A$ is its boundary, $i:\partial \mathbb{O}_A \hookrightarrow \overline{\mathbb{O}_A}$ is the closed embedding, and $j:\mathbb{O}_A\hookrightarrow \overline{\mathbb{O}_A}$ is the open embedding.
%\end{lemma}

\begin{proof}
Set $U=M\setminus N.$ Then one gets the long exact sequence (see Equation \ref{equation: long exact sequence}):
\begin{equation*}
\begin{tikzcd}
\ldots\arrow[r] & K^G_0(N)\arrow[r,"i_*"] & K^G_0(M) \arrow[r,"j^*"] &K^G_0(U) \arrow[r] & 0,
\end{tikzcd}
\end{equation*}
Therefore, it suffices to prove that the map $i_*:K^G(N)\to K^G(M)$ is injective. For the sake of finding a contradiction, suppose that $i_*$ is not injective and let $\beta\neq 0$ be such that $i_* \beta =0.$ 

Let $\mathbb{O}_{B_1}\subset N$ be an open orbit in $N,$ and suppose that $\beta|_{\mathbb{O}_{B_1}}=0.$ Let $S_1:=N\setminus \mathbb{O}_{B_1}$ be the complement. Then it follows from the long exact sequence for $\mathbb{O}_{B_1}\subset N$ that there exists $\beta_1\in K^G(S_1)$ such that $i_{1*}\beta_1=\beta,$ where $i_1:S_1\hookrightarrow
 N$ is the closed embedding. Similarly, let $\mathbb{O}_{B_2}\subset S_1$ be an open orbit in $S_1$ such that $\beta_1|_{\mathbb{O}_{B_2}}=0,$ and let $S_2:=S_1\setminus U_1$ be the complement. Then it follows from the long exact sequence for $\mathbb{O}_{B_2}\subset S_1$ that there exists $\beta_2\in K^G(S_2)$ such that $i_{2*}\beta_2=\beta_1,$ where $i_2:S_2\hookrightarrow
 S_1$ is the closed embedding. Since there are only finitely many orbits in $\Fnd\times \Fnd,$ after repeating the above step several times one gets a closed invariant subvariety $i_S:S\hookrightarrow
 N,$ and a class $\alpha\in K^G(S)$ such that $i_{S*}\alpha=\beta$ and for any orbit $\mathbb{O}_C\subset S$ open in $S$ one has $\alpha|_{\mathbb{O}_C}\neq 0.$
 
Let $L\subset S$ be the union of all orbits open in $S,$ $L=\mathbb{O}_1\sqcup\ldots\sqcup \mathbb{O}_k,$ so that $\overline{L}=S.$ Let also $W:=\Fnd\times \Fnd\setminus \partial L.$ Then one gets the following Cartesian square:

\begin{equation*}
\begin{tikzcd} 
L \arrow[r, hook, "i_L"] \arrow[d, hook,"j_L"]& W \arrow[d, hook,"\hat{j}"] \\ 
S \arrow[r, hook, "\hat{i}"] &  \Fnd\times \Fnd 
\end{tikzcd}
\end{equation*}

Note that $\hat{i}: S\hookrightarrow \Fnd\times \Fnd$ factors as 

\begin{equation*}
\begin{tikzcd} 
S \arrow[r, hook, "i_S"] & N \arrow[r, hook, "i"] & M \arrow[r, hook, "\tilde{i}"] & \Fnd\times \Fnd. 
\end{tikzcd}
\end{equation*}

Hence one gets

\begin{equation*}
\hat{i}_* \alpha=\tilde{i}_*i_*i_{S*} \alpha=\tilde{i}_*i_* \beta=0.
\end{equation*}

Using base change, one gets

\begin{equation*}
0=\hat{j}^*\hat{i}_* \alpha = i_{L*}j_{L}^* \alpha=i_{L*} \alpha|_L.
\end{equation*}

Since $L=\mathbb{O}_1\sqcup\ldots\sqcup \mathbb{O}_k,$ one gets 

\begin{equation*}
K^G(L)=\bigoplus\limits_m K^G(\mathbb{O}_m),\ \ \alpha|_L=(\alpha_1,\ldots,\alpha_k),
\end{equation*}
and for every $1\le m\le k,$ $0\neq\alpha_m\in K^G(\mathbb{O}_m).$ Furthermore, for every $1\le m\le k,$ $\mathbb{O}_m\hookrightarrow W$ is a smooth closed embedding, therefore according to Equation \ref{equation: conormal bundle},

\begin{equation*}
i_L^*i_{L*}\alpha|_L=(\alpha_1\otimes\Lambda(T^*_{\mathbb{O}_1} W),\ldots,\alpha_k\otimes\Lambda(T^*_{\mathbb{O}_k} W)).
\end{equation*}

Since the torus fixed points are isolated it follows that $\Lambda(T^*_O W)\neq 0$ for every orbit $O.$ Also, it follows from equation \ref{equation: K^G(O_M)} that there are no zero divisors in $K^G(O).$ Therefore, 

\begin{equation*}
\alpha_m\otimes\Lambda(T^*_{\mathbb{O}_m} W)\neq 0
\end{equation*}
for every $1\le m\le k.$ In particular, $i_{L*}\alpha|_L\neq 0.$ Contradiction.

%\begin{equation*}
%S_k\subset S_{k-1}\subset\ldots\subset S_1\subset \partial \mathbb{O}_A
%\end{equation*}
%of closed embeddings of invariant subvarieties, and a sequence of elements

\end{proof}

\begin{corollary}\label{corollary: filtration}
For every closed invariant subvariety $M\subset \Fnd\times \Fnd$ one gets an injective morphism $K^G(M)\hookrightarrow K^G(\Fnd\times \Fnd),$ and the image is a subalgebra under tensor product.
\end{corollary}

We usually abuse notations by simply writing $K^G(M)\subset K^G(\Fnd\times \Fnd).$

\begin{corollary}\label{corollary: K(O) as quotient}
Furthermore, for every orbit $\mathbb{O}_A\subset \Fnd\times \Fnd$ one gets 
\begin{equation*}
K^G(\partial \mathbb{O}_A)\subset K^G(\overline{\mathbb{O}_A})\subset K^G(\Fnd\times \Fnd), 
\end{equation*}
where $\overline{\mathbb{O}_A}$ is the closure of the orbit $\mathbb{O}_A$ and $\partial \mathbb{O}_A=\overline{\mathbb{O}_A}\setminus \mathbb{O}_A$ is its boundary, and
\begin{equation*}
K^G(\mathbb{O}_A)=K^G(\overline{\mathbb{O}_A})/K^G(\partial \mathbb{O}_A).
\end{equation*}
\end{corollary}

It then follows that classes in $K^G(\Fnd\times \Fnd)$ have well defined supports:

\begin{corollary}\label{corollary: support}
For every class $\alpha\in K^G(\Fnd\times \Fnd)$ there exists a unique minimal (under inclusion) closed invariant subvariety $\supp(\alpha)\subset \Fnd\times \Fnd$ such that $\alpha\in K^G(\supp(\alpha))$ and for every orbit $O\subset \supp(\alpha)$ open in $\supp(\alpha)$ one has $\alpha|_O\neq 0.$
\end{corollary}

\begin{theorem}\label{theorem: Snd generated by diag and almost diag}
The convolution algebra $\Snd$ is generated by the classes supported on the orbits $\mathbb{O}_A$ where $A$ is either a diagonal matrix, or an almost diagonal matrix.
\end{theorem}

\begin{proof} Similar to Corollary \ref{corollary:matrix factorization}, the proof proceeds by induction in the diagonal norm $\dnorm{M}$ (see Definition \ref{definition: dnorm}).

%\begin{equation*}
%DN(M):=\sum\limits_{1\le i,j\le n} |i-j|m_{ij}.
%\end{equation*}

More precisely, we will use induction in $\max\{\dnorm{N}|\mathbb{O}_N\subset\supp(f)\}$ for $f\in K^G(\Fnd\times \Fnd).$ If $\max\{\dnorm{N}|\mathbb{O}_N\subset\supp(f)\}\le 1$ then $f$ is supported on a union of closed orbits corresponding to diagonal and almost diagonal matrices, therefore the base case of the induction is immediate.    

Note that whenever $N\lessdot M$ one has $\dnorm{N}\le \dnorm{M}$ (this can be checked by a direct computation). It follows that one can choose an orbit $\mathbb{O}_M\subset \supp(f)$ such that $\dnorm{M}=\max\{\dnorm{N}|\mathbb{O}_N\subset \supp(f)\}$ and $\mathbb{O}_M$ is open in $\supp(f).$ We will show that there exists a class $f'\in K^G(\Fnd\times \Fnd)$ such that $\supp(f')\subset\supp(f),$ $\supp(f-f')\subset\supp(f)\setminus \mathbb{O}_M,$ and $f'$ is generated by classes supported on orbits with lower diagonal norm. Then one can finish the proof by applying this argument several times and using the inductive assumption.

Applying Corollary \ref{corollary:strictly monotone factorization} to $M$ we obtain a factorization $M=D\circ N$ such that $D$ is almost diagonal, $\dnorm{N}<\dnorm{M},$ and for any $N'\prec N$ one has $D\circ N'\prec M.$ Since $\mathbb{O}_D$ is a closed orbit, Lemma \ref{lemma: triple intersection} implies that the triple intersection 
\begin{equation*}
\mathbb{O}:=\pi_{12}^{-1}(\mathbb{O}_D)\cap\pi_{13}^{-1}(\mathbb{O}_M)\cap\pi_{23}^{-1}(\mathbb{O}_N)\subset \Fnd\times \Fnd\times \Fnd
\end{equation*}
is a single orbit. According to Lemma \ref{lemma: orbits in preimage of closed orbits}, the orbit $\mathbb{O}\subset \Fnd\times \Fnd\times \Fnd$ is isomorphic to an orbit $\mathbb{O}_{L}\subset \mathcal{F}_{n+1}^d\times \Fnd,$ where $L$ is an $(n+1)\times n$ matrix. Let us first focus on the convolution map on the orbits:
\begin{equation*}
\tilde{\star}:K^G(\mathbb{O}_D)\otimes K^G(\mathbb{O}_N) \xrightarrow{\tilde{\pi}_{12}^*\otimes\tilde{\pi}_{23}^*} K^G(\mathbb{O})\xrightarrow{\tilde{\pi}_{13*}} K^G(\mathbb{O}_M),
\end{equation*}
where $\tilde{\pi}_{ij}=\pi_{ij}|_{\mathbb{O}}.$ We will prove that this multiplication map is surjective. The push-forward map $\tilde{\pi}_{13*}: K^G(\mathbb{O}_L)\simeq K^G(\mathbb{O})\to K^G(\mathbb{O}_M)$ is surjective by Lemma \ref{lemma: push-forward surjectivity}, therefore it will suffice to prove that $\tilde{\pi}_{12}^*\otimes\tilde{\pi}_{23}^*$ is surjective as well.

Suppose that $D$ is upper-triangular, i.e. $D=E(\mu,k)$ where $1\le k<n$ and $\mu=\textbf{R}^M-\textbf{e}_k$ (the lower-triangular case is done in a similar way). Here, according to Theorem \ref{theorem: matrix mult by almost diag}, Lemma \ref{lemma: non-unique division}, and Corollary \ref{corollary:strictly monotone factorization}, the choice of $k$ satisfies the following condition: there exists $l$ such that 
\begin{equation*}
m_{k,l+1}=\ldots=m_{k,n}=m_{k+1,l+1}=\ldots=m_{k+1,n}=0,
\end{equation*}
and $m_{kl}>0.$ Then the matrix $L$ is obtained from $M$ by reducing $m_{kl}$ by one and inserting the row $(0,\ldots,0,1,0,\ldots,0)$ between $k$th and $(k+1)$th rows (here $1$ is in the $l$th position): 
\begin{equation*}
L=\left(
\begin{array}{ccccccc}
m_{11}    & \ldots & m_{1,l-1}  & m_{1l}   & m_{1,l+1}& \ldots & m_{1n} \\
\vdots    &        & \vdots     & \vdots   & \vdots   &        & \vdots \\
m_{k1}    & \ldots & m_{k,l-1}  & m_{kl}-1 & 0        & \ldots & 0      \\
0         & \ldots & 0          &    1     & 0        & \ldots & 0      \\
m_{k+1,1} & \ldots & m_{k+1,l-1}& m_{k+1,l}& 0        & \ldots & 0      \\
\vdots    &        & \vdots     & \vdots   & \vdots   &        & \vdots   \\
m_{n1}    & \ldots & m_{n,l-1}  & m_{nl}   & m_{n,l+1}& \ldots & m_{nn} \\
\end{array}
\right),
\end{equation*} 
and the matrix $N$ is obtained from $L$ by adding this inserted row to the row below it: 
\begin{equation*}
N=\left(
\begin{array}{ccccccc}
m_{11}    & \ldots & m_{1,l-1}  & m_{1l}     & m_{1,l+1}& \ldots & m_{1n} \\
\vdots    &        & \vdots     & \vdots     & \vdots   &        & \vdots \\
m_{k1}    & \ldots & m_{k,l-1}  & m_{kl}-1   & 0        & \ldots & 0      \\
m_{k+1,1} & \ldots & m_{k+1,l-1}& m_{k+1,l}+1& 0        & \ldots & 0      \\
\vdots    &        & \vdots     & \vdots     & \vdots   &        & \vdots   \\
m_{n1}    & \ldots & m_{n,l-1}  & m_{nl}     & m_{n,l+1}& \ldots & m_{nn} \\
\end{array}
\right).
\end{equation*} 

%The push-forward map $\tilde{\pi}_{13*}:K^G(\mathbb{O}_L)\to K^G(\mathbb{O}_M)$ is surjective according to Lemma \ref{lemma: push-forward surjectivity}. 
Pick torus fixed points $p_M\in \mathbb{O}_M,\ p_D\in \mathbb{O}_D,$ and $p_L\in O\simeq \mathbb{O}_L$ so that $\tilde{\pi}_{12}(p_L)=p_D$ and $\tilde{\pi}_{23}(p_L)=p_D.$ According to equation \ref{equation: K^G(O_M)}, this produces isomorphisms 
\begin{align*}
&K^G(\mathbb{O}_N)\simeq \mathbb{C}[x_1^{\pm 1},\ldots,x_d^{\pm 1}]^{S_N}\simeq \left(\bigotimes_{(i,j)\notin\{(k,l),(k+1,l)\}} \Lambda^{\pm}_{m_{ij}}\right)\otimes \Lambda^{\pm}_{m_{kl}-1}\otimes \Lambda^{\pm}_{m_{k+1,l}+1},\\
&K^G(\mathbb{O}_D)\simeq \mathbb{C}[x_1^{\pm 1},\ldots,x_d^{\pm 1}]^{S_D}\simeq \left(\bigotimes_i \Lambda^{\pm}_{\mu_i}\right)\otimes \Lambda^{\pm}_1,\\
&K^G(\mathbb{O}_L)\simeq \mathbb{C}[x_1^{\pm 1},\ldots,x_d^{\pm 1}]^{S_L}\simeq \left(\bigotimes_{(i,j)\neq(k,l)} \Lambda^{\pm}_{m_{ij}}\right)\otimes \Lambda^{\pm}_{m_{kl}-1}\otimes \Lambda^{\pm}_1.
\end{align*}

The morphisms $\tilde{\pi}_{12}^*: K^G(\mathbb{O}_D)\to K^G(\mathbb{O}_L)$ and $\tilde{\pi}_{23}^*: K^G(\mathbb{O}_N)\to K^G(\mathbb{O}_L)$ are the natural embeddings, where the factors $\Lambda^{\pm}_{m_{ij}},$ $(i,j)\notin \{(k,l),(k+1,l)\},$  and $\Lambda^{\pm}_{m_{kl}-1}$ of $K^G(\mathbb{O}_N)$ are mapped isomorphically to the corresponding factors of $K^G(\mathbb{O}_L),$ the factor $\Lambda^{\pm}_{m_{k+1,l}+1}$ of $K^G(\mathbb{O}_N)$ is mapped injectively to $\Lambda^{\pm}_{m_{k+1,l}}\otimes \Lambda^{\pm}_1,$ the factors $\Lambda^{\pm}_{\mu_i},$ $i\neq k$  of $K^G(\mathbb{O}_D)$ are mapped injectively to $\bigotimes_j \Lambda^{\pm}_{m_{ij}}$ (recall that for $i\neq k$ one has $\mu_i=\sum_j m_{ij}$), the factor $\Lambda^{\pm}_{\mu_k}$ is mapped injectively to $\left(\bigotimes_{j\neq l} \Lambda^{\pm}_{m_{kj}}\right)\otimes \Lambda^{\pm}_{m_{kl}-1}$ (recall that $\mu_k=\sum_j m_{kj}-1$), and the factor $\Lambda^{\pm}_1$ of $K^G(\mathbb{O}_D)$ is mapped isomorphically to the corresponding factor of $K^G(\mathbb{O}_L).$

It remains to show that the factor $\Lambda^{\pm}_{m_{k+1,l}}$ is also generated by the images of $K^G(\mathbb{O}_N)$ and $K^G(\mathbb{O}_D),$ which follows from the elementary fact that the multiplication map

\begin{equation*}
\Lambda^{\pm}_{m_{k+1,l}+1}\otimes \Lambda^{\pm}_1 \to \Lambda^{\pm}_{m_{k+1,l}}\otimes \Lambda^{\pm}_1
\end{equation*}
is surjective (here the factor $\Lambda^{\pm}_{m_{k+1,l}+1}$ is mapped injectively to $\Lambda^{\pm}_{m_{k+1,l}}\otimes \Lambda^{\pm}_1,$ and $\Lambda^{\pm}_1$ is mapped isomorphically to $\Lambda^{\pm}_1$).

We have  $f\in K^G(\supp(f))\subset K^G(\Fnd\times\Fnd)$ and $\mathbb{O}_M\subset \supp(f),$ where $\mathbb{O}_M$ is open in $\supp(X).$ According to the above, the restriction $\tilde{f}:=f|_{\mathbb{O}_M}$ can be written as 
\begin{equation*}
\tilde{f}=\sum\limits_{i=1}^m g_i\tilde{\star} \tilde{h}_i,
\end{equation*}
where $m$ is an integer, $g_1,\ldots,g_m\in K^G(\mathbb{O}_D),$ and $\tilde{h}_1,\ldots,\tilde{h}_m\in K^G(\mathbb{O}_N).$ %According to Lemma \ref{lemma: short exact sequence}, one gets a short exact sequence

%\begin{equation*}
%\begin{tikzcd}
%0 \arrow[r] & K^G(\partial \mathbb{O}_N)\arrow[r] & K^G(\overline{\mathbb{O}_N})\arrow[r] & K^G(\mathbb{O}_N)\arrow[r] & 0.
%\end{tikzcd}
%\end{equation*}

According to Corollary \ref{corollary: K(O) as quotient}, the classes $\tilde{h}_1,\ldots,\tilde{h}_m$ can be extended to classes $h_1,\ldots h_m\in K^G(\overline{\mathbb{O}_N})\subset K^G(\Fnd\times \Fnd).$ Since $\mathbb{O}_D$ is closed, one also gets $g_1,\ldots,g_m\in K^G(\mathbb{O}_D)\subset K^G(\Fnd\times \Fnd).$ Consider
\begin{equation*}
f':=\sum\limits_{i=1}^m g_i\star h_i={\pi_{13}}_*\left(\sum\limits_{i=1}^m \pi_{12}^*(g_i)\otimes\pi_{23}^*(h_i)\right).
\end{equation*}

According to Corollary \ref{corollary: irreducibility with closure}, the intersection $(\mathbb{O}_D\times\Fnd)\cap(\Fnd\times\overline{\mathbb{O}_N})$ is irreducible. Furthermore, $\mathbb{O}$ is an open orbit in it, therefore $(\mathbb{O}_D\times\Fnd)\cap(\Fnd\times\overline{\mathbb{O}_N})=\overline{\mathbb{O}}.$ For every $i$  the tensor product $\pi_{12}^*(g_i)\otimes\pi_{23}^*(h_i)$ is supported on $\overline{\mathbb{O}},$ i.e. belongs to the image $j_*(K^G(\overline{\mathbb{O}}))$ where $j:\overline{\mathbb{O}}\hookrightarrow \Fnd\times\Fnd\times\Fnd$ is the embedding, and $f'$ is supported on $\pi_{13}(\overline{\mathbb{O}})=\overline{\mathbb{O}_M}.$ According to Lemma \ref{corollary: filtration}, we can then say that 
\begin{equation*}
f'\in K^G(\overline{\mathbb{O}_M})\subset K^G(\Fnd\times\Fnd).
\end{equation*}

Let $U:=\Fnd\times\Fnd\setminus\partial\mathbb{O}_M.$ The use of Corollary \ref{corollary:strictly monotone factorization} in the choice of the matrix $N$ in the decomposition $M=D\circ N$ and Lemma \ref{lemma: triple intersection} guarantee that 

$$
\pi_{12}^{-1}(\mathbb{O}_D)\cap\pi_{23}^{-1}(\overline{\mathbb{O}_N})\cap\pi_{13}^{-1}(U)=\pi_{12}^{-1}(\mathbb{O}_D)\cap\pi_{23}^{-1}(\mathbb{O}_N)\cap\pi_{13}^{-1}(U)=\mathbb{O}.
$$

Hence we can apply Corollary \ref{corollary: convolution with restriction} and conclude that 

$$
f'|_U=\tilde{f}=f|_{\mathbb{O}_M},
$$
where on the left we use restriction with supports (Equation \ref{equation: restriction with supports}).
Therefore, $$\supp(f-f')\subset \supp(f)\setminus \mathbb{O}_M,$$ which concludes our proof.
\end{proof}

\subsection{Local generators}\label{subsection: local generators}

Let $\mu=(\mu_1,\ldots,\mu_n)$ be a composition of $d-1$ of length $n$ and $1\le k<n$ be an integer. Consider the matrix $\textbf{E}(\mu,k)=Diag(\mu)+E_{k,k+1}$ and the corresponding (closed) orbit $\mathbb{O}_{\textbf{E}(\mu,k)}\subset \Fnd\times \Fnd.$ We have an isomorphism $\mathbb{O}_{\textbf{E}(\mu,k)}\simeq \mathcal{F}_{\tilde{\mu}},$ where $\tilde{\mu}=(\mu_1,\ldots,\mu_k,1,\mu_{k+1},\ldots,\mu_n).$ Let $\mathcal{T}\to \mathbb{O}_{\textbf{E}(\mu,k)}$ be the tautological line bundle corresponding to the codimension one step in the flag.

\begin{definition}
For any integer $p$ define the element $\mathfrak{e}_{\mu,k}(p)\in \Snd$ by
\begin{equation*}
\mathfrak{e}_{\mu,k}(p):=[\mathcal{T}]^p\in K^G(\mathbb{O}_{\textbf{E}(\mu,k)})\subset \Snd.
\end{equation*}
\end{definition} 

Similarly, for the matrix $\textbf{F}(\mu,k)=Diag(\mu)+E_{k+1,k}$ the corresponding (closed) orbit $\mathbb{O}_{\textbf{F}(\mu,k)}\subset \Fnd\times \Fnd$ is also isomorphic to $\mathcal{F}_{\tilde{\mu}}.$ Let $\mathcal{S}\to \mathbb{O}_{\textbf{F}(\mu,k)}$ be the tautological line bundle corresponding to the codimension one step in the flag.

\begin{definition}
For any integer $p$ define the element $\mathfrak{f}_{\mu,k}(p)\in \Snd$ by
\begin{equation*}
\mathfrak{f}_{\mu,k}(p):=[\mathcal{S}]^p\in K^G(\mathbb{O}_{\textbf{F}(\mu,k)})\subset \Snd.
\end{equation*}
\end{definition}

From considering supports, we immediately get the following Lemma:

\begin{lemma}\label{lemma: local gens prod=0 unless}
One has $\mathfrak{e}_{\mu,k}(p)\star\mathfrak{e}_{\nu,l}(q)=0$ unless
\begin{equation*}
\textbf{C}^{\textbf{E}(\mu,k)}=\textbf{R}^{\textbf{E}(\nu,l)} \Leftrightarrow \mu+\textbf{e}_{k+1}=\nu+\textbf{e}_l.
\end{equation*}
Similarly
\begin{equation*}
\mu+\textbf{e}_k\neq \nu+\textbf{e}_{l+1} \Rightarrow \mathfrak{f}_{\mu,k}(p)\star\mathfrak{f}_{\nu,l}(q)=0,
\end{equation*}
\begin{equation*}
\mu+\textbf{e}_k\neq \nu+\textbf{e}_l \Rightarrow \mathfrak{f}_{\mu,k}(p)\star\mathfrak{e}_{\nu,l}(q)=0,
\end{equation*}
and
\begin{equation*}
\mu+\textbf{e}_{k+1}\neq \nu+\textbf{e}_{l+1} \Rightarrow\mathfrak{e}_{\mu,k}(p)\star\mathfrak{f}_{\nu,l}(q)=0.
\end{equation*}
\end{lemma}

\begin{lemma}\label{lemma: local E F comm}
Let $1\le k,l<n$ be such that $|k-l|>1,$ and $\mu$ be a composition of $d-1$ of length $n$. Then
\begin{equation*}
\mathfrak{e}_{\mu+\textbf{e}_l,k}(p)\star\mathfrak{e}_{\mu+\textbf{e}_{k+1},l}(q)=\mathfrak{e}_{\mu+\textbf{e}_k,l}(q)\star\mathfrak{e}_{\mu+\textbf{e}_{l+1},k}(p),
\end{equation*}
and
\begin{equation*}
\mathfrak{f}_{\mu+\textbf{e}_{l+1},k}(p)\star\mathfrak{f}_{\mu+\textbf{e}_k,l}(q)=\mathfrak{f}_{\mu+\textbf{e}_{k+1},l}(q)\star\mathfrak{f}_{\mu+\textbf{e}_l,k}(p).
\end{equation*}
Also, for any $1\le k,l<n$ such that $k\neq l$ one has
\begin{equation*}
\mathfrak{e}_{\mu+\textbf{e}_{l+1},k}(p)\star\mathfrak{f}_{\mu+\textbf{e}_{k+1},l}(q)=\mathfrak{f}_{\mu+\textbf{e}_k,l}(q)\star\mathfrak{e}_{\mu+\textbf{e}_l,k}(p).
\end{equation*}
\end{lemma}

\begin{proof}
All three formulas are proved very similarly. Let us focus on the first one. On the level of supports one immediately gets
\begin{equation*}
\textbf{E}(\mu+\textbf{e}_l,k)\circ\textbf{E}(\mu+\textbf{e}_{k+1},l)=\textbf{E}(\mu+\textbf{e}_k,l)\circ\textbf{E}(\mu+\textbf{e}_{l+1},k)=Diag(\mu)+E_{k,k+1}+E_{l,l+1}.
\end{equation*}
Note that $\mathbb{O}_{Diag(\mu)+E_{k,k+1}+E_{l,l+1}}$ is a closed orbit, isomorphic to $\mathcal{F}_{(\mu_1,\ldots,\mu_k,1,\mu_{k+1},\ldots,\mu_l,1,\mu_{l+1},\ldots,\mu_n)}.$ Furthermore, the intersections $$\pi_{12}^{-1}(\mathbb{O}_{\textbf{E}(\mu+\textbf{e}_l,k})\cap \pi_{23}^{-1}(\mathbb{O}_{\textbf{E}(\mu+\textbf{e}_{k+1},l)})\subset \Fnd\times \Fnd\times \Fnd$$ and $$\pi_{12}^{-1}(\mathbb{O}_{\textbf{E}(\mu+\textbf{e}_k,l})\cap \pi_{23}^{-1}(\mathbb{O}_{\textbf{E}(\mu+\textbf{e}_{l+1},k)})\subset \Fnd\times \Fnd\times \Fnd,$$ are also isomorphic to $\mathcal{F}_{(\mu_1,\ldots,\mu_k,1,\mu_{k+1},\ldots,\mu_l,1,\mu_{l+1},\ldots,\mu_n)},$ and the restrictions of $\pi_{13}$ to these orbits are isomorphisms to $\mathbb{O}_{Diag(\mu)+E_{k,k+1}+E_{l,l+1}}$. Note that according to Lemma \ref{lemma: transversality}, both these intersections are transversal, and one can apply Corollary \ref{corollary: convolution with restriction} with $U=\Fnd\times\Fnd$. Finally, both the left hand side and the right hand side of the equation are equal to 
\begin{equation*}
z_k^pz_l^q\in K^G(\mathbb{O}_{Diag(\mu)+E_{k,k+1}+E_{l,l+1}})\subset K^G(\Fnd\times \Fnd),
\end{equation*}
where $z_k$ and $z_l$ are the classes of the tautological line bundles corresponding to parts of size $1$ in the composition $(\mu_1,\ldots,\mu_k,1,\mu_{k+1},\ldots,\mu_l,1,\mu_{l+1},\ldots,\mu_n)$ between $\mu_k$ and $\mu_{k+1},$ and between $\mu_l$ and $\mu_{l+1}$ correspondingly.
\end{proof}

In the remaining cases we get the following relations.

\begin{lemma}\label{Lemma: local relation E two terms}
Let $1\le k<n$ be an integer and $\mu$ be a composition of $d$ of length $n.$ Then one has 
\begin{equation*}
\mathfrak{e}_{\mu-\textbf{e}_{k+1},k}(p)\star\mathfrak{e}_{\mu-\textbf{e}_k,k}(q)=-\mathfrak{e}_{\mu-\textbf{e}_{k+1},k}(q-1)\star\mathfrak{e}_{\mu-\textbf{e}_k,k}(p+1).
\end{equation*}
\end{lemma}

\begin{proof}
Consider the intersection of preimages 
\begin{equation*}
X:=\pi_{12}^{-1}(\mathbb{O}_{\textbf{E}(\mu-\textbf{e}_{k+1},k)})\cap\pi_{23}^{-1}(\mathbb{O}_{\textbf{E}(\mu-\textbf{e}_k,k)})\subset \Fnd\times \Fnd\times \Fnd.
\end{equation*}
Note that according to Lemma \ref{lemma: transversality}, the intersection is transversal and one can apply Corollary \ref{corollary: convolution with restriction} with $U=\Fnd\times\Fnd$. Furthermore, the intersection consists of triples of flags $(\textbf{U},\textbf{V},\textbf{W})$ satisfying the following conditions. 
\begin{itemize}
\item Since $\textbf{V}\in \mathcal{F}_{\textbf{C}^{\textbf{E}(\mu-\textbf{e}_{k+1},k)}}=\Fmu,$ we get $\dim V_i/V_{i-1}=\mu_i$ for all $1\le i<n.$
\item Since $(\textbf{U},\textbf{V})\in \mathbb{O}_{\textbf{E}(\mu-\textbf{e}_{k+1},k)},$ we get $U_i=V_i$ for $i\neq k$ and $V_k\subset U_k,$ $\dim U_k/V_k=1.$ 
\item Since $(\textbf{V},\textbf{W})\in \mathbb{O}_{\textbf{E}(\mu-\textbf{e}_k,k)},$ we get $V_i=W_i$ for $i\neq k$ and $W_k\subset V_k,$ $\dim V_k/W_k=1.$ 
\end{itemize}
One concludes that $X\simeq \mathcal{F}_{(\mu_1,\ldots,\mu_k-1,1,1,\mu_{k+1}-1,\ldots,\mu_n)}.$ Moreover, one gets
\begin{equation*}
\pi_{12}^*(\mathfrak{e}_{\mu-\textbf{e}_{k+1},k}(p))=l_2^p,
\end{equation*}  
and
\begin{equation*}
\pi_{23}^*\mathfrak{e}_{\mu-\textbf{e}_k,k}(q)=l_1^q,
\end{equation*}
where $l_1$ and $l_2$ are the classes of the tautological line bundles on $X$ corresponding to the two one dimensional steps, left to right correspondingly. One concludes that 
\begin{equation*}
\mathfrak{e}_{\mu-\textbf{e}_{k+1},k}(p)\star\mathfrak{e}_{\mu-\textbf{e}_k,k}(q)=\pi_{13*}(l_1^ql_2^p),
\end{equation*}
and
\begin{equation*}
\mathfrak{e}_{\mu-\textbf{e}_{k+1},k}(q-1)\star\mathfrak{e}_{\mu-\textbf{e}_k,k}(p+1)=\pi_{13*}(l_1^{p+1}l_2^{q-1}).
\end{equation*}
The image $\pi_{13}(X)\subset \Fnd\times \Fnd$ is obtained from $X\simeq \mathcal{F}_{(\mu_1,\ldots,\mu_k-1,1,1,\mu_{k+1}-1,\ldots,\mu_n)}$ by forgetting the middle flag $\textbf{V}.$  Therefore, one gets $\pi_{13}(X)\simeq \mathcal{F}_{(\mu_1,\ldots,\mu_k-1,2,\mu_{k+1}-1,\ldots,\mu_n)},$ and the projection $\pi_{13}|_X$ is a natural projection between varieties of partial flags. The Lemma now follows immediately from Formulas (\ref{formula: push-forward Schur}) and (\ref{formula: general schur}).
\end{proof}

\begin{lemma}\label{Lemma: local relation E three terms}
Let $1\le k<n-1$ be an integer and $\mu$ be a composition of $d-1$ of length $n.$ Then one has 
\begin{equation*}
\mathfrak{e}_{\mu+\textbf{e}_k,k+1}(p)\star\mathfrak{e}_{\mu+\textbf{e}_{k+2},k}(q)=\mathfrak{e}_{\mu+\textbf{e}_{k+1},k}(q)\star\mathfrak{e}_{\mu+\textbf{e}_{k+1},k+1}(p)-\mathfrak{e}_{\mu+\textbf{e}_{k+1},k}(q+1)\star\mathfrak{e}_{\mu+\textbf{e}_{k+1},k+1}(p-1).
\end{equation*}
\end{lemma}

\begin{proof}
Consider the intersection of preimages 
\begin{equation*}
X:=\pi_{12}^{-1}(\mathbb{O}_{\textbf{E}(\mu+\textbf{e}_{k+1},k)})\cap\pi_{23}^{-1}(\mathbb{O}_{\textbf{E}(\mu+\textbf{e}_{k+1},k+1)})\subset \Fnd\times \Fnd\times \Fnd.
\end{equation*}
Note that according to Lemma \ref{lemma: transversality}, the intersection is transversal and one can apply Corollary \ref{corollary: convolution with restriction} with $U=\Fnd\times\Fnd$. Furthermore, the intersection consists of triples of flags $(\textbf{U},\textbf{V},\textbf{W})$ satisfying the following conditions:
\begin{itemize}
\item Since $\textbf{U}\in \mathcal{F}_{\textbf{R}^{\textbf{E}(\mu+\textbf{e}_{k+1},k)}}=\mathcal{F}_{\mu+\textbf{e}_{k+1}+\textbf{e}_k},$ we get 
$$
\dim U_i/U_{i-1}=\mu_i\ \ \textrm{for}\ \ 1\le i<k\ \ \textrm{and}\ \ k+1<i<n,
$$ 
$$
\dim U_{k}/U_{k-1}=\mu_k+1,\ \ \textrm{and}\ \ \dim U_{k+1}/U_k=\mu_{k+1}+1.
$$
\item Since $(\textbf{U},\textbf{V})\in \mathbb{O}_{\textbf{E}(\mu+\textbf{e}_{k+1},k)}$ we get $U_i=V_i$ for $i\neq k$ and $V_k\subset U_k,$ $\dim U_k/V_k=1.$
\item Since $(\textbf{V},\textbf{W})\in \mathbb{O}_{\textbf{E}(\mu+\textbf{e}_{k+1},k+1)}$ we get $V_i=W_i$ for $i\neq k+1$ and $W_{k+1}\subset V_{k+1},$ $\dim V_{k+1}/W_{k+1}=1.$
\end{itemize}
We use the following configuration space notation for this intersection:

\begin{equation*}
\begin{tikzpicture}
\draw (-3.2,0) node {$X=$};
\draw (-0.7,0) node {$\{U_0\subset\ldots\subset U_{k-1}\subset V_k $};
\draw (4.8,0) node {$U_{k+1}\subset\ldots\subset U_n\}.$};
\draw (2.23,0.5) node {$U_k$};
\draw (2.23,-0.5) node {$W_{k+1}$};

\draw (1.6,0.35) node[rotate=30] {$\subset$};
\draw (2.8,0.35) node[rotate=-30] {$\subset$};
\draw (1.5,-0.35) node[rotate=-30] {$\subset$};
\draw (2.9,-0.35) node[rotate=30] {$\subset$};
\end{tikzpicture} 
\end{equation*}
Moreover, one gets
\begin{equation*}
\pi_{12}^*(\mathfrak{e}_{\mu+\textbf{e}_{k+1},k}(q))=l^q,
\end{equation*}  
and
\begin{equation*}
\pi_{23}^*(\mathfrak{e}_{\mu+\textbf{e}_{k+1},k+1}(p))=m^p,
\end{equation*}
where $l=[U_k/V_k]\in K^G(X)$ and $m=[U_{k+1}/W_{k+1}]\in K^G(X).$ Note that $X$ is smooth. In fact, it is isomorphic to a $\mathbb{P}_{\mu_{k+1}+1}$-bundle over a variety of partial flags.
%$\mathcal{F}_{(\mu_1,\ldots,\mu_k,\mu_{k+1}+1,1,\mu_{k+2},\ldots,\mu_n)}$
The variety $X$ consists of two orbits: the closed orbit $D$ where $U_k\subset W_{k+1},$ and the open orbit where $U_k\cap W_{k+1}=V_k.$ Note that the closed orbit $D\subset X$ has codimension one. We conclude that 
\begin{align*}
\pi_{12}^*(\mathfrak{e}_{\mu+\textbf{e}_{k+1},k}(q))\pi_{23}^*(\mathfrak{e}_{\mu+\textbf{e}_{k+1},k+1}(p))&-\pi_{12}^*(\mathfrak{e}_{\mu+\textbf{e}_{k+1},k}(q+1))\pi_{23}^*(\mathfrak{e}_{\mu+\textbf{e}_{k+1},k+1}(p-1))\\
&=l^qm^p-l^{q+1}m^{p-1}=l^qm^p(1-lm^{-1})
\end{align*}

The key observation is that the class $l^qm^p(1-lm^{-1})$ is supported on the closed orbit $D\subset X.$ Indeed, consider the short exact sequence 

\begin{center}
\begin{tikzcd}
0 \arrow[r]& \mathcal{L}^{-1}\arrow[r] & \mathcal{O}_X \arrow[r] & \mathcal{O}_D \arrow[r] & 0.\\
\end{tikzcd} 
\end{center}
where $\mathcal{L}$ is the line bundle of the divisor $D\subset X.$ Then, in the $K$-theory we get 

\begin{equation*}
[\mathcal{O}_D]=1-[\mathcal{L}]^{-1},
\end{equation*}

Consider the natural map $\phi: U_k/V_k\to U_{k+1}/W_{k+1}$ defined as the composition of the embedding $U_k/V_k\hookrightarrow U_{k+1}/V_k$ and the projection $U_{k+1}/V_k\twoheadrightarrow U_{k+1}/W_{k+1}.$ Note that $D\subset X$ is the divisor of zeros of $\phi$. We conclude that
\begin{equation*}
[\mathcal{L}]=ml^{-1},
\end{equation*}
and
\begin{equation*}
[\mathcal{O}_D]=1-lm^{-1}.
\end{equation*}
Let $i:D\to X$ be the inclusion map. Denote $\bar{l}:=i^*l$ and $\bar{m}:=i^*m.$
Note that one also has $i^*[\mathcal{L}]=\bar{l}/\bar{m}$ and $i_*i^* \beta=[\mathcal{O}_D]\beta$ for any $\beta\in K^G(X).$ One gets
\begin{align*}
\mathfrak{e}_{\mu+\textbf{e}_{k+1},k}(q)&\star\mathfrak{e}_{\mu+\textbf{e}_{k+1},k+1}(p)-\mathfrak{e}_{\mu+\textbf{e}_{k+1},k}(q+1)\star\mathfrak{e}_{\mu+\textbf{e}_{k+1},k+1}(p-1)\\
&=\pi_{13*}\left(l^qm^p(1-lm^{-1})\right)=\pi_{13*}\left(l^qm^p[\mathcal{O}_D]\right)\\
&=\pi_{13*}\left(i_*i^*(l^qm^p)\right)=(\pi_{13}|_D)_*(\bar{l}^q\bar{m}^p).
\end{align*}

Let us now consider the convolution product $\mathfrak{e}_{\mu+\textbf{e}_k,k+1}(p)\star\mathfrak{e}_{\mu+\textbf{e}_{k+2},k}(q).$ The intersection of the preimages 
\begin{equation*}
Y:=\pi_{12}^{-1}(\mathbb{O}_{\textbf{E}(\mu+\textbf{e}_k,k+1)})\cap\pi_{23}^{-1}(\mathbb{O}_{\textbf{E}(\mu+\textbf{e}_{k+2},k)})\subset \Fnd\times \Fnd\times \Fnd.
\end{equation*}
Note that according to Lemma \ref{lemma: transversality}, the intersection is transversal and one can apply Corollary \ref{corollary: convolution with restriction} with $U=\Fnd\times\Fnd$. Furthermore, the intersection consists of triples of flags $(\textbf{U},\textbf{V},\textbf{W})$ satisfying the following conditions:
\begin{itemize}
\item Since $\textbf{U}\in \mathcal{F}_{\textbf{R}^{\textbf{E}(\mu+\textbf{e}_k,k+1)}}=\mathcal{F}_{\mu+\textbf{e}_{k+1}+\textbf{e}_k},$ we get 
$$
\dim U_i/U_{i-1}=\mu_i\ \ \textrm{for}\ \ 1\le i<k\ \ \textrm{and}\ \ k+1<i<n,
$$ 
$$
\dim U_{k}/U_{k-1}=\mu_k+1,\ \ \textrm{and}\ \ \dim U_{k+1}/U_k=\mu_{k+1}+1.
$$
\item Since $(\textbf{U},\textbf{V})\in \mathbb{O}_{\textbf{E}(\mu+\textbf{e}_k,k+1)}$ we get $U_i=V_i$ for $i\neq k+1$ and $V_{k+1}\subset U_{k+1},$ $\dim U_{k+1}/V_{k+1}=1.$
\item Since $(\textbf{V},\textbf{W})\in \mathbb{O}_{\textbf{E}(\mu+\textbf{e}_{k+2},k)}$ we get $V_i=W_i$ for $i\neq k$ and $W_k\subset V_k,$ $\dim V_k/W_k=1.$
\end{itemize}

In other words, one gets
\begin{equation*}
Y=\{U_0\subset\ldots\subset U_{k-1}\subset W_k\subset V_k=U_k\subset V_{k+1}\subset U_{k+1}\subset\ldots\subset U_n\},
\end{equation*}
and
\begin{equation*}
\pi_{12}^*(\mathfrak{e}_{\mu+\textbf{e}_k,k+1}(p))\pi_{23}^*(\mathfrak{e}_{\mu+\textbf{e}_{k+2},k}(q))=\hat{m}^p\hat{l}^q,
\end{equation*}
where $\hat{m}=[U_{k+1}/V_{k+1}]\in K^G(Y)$ and $\hat{l}=[V_k/W_k]\in K^G(Y).$
Note that both $Y$ and the closed orbit $D\subset X$ are varieties of partial flags. In fact, one gets
\begin{equation*}
Y\simeq D\simeq \mathcal{F}_{(\mu_1,\ldots,\mu_k,1,\mu_{k+1},1,\mu_{k+2},\ldots,\mu_n)},
\end{equation*} 
and the classes $\bar{l}$ and $\bar{m}$ correspond to $\hat{l}$ and $\hat{m}$ under the natural isomorphism. Moreover, the isomorphism commutes with the restrictions of the projection $\pi_{13}: \Fnd\times \Fnd\times \Fnd \to \Fnd\times \Fnd$ to $Y$ and to $D$ respectively. One concludes that
\begin{equation*}
(\pi_{13}|_D)_*(\bar{l}^q\bar{m}^p)=(\pi_{13}|_Y)_*(\hat{l}^q\hat{m}^p),
\end{equation*}
which completes the proof.
\end{proof}

\begin{lemma}\label{lemma: local E three term special}
Let $1\le k<n$ and let $\mu$ be a composition of $d$ of length $n$ such that $\mu_{k+1}=0.$ One has
\begin{equation*}
0=\mathfrak{e}_{\mu,k}(q)\star\mathfrak{e}_{\mu,k+1}(p)-\mathfrak{e}_{\mu,k}(q+1)\star\mathfrak{e}_{\mu,k+1}(p-1),
\end{equation*}   
for any $p,q\in\mathbb{Z}.$
\end{lemma}

\begin{remark}
This Lemma can considered as a special case of Lemma \ref{Lemma: local relation E three terms}: formally speaking, in this case the left hand side of Lemma \ref{Lemma: local relation E three terms} is undefined, as it would require flags with steps of negative codimension, but the right hand side is well defined and is equal to zero.
\end{remark}

\begin{proof}
Consider the preimage 

\begin{equation*}
X:=\pi_{12}^{-1}(\mathbb{O}_{\textbf{E}(\mu,k)})\cap\pi_{23}^{-1}(\mathbb{O}_{\textbf{E}(\mu,k+1)})\subset \Fnd\times \Fnd\times \Fnd.
\end{equation*}
Note that according to Lemma \ref{lemma: transversality}, the intersection is transversal. Furthermore, the intersection consists of triples of flags $(\textbf{U},\textbf{V},\textbf{W})$ satisfying conditions similar to those in the proof of Lemma \ref{Lemma: local relation E three terms}. The main difference is that this time $\dim U_{k+1}/U_k=\mu_{k+1}=0,$ which forces $U_k=U_{k+1}.$ One also gets $\dim W_{k+1}/W_k=0,$ or $W_{k+1}=W_k=V_k:$
\begin{equation*}
X=\{U_0\subset\ldots\subset U_{k-1}\subset V_k=W_k=W_{k+1}\subset U_k=U_{k+1}=V_{k+1}\subset U_{k+2}\subset\ldots\subset U_n\}.
\end{equation*}
In particular,
\begin{equation*}
\pi_{12}^*(\mathfrak{e}_{\mu,k}(q))=l^q,\ \textrm{and}\ \pi_{23}^*(\mathfrak{e}_{\mu,k+1}(p))=l^p,
\end{equation*}
where $l=[U_k/V_k]=[V_{k+1}/W_{k+1}]\in K^G(X).$ Therefore,
\begin{equation*}
\mathfrak{e}_{\mu,k}(q)\star\mathfrak{e}_{\mu,k+1}(p)-\mathfrak{e}_{\mu,k}(q+1)\star\mathfrak{e}_{\mu,k+1}(p-1)=\pi_{13*}(l^{q+p}-l^{(q+1)+(p-1)})=0.
\end{equation*}
\end{proof}

Similarly, for the $\mathfrak{f}$ classes, one gets

\begin{lemma} Let $1\le k<n$ be an integer and $\mu$ be a composition of $d-1$ of length $n.$ Then one has 
\begin{equation*}
\mathfrak{f}_{\mu+\textbf{e}_{k+1},k}(p)\star\mathfrak{f}_{\mu+\textbf{e}_k,k}(q)=-\mathfrak{f}_{\mu+\textbf{e}_{k+1},k}(q+1)\star\mathfrak{f}_{\mu+\textbf{e}_k,k}(p-1).
\end{equation*}
\end{lemma}

\begin{proof}
We skip the proof as it is very similar to the proof of Lemma \ref{Lemma: local relation E two terms}.
\end{proof}

\begin{lemma} Let $1\le k<n-1$ be an integer and $\mu$ be a composition of $d-1$ of length $n.$ Then one has 
\begin{equation*}
\mathfrak{f}_{\mu+\textbf{e}_{k+2},k}(p)\star\mathfrak{f}_{\mu+\textbf{e}_k,k+1}(q)=\mathfrak{f}_{\mu+\textbf{e}_{k+1},k+1}(q)\star\mathfrak{f}_{\mu+\textbf{e}_{k+1},k}(p)-\mathfrak{f}_{\mu+\textbf{e}_{k+1},k+1}(q-1)\star\mathfrak{f}_{\mu+\textbf{e}_{k+1},k}(p+1).
\end{equation*}
\end{lemma}

\begin{proof}
We skip the proof as it is very similar to the proof of Lemma \ref{Lemma: local relation E three terms}.
\end{proof}

\begin{lemma}
Let $1\le k<n$ and let $\mu$ be a composition of $d$ of length $n$ such that $\mu_{k+1}=0.$ One has
\begin{equation*}
0=\mathfrak{f}_{\mu,k+1}(q)\star\mathfrak{f}_{\mu,k}(p)-\mathfrak{f}_{\mu,k+1}(q-1)\star\mathfrak{f}_{\mu,k}(p+1),
\end{equation*}   
for any $p,q\in\mathbb{Z}.$
\end{lemma}

\begin{proof}
We skip the proof as it is very similar to the proof of Lemma \ref{lemma: local E three term special}.
\end{proof}

Finally, between the $\mathfrak{e}$ and $\mathfrak{f}$ classes one gets

\begin{lemma}\label{lemma: local H as commutator} 
Let $1\le k<n$ be an integer and $\mu$ be a composition of $d$ of length $n$ such that $\mu_k>0$ and $\mu_{k+1}>0.$ Then the class 
\begin{equation*}
\mathfrak{h}_{\mu,k}(p+q):=\mathfrak{e}_{\mu-\textbf{e}_k,k}(p)\star\mathfrak{f}_{\mu-\textbf{e}_k,k}(q)-\mathfrak{f}_{\mu-\textbf{e}_{k+1},k}(q)\star\mathfrak{e}_{\mu-\textbf{e}_{k+1},k}(p)
\end{equation*}
is supported on $\mathbb{O}_{Diag(\mu)}\subset \Fnd\times \Fnd$ and only depends on the sum $p+q.$ Furthermore, 

\begin{equation*}
\mathfrak{h}_{\mu,k}(p+q)=
\left\{
\begin{array}{rl}
(-1)^{\mu_k-1}[\Lambda_{\mu_k}(\mathcal{T}^*_k)][Sym_{-(p+q)-\mu_k}(\mathcal{T}^*_{k+1}\oplus \mathcal{T}^*_k)],& p+q\le -\mu_k,\\
0,& -\mu_k<p+q<\mu_{k+1},\\
(-1)^{\mu_{k+1}}[\Lambda_{\mu_{k+1}}(\mathcal{T}_{k+1})][Sym_{p+q-\mu_{k+1}}(\mathcal{T}_{k+1}\oplus \mathcal{T}_k)],& p+q\ge \mu_{d+1},
\end{array}
\right.
\end{equation*}
where $\mathcal{T}_k$ and $\mathcal{T}_{k+1}$ are the tautological bundles over $\mathbb{O}_{Diag(\mu)}\simeq \Fmu$ corresponding to the $k$th and $(k+1)$th steps respectively.
\end{lemma}

\begin{proof}
Consider first the product $\mathfrak{e}_{\mu-\textbf{e}_k,k}(p)\star\mathfrak{f}_{\mu-\textbf{e}_k,k}(q).$ The intersection of preimages 
\begin{equation*}
X:=\pi_{12}^{-1}(\mathbb{O}_{\textbf{E}(\mu-\textbf{e}_k,k)})\cap\pi_{23}^{-1}(\mathbb{O}_{\textbf{F}(\mu-\textbf{e}_k,k)})\subset \Fnd\times \Fnd\times \Fnd
\end{equation*}
is transversal according to Lemma \ref{lemma: transversality} and consists of triples of flags $(\textbf{U},\textbf{V},\textbf{W})$ satisfying the following conditions:
\begin{itemize}
\item Since $\textbf{U}\in \mathcal{F}_{\textbf{R}^{\textbf{E}(\mu-\textbf{e}_k,k)}}=\Fmu,$ we get $\dim U_i/U_{i-1}=\mu_i$ for all $1<i\le n.$
\item Since $(\textbf{U},\textbf{V})\in \mathbb{O}_{\textbf{E}(\mu-\textbf{e}_k,k)}$ we get $U_i=V_i$ for $i\neq k$ and $V_k\subset U_k,$ $\dim U_k/V_k=1.$
\item Since $(\textbf{V},\textbf{W})\in \mathbb{O}_{\textbf{F}(\mu-\textbf{e}_k,k)}$ we get $V_i=W_i$ for $i\neq k$ and $V_k\subset W_k,$ $\dim W_k/V_k=1.$
\end{itemize}
We use the following configuration space notation for this intersection:

\begin{equation*}
\begin{tikzpicture}
\draw (-3.9,0) node {$X
%:=\pi_{12}^{-1}(\textbf{F}(\mu-\textbf{e}_{k+1},k))\cap\pi_{23}^{-1}(\textbf{E}(\mu-\textbf{e}_{k+1},k))
=$};
\draw (-1.2,0) node {$\{U_0\subset\ldots\subset U_{k-1}\subset V_k$};
\draw (4.2,0) node {$U_{k+1}\subset\ldots\subset U_n\}.$};
\draw (1.73,0.5) node {$U_k$};
\draw (1.73,-0.5) node {$W_k$};

\draw (1.1,0.35) node[rotate=30] {$\subset$};
\draw (2.3,0.35) node[rotate=-30] {$\subset$};
\draw (1.1,-0.35) node[rotate=-30] {$\subset$};
\draw (2.3,-0.35) node[rotate=30] {$\subset$};
\end{tikzpicture} 
\end{equation*}

Using Corollary \ref{corollary: convolution with restriction} with $U=\Fnd\times\Fnd$ we immediately get
\begin{equation*}
\pi_{12}^*(\mathfrak{e}_{\mu-\textbf{e}_k,k}(p))\pi_{23}^*(\mathfrak{f}_{\mu-\textbf{e}_k,k}(q))=l_1^pl_2^q\in K^G(X)\subset K^G(\Fnd\times \Fnd\times \Fnd),
\end{equation*}
where $l_1$ and $l_2$ are the classes of the line bundles $U_k/V_k$ and $W_k/V_k$ respectively. Note that $X$ is a smooth variety. In fact, it is isomorphic to a $\mathbb{P}_{\mu_{k+1}}$-bundle over the variety of partial flags $\mathcal{F}_{(\mu_1,\ldots,\mu_k-1,1,\mu_{k+1},\ldots,\mu_n)}.$ The variety $X$ consists of two orbits: the open orbit, where $U_k\cap W_k=V_k$ and the closed orbit, where $U_k=W_k.$ Note that the closed orbit has codimention $\mu_{k+1}$ in $X.$ 

The image $Q:=\pi_{13}(X)\subset \Fnd\times \Fnd$ of $X$ is obtained by forgetting the middle flag $\textbf{V}.$ One gets

\begin{center}
\begin{tikzpicture}
\draw (-2.9,0) node {$Q:=$};
\draw (-0.75,0) node {$\{U_0\subset\ldots\subset U_{k-1}$};
\draw (4.3,0) node {$U_{k+1}\subset\ldots\subset U_n\}.$};
\draw (1.73,0.5) node {$U_k$};
\draw (1.73,-0.5) node {$W_k$};

\draw (1.1,0.35) node[rotate=30] {$\subset$};
\draw (2.3,0.35) node[rotate=-30] {$\subset$};
\draw (1.1,-0.35) node[rotate=-30] {$\subset$};
\draw (2.3,-0.35) node[rotate=30] {$\subset$};
\end{tikzpicture}
\end{center}
Here as before $U_i=W_i$ for $i\neq k,$ $\dim U_i/U_{i-1}=\mu_i,$ $\dim W_k=\dim U_k,$ and either $W_k=U_k$ or $W_k\cap U_k$ has codimension $1$ in both $W_k$ and $U_k.$ Note that $Q\subset \Fnd\times \Fnd$ also consists of two orbits: the closed orbit where $U_k=W_k$ and the complement to it. However, $Q$ is not necessarily smooth.

Similarly, for the product $\mathfrak{f}_{\mu-\textbf{e}_{k+1},k}(q)\star\mathfrak{e}_{\mu-\textbf{e}_{k+1},k}(p)$ we get that the intersection of preimages 
\begin{equation*}
Y:=\pi_{12}^{-1}(\mathbb{O}_{\textbf{F}(\mu-\textbf{e}_{k+1},k)})\cap\pi_{23}^{-1}(\mathbb{O}_{\textbf{E}(\mu-\textbf{e}_{k+1},k)})\subset \Fnd\times \Fnd\times \Fnd
\end{equation*}
is also transversal according to Lemma \ref{lemma: transversality} and consists of triples of flags $(\textbf{U},\textbf{V}',\textbf{W})$ satisfying the following conditions:
\begin{itemize}
\item Since $\textbf{U}\in \mathcal{F}_{\textbf{R}^{\textbf{F}(\mu-\textbf{e}_{k+1},k)}}=\Fmu,$ we get $\dim U_i/U_{i-1}=\mu_i$ for all $1<i\le n.$
\item Since $(\textbf{U},\textbf{V}')\in \mathbb{O}_{\textbf{F}(\mu-\textbf{e}_{k+1},k)}$ we get $U_i=V'_i$ for $i\neq k$ and $U_k\subset V'_k,$ $\dim V'_k/U_k=1.$
\item Since $(\textbf{V}',\textbf{W})\in \mathbb{O}_{\textbf{E}(\mu-\textbf{e}_{k+1},k)}$ we get $V'_i=W_i$ for $i\neq k$ and $W_k\subset V'_k,$ $\dim V'_k/W_k=1.$
\end{itemize}
In other words, using Corollary \ref{corollary: convolution with restriction} with $U=\Fnd\times\Fnd$ we get:
\begin{equation*}
\begin{tikzpicture}
\draw (-2.8,0) node {$Y
%:=\pi_{12}^{-1}(\textbf{E}(\mu-\textbf{e}_k,k))\cap\pi_{23}^{-1}(\textbf{F}(\mu-\textbf{e}_k,k))
=$};
\draw (-0.7,0) node {$\{U_0\subset\ldots\subset U_{k-1}$};
\draw (4.8,0) node {$V'_k\subset U_{k+1}\subset\ldots\subset U_n\}.$};
\draw (1.73,0.5) node {$U_k$};
\draw (1.73,-0.5) node {$W_k$};

\draw (1.1,0.35) node[rotate=30] {$\subset$};
\draw (2.3,0.35) node[rotate=-30] {$\subset$};
\draw (1.1,-0.35) node[rotate=-30] {$\subset$};
\draw (2.3,-0.35) node[rotate=30] {$\subset$};
\end{tikzpicture} 
\end{equation*}

Here we abuse notations by using $\textbf{U}$ and $\textbf{W}$ both in $X$ and $Y.$ This is justified because the corresponding flags are identified by the restrictions of the projection $\pi_{13}.$ Indeed, the image of $Y$ under $\pi_{13}$ is obtained by forgetting the flag $\textbf{V}'$ and coincides with $Q\subset \Fnd\times \Fnd,$ same as the image of $X.$ We also immediately get
\begin{equation*}
\pi_{12}^*(\mathfrak{f}_{\mu-\textbf{e}_{k+1},k}(q))\pi_{23}^*(\mathfrak{e}_{\mu-\textbf{e}_{k+1},k}(p))=m_1^q m_2^p\in K^G(Y)\subset K^G(\Fnd\times \Fnd\times \Fnd),
\end{equation*}
where $m_1$ and $m_2$ are the classes of the line bundles $V'_k/U_k$ and $V'_k/W_k$ respectively. Note that $Y$ is also smooth and isomorphic to a $\mathbb{P}_{\mu_k}$-bundle over the variety of partial flags $\mathcal{F}_{(\mu_1,\ldots,\mu_k,1,\mu_{k+1}-1,\ldots,\mu_n)}.$ The variety $Y$ also consists of two orbits: the open orbit, where $U_k+W_k=V'_k$ and the closed orbit, where $U_k=W_k.$ The closed orbit has codimention $\mu_k$ in $Y.$

%The images of $X$ and $Y$ under the projection $\pi_{13}$ coincide: one obtains pairs of flags $(\textbf{L,M}),$ where $\dim L_i/L_{i-1}=\dim M_i/M_{i-1}=\mu_i$ for all $0<i\le n,$ $L_i=M_i$ for $i\neq k,$ and $L_k\cap M_k$ has codimension at most $1$ in $L_k$ (and $M_k$). In other words, we obtain

%\begin{center}
%\begin{tikzpicture}
%\draw (-2.9,0) node {$T:=$};
%\draw (-0.75,0) node {$\{L_0\subset\ldots\subset L_{k-1}$};
%\draw (4.3,0) node {$L_{k+1}\subset\ldots\subset L_n\},$};
%\draw (1.73,0.5) node {$L_k$};
%\draw (1.73,-0.5) node {$M_k$};

%\draw (1.1,0.35) node[rotate=30] {$\subset$};
%\draw (2.3,0.35) node[rotate=-30] {$\subset$};
%\draw (1.1,-0.35) node[rotate=-30] {$\subset$};
%\draw (2.3,-0.35) node[rotate=30] {$\subset$};
%\end{tikzpicture}
%\end{center}
%with the extra condition $\dim L_k/(L_k\cap M_k)\le 1.$ Note that $T\subset \Fnd\times \Fnd$ also consists of two orbits: the open orbit where $\dim L_k/(L_k\cap M_k)=1$ and the closed orbit where $L_k=M_k.$ However, $T$ is not necessarily smooth. 

Consider the configuration space

\begin{center}
\begin{tikzpicture}
\draw (-3.9,0) node {$Z:=$};
\draw (-1.3,0) node {$\{U_0\subset\ldots\subset U_{k-1}\subset V_k$};
\draw (4.8,0) node {$V'_k\subset U_{k+1}\subset\ldots\subset U_n\},$};
\draw (1.73,0.5) node {$U_k$};
\draw (1.73,-0.5) node {$W_k$};

\draw (1.1,0.35) node[rotate=30] {$\subset$};
\draw (2.3,0.35) node[rotate=-30] {$\subset$};
\draw (1.1,-0.35) node[rotate=-30] {$\subset$};
\draw (2.3,-0.35) node[rotate=30] {$\subset$};
\end{tikzpicture}
\end{center}

\noindent where as before $\dim U_i/U_{i-1}=\mu_i$ for $1<i\le n,$ and 
\begin{equation*}
\dim U_k/V_k=\dim W_k/V_k=\dim V'_k/U_k=\dim V'_k/W_k=1.
\end{equation*}

Similar to $X$ and $Y,$ $Z$ is a smooth variety, isomorphic to a $\mathbb{P}_1$-bundle over the variety of partial flags $\mathcal{F}_{(\mu_1,\ldots,\mu_k-1,1,1,\mu_{k+1}-1,\ldots,\mu_n)}.$ It also consists of two orbits: $U_k\neq W_k$ and $U_k=W_k.$ However, in this case the closed orbit is a hypersurface.

We are abusing notations again by recycling notations for $\textbf{U},\textbf{W},\textbf{V},$ and $\textbf{V}'.$ This is justified by the identifications from the following natural projections
\begin{center}
\begin{tikzpicture}
\node (T1) at (2,3.7) {};
\node (T2) at (6,3.7) {};

\node (T) at (4,4) {
\begin{tikzpicture}
\draw (-1.8,0) node {$Z=\{$};
\draw (-0.1,0) node {\scriptsize $\ldots\subset U_{k-1}\subset V_k$};
\draw (3.65,0) node {{\scriptsize $V'_k\subset U_{k+1}\subset\ldots$}$\}$};
\draw (1.73,0.4) node {\scriptsize $U_k$};
\draw (1.73,-0.4) node {\scriptsize $W_k$};

\draw (1.3,0.3) node[rotate=30] {\scriptsize $\subset$};
\draw (2.1,0.3) node[rotate=-30] {\scriptsize $\subset$};
\draw (1.3,-0.3) node[rotate=-30] {\scriptsize $\subset$};
\draw (2.1,-0.3) node[rotate=30] {\scriptsize $\subset$};
\end{tikzpicture}
};
\node (M1) at (0,2) {
\begin{tikzpicture}
\draw (-1.82,0) node {$X=\{$};
\draw (0,0) node {\scriptsize $\ldots\subset U_{k-1}\subset V_k$};
\draw (3.2,0) node {{\scriptsize $U_{k+1}\subset\ldots$}$\}$};
\draw (1.73,0.4) node {\scriptsize $U_k$};
\draw (1.73,-0.4) node {\scriptsize $W_k$};

\draw (1.3,0.3) node[rotate=30] {\scriptsize $\subset$};
\draw (2.1,0.3) node[rotate=-30] {\scriptsize $\subset$};
\draw (1.3,-0.3) node[rotate=-30] {\scriptsize $\subset$};
\draw (2.1,-0.3) node[rotate=30] {\scriptsize $\subset$};
\end{tikzpicture} 
};
\node (M2) at (8,2) {
\begin{tikzpicture}
\draw (-0.95,0) node {$Y=\{$};
\draw (0.4,0) node {\scriptsize $\ldots\subset U_{k-1}$};
\draw (3.6,0) node {{\scriptsize $V'_k\subset U_{k+1}\subset\ldots$}$\}$};
\draw (1.73,0.4) node {\scriptsize $U_k$};
\draw (1.73,-0.4) node {\scriptsize $W_k$};

\draw (1.3,0.3) node[rotate=30] {\scriptsize $\subset$};
\draw (2.1,0.3) node[rotate=-30] {\scriptsize $\subset$};
\draw (1.3,-0.3) node[rotate=-30] {\scriptsize $\subset$};
\draw (2.1,-0.3) node[rotate=30] {\scriptsize $\subset$};
\end{tikzpicture} 
};
\draw
(T1) edge[->,>=angle 90] node[left] {$b_1$} (M1)
(T2) edge[->,>=angle 90] node[right] {$b_2$} (M2);
\end{tikzpicture}
\end{center}
where $b_1$ forgets $V'_k$ and $b_2$ forgets $V_k.$ The restrictions of $b_1$ and $b_2$ to the open orbit $\{U_k\neq W_k\}$ are both isomorphisms to the open orbits in $X$ and $Y$ respectively. Indeed, if $U_k\neq W_k,$ then $V_k=U_k\cap W_k$ and $V'_k=U_k+W_k.$ In fact, these maps are blow ups of $X$ and $Y$ respectively along the closed orbits. It follows that both $b_1$ and $b_2$ are proper birational maps. In particular, $(b_j)_*(b_j)^*=id$ in $K$-theory for $j=1,2.$ 

The compositions $\pi_{13}\circ b_1:Z\to Q\subset \Fnd\times \Fnd$ and $\pi_{13}\circ b_2:Z\to Q\subset \Fnd\times \Fnd$ coincide. Indeed, both maps forget both $V'_k$ and $V_k.$ Denote $p:=\pi_{13}\circ b_1=\pi_{13}\circ b_2:Z\to \Fnd\times \Fnd.$ Also, let

\begin{align*}
\hat{l}_1&:=b_1^*(l_1)=[U_k/V_k]\in K^G(Z),\\
\hat{l}_2&:=b_1^*(l_2)=[W_k/V_k]\in K^G(Z),\\
\hat{m}_1&:=b_2^*(m_1)=[V'_k/U_k]\in K^G(Z),\\
\hat{m}_2&:=b_2^*(m_2)=[V'_k/W_k]\in K^G(Z).
\end{align*}

Combining the above, we obtain:

\begin{align*}
&\mathfrak{e}_{\mu-\textbf{e}_k,k}(p)\star\mathfrak{f}_{\mu-\textbf{e}_k,k}(q)-\mathfrak{f}_{\mu-\textbf{e}_{k+1},k}(q)\star\mathfrak{e}_{\mu-\textbf{e}_{k+1},k}(p)=\pi_{13*}(l_1^pl_2^q-m_1^q m_2^p)\\
&=\pi_{13*}\left[b_{1*}(\hat{l}_1^p\hat{l}_2^q)-b_{2*}(\hat{m}_1^q\hat{m}_2^p)\right]=p_*\left(\hat{l}_1^p\hat{l}_2^q-\hat{m}_1^q\hat{m}_2^p\right).
\end{align*}

The key observation is that the class $\hat{l}_1^p\hat{l}_2^q-\hat{m}_1^q\hat{m}_2^p\in K^G(Z)$ is supported on the closed orbit:

\begin{equation*}
D:=\{\ldots\subset U_{k-1}\subset V'_k\subset U_k=W_k\subset V_k\subset U_{k+1}\subset\ldots\}\subset Z.
\end{equation*}

Indeed, consider the short exact sequence 

\begin{center}
\begin{tikzcd}
0 \arrow[r]& \mathcal{L}^{-1}\arrow[r] & \mathcal{O}_Z \arrow[r] & \mathcal{O}_D \arrow[r] & 0,\\
\end{tikzcd} 
\end{center}
where $\mathcal{L}$ is the line bundle of the divisor $D\subset Z.$ Then, in the $K$-theory we get 

\begin{equation*}
[\mathcal{O}_D]=1-[\mathcal{L}]^{-1},
\end{equation*}

Consider the natural map $\phi_1: U_k/V_k\to V'_k/W_k$ defined as the composition of the embedding $U_k/V_k\hookrightarrow W_k/V_k$ and the projection $W_k/V_k\twoheadrightarrow V'_k/W_k.$ The natural map $\phi_2:W_k/V_k\to V'_k/U_k$ is constructed in a similar manner. Note that $D\subset Z$ 
is the divisor of zeros for both maps. We conclude that
\begin{equation*}
[\mathcal{L}]=\hat{m}_2\hat{l}_1^{-1}=\hat{m}_1\hat{l}_2^{-1},
\end{equation*}
and 
\begin{equation*}
\hat{l}_1^p\hat{l}_2^q-\hat{m}_1^q\hat{m}_2^p=\hat{l}_1^p\hat{l}_2^q\Bigl(1-\left(\frac{\hat{m}_1}{\hat{l}_2}\right)^q\left(\frac{\hat{m}_2}{\hat{l}_1}\right)^p\Bigr)=\hat{l}_1^p\hat{l}_2^q\left(1-[\mathcal{L}]^{p+q}\right).
\end{equation*}

Finally, $1-\mathcal{L}^{p+q}$ is always divisible by $[\mathcal{O}_D]=1-[\mathcal{L}]^{-1}$ in $K^G(Z).$ Let $i:D\to Z$ be the inclusion map. Denote

\begin{align*}
%&\alpha_r:=\frac{1-l^{-r}}{1-l^{-1}}\in K^G(Z),\\
%&\bar{\alpha}_r:=i^*\alpha_r,\\
%&\bar{p}:=p\circ i,\\
&l:=i^*\hat{l}_1=i^*\hat{l}_2,\\
&m:=i^*\hat{m}_1=i^*\hat{m}_2,\\
&\pi:=p|_D=p\circ i.
\end{align*}

Note that one also has $i^*[\mathcal{L}]=m/l$ and $i_*i^* \beta=[\mathcal{O}_D]\beta$ for any $\beta\in K^G(Z).$ One gets
\begin{align*}
p_*\left(\hat{l}_1^p\hat{l}_2^q\left(1-[\mathcal{L}]^{p+q}\right)\right)
&=p_*\left(\hat{l}_1^p\hat{l}_2^q[\mathcal{O}_D]\frac{1-[\mathcal{L}]^{p+q}}{1-[\mathcal{L}]^{-1}}\right)\\
&=p_*\left(i_*i^*\left[\hat{l}_1^p\hat{l}_2^q\frac{1-[\mathcal{L}]^{-(p+q)}}{1-[\mathcal{L}]^{-1}}\right]\right)\\
&=(p_*i_*)\left(l^{p+q}\frac{1-(m/l)^{p+q}}{1-(m/l)^{-1}}\right)\\
&=\pi_*\left(-m\frac{l^{p+q}-m^{p+q}}{l-m}\right).
\end{align*}

Note that the image $\pi(D)\subset Q\subset \Fnd\times \Fnd$ is the closed orbit 

\begin{equation*}
\pi(D)=\{\ldots\subset U_{k-1}\subset U_k\subset U_{k+1}\subset\ldots\}=\mathbb{O}_{Diag(\mu)}\simeq \Fmu,
\end{equation*}
and
\begin{equation*}
D=\{\ldots\subset U_{k-1}\subset V_k\subset U_k\subset V'_k\subset U_{k+1}\subset\ldots\}\simeq \mathcal{F}_{(\mu_1,\ldots,\mu_k-1,1,1,\mu_{k+1}-1,\ldots,\mu_n)},
\end{equation*}
with the map $\pi$ forgetting both $V_k$ and $V'_k.$ One can factor $\pi=\pi^-\circ \pi^+$ where $\pi^-$ forgets $V_k$ and and $\pi^+$ forgets $V'_k,$ and apply Corollaries \ref{corollary: pi_*^-} and \ref{corollary: pi_*^+} and Lemma \ref{lemma: push-forward for other steps} to compute $\pi_*(l^r)$ and $\pi_*(m^r)$ for all $r\in\mathbb{Z}.$ In order to complete the computations, we need to consider the cases $p+q=0,$ $p+q>0,$ and $p+q<0$ separately. Recall that $\mathcal{T}_k$ and $\mathcal{T}_{k+1}$ be the tautological bundles over $\mathbb{O}_{Diag(\mu)}\simeq \Fmu$ corresponding to the $k$th and $(k+1)$th steps. 

\begin{enumerate}
\item[(0)] Let $p+q=0.$ In this case one gets $-m\frac{l^{p+q}-m^{p+q}}{l-m}=0,$ therefore
\begin{equation*}
\mathfrak{e}_{\mu-\textbf{e}_k,k}(p)\star\mathfrak{f}_{\mu-\textbf{e}_k,k}(q)-\mathfrak{f}_{\mu-\textbf{e}_{k+1},k}(q)\star\mathfrak{e}_{\mu-\textbf{e}_{k+1},k}(p)=0.
\end{equation*}
 
\item Let $p+q>0,$ then one gets
\begin{align*}
\pi_*&\left(-m\frac{l^{p+q}-m^{p+q}}{l-m}\right)=-\pi_*\left(l^{p+q-1}m+l^{p+q-2}m^2+\ldots+lm^{p+q-1}+m^{p+q}\right)\\
&=(-1)^{\mu_{k+1}}[\Lambda_{\mu_{k+1}}(\mathcal{T}_{k+1})]\sum\limits_{j=0}^{p+q-\mu_{k+1}}[Sym_j(\mathcal{T}_k)][Sym_{p+q-\mu_{k+1}-j}(\mathcal{T}_{k+1})]\\
&=(-1)^{\mu_{k+1}}[\Lambda_{\mu_{k+1}}(\mathcal{T}_{k+1})][Sym_{p+q-\mu_{k+1}}(\mathcal{T}_{k+1}\oplus \mathcal{T}_k)].
\end{align*}

\item Let $p+q<0$ and set $s=-(p+q)>0.$ One gets
\begin{align*}
\pi_*&\left(-m\frac{l^{p+q}-m^{p+q}}{l-m}\right)=\pi_*\left(-m\frac{\frac{1}{l^s}-\frac{1}{m^s}}{l-m}\right)=\pi_*\left(-\frac{m^s-l^s}{l^s m^{s-1}(l-m)}\right)\\
&=\pi_*(\frac{l^{s-1}+l^{s-2}m+\ldots+m^{s-1}}{l^s m^{s-1}})\\
&=\pi_*(l^{-1}m^{1-s}+l^{-2}m^{2-s}+\ldots+l^{1-s}m^{-1}+l^{-s})\\
&=(-1)^{\mu_k-1}[\Lambda_{\mu_k}(\mathcal{T}^*_k)]\sum\limits_{j=0}^{s-\mu_k}[Sym_j(\mathcal{T}^*_{k+1})][Sym_{s-\mu_k-j}(\mathcal{T}^*_k)]\\
&=(-1)^{\mu_k-1}[\Lambda_{\mu_k}(\mathcal{T}^*_k)][Sym_{s-\mu_k}(\mathcal{T}^*_{k+1}\oplus \mathcal{T}^*_k)]\\
&=(-1)^{\mu_k-1}[\Lambda_{\mu_k}(\mathcal{T}^*_k)][Sym_{-(p+q)-\mu_k}(\mathcal{T}^*_{k+1}\oplus \mathcal{T}^*_k)].
\end{align*}

\end{enumerate} 

Note that for $0<p+q<\mu_{k+1}$ one has $[Sym_{p+q-\mu_{k+1}}(\mathcal{T}_{k+1}\oplus \mathcal{T}_k)]=0,$ and for $-\mu_k<p+q<0$ one has $[Sym_{-(p+q)-\mu_k}(\mathcal{T}^*_{k+1}\oplus \mathcal{T}^*_k)]=0,$ which concludes the proof. 

\end{proof}

The next two Lemmas can be considered as special cases of Lemma \ref{lemma: local H as commutator}:

\begin{lemma}\label{lemma: local H as commutator spec 1} 
Let $1\le k<n$ be an integer and $\mu$ be a composition of $d$ of length $n$ such that $\mu_k=0$ and $\mu_{k+1}>0.$ Then the class 
\begin{equation*}
\mathfrak{h}_{\mu,k}(p+q):=-\mathfrak{f}_{\mu-\textbf{e}_{k+1},k}(q)\star\mathfrak{e}_{\mu-\textbf{e}_{k+1},k}(p)
\end{equation*}
is supported on $\mathbb{O}_{Diag(\mu)}\subset \Fnd\times \Fnd$ and only depends on the sum $p+q.$ Furthermore, 

\begin{equation*}
\mathfrak{h}_{\mu,k}(p+q)=
\left\{
\begin{array}{rl}
-[Sym_{-(p+q)}(\mathcal{T}^*_{k+1})],& p+q\le 0,\\
0,& 0<p+q<\mu_{k+1},\\
(-1)^{\mu_{k+1}}[\Lambda_{\mu_{k+1}}(\mathcal{T}_{k+1})][Sym_{p+q-\mu_{k+1}}(\mathcal{T}_{k+1})],& p+q\ge \mu_{k+1},
\end{array}
\right.
\end{equation*}
where $\mathcal{T}_{k+1}$ is the tautological bundles over $\mathbb{O}_{Diag(\mu)}\simeq \Fmu$ corresponding to the $(k+1)$th step of the flag.
\end{lemma}

\begin{remark}
Note that in the case when $\mu_k=0$ the first term on the right of the relation in Lemma \ref{lemma: local H as commutator} is not well defined. Lemma \ref{lemma: local H as commutator spec 1} is obtained by replacing it with zero.
\end{remark}

\begin{proof}
Consider the intersection of preimages 
\begin{equation*}
Y:=\pi_{12}^{-1}(\mathbb{O}_{\textbf{F}(\mu-\textbf{e}_{k+1},k)})\cap\pi_{23}^{-1}(\mathbb{O}_{\textbf{E}(\mu-\textbf{e}_{k+1},k)})\subset \Fnd\times \Fnd\times \Fnd.
\end{equation*}
Note that according to Lemma \ref{lemma: transversality}, the intersection is transversal. Similar to Lemma \ref{lemma: local H as commutator}, we get triples of flags $(\textbf{U},\textbf{V}',\textbf{W})$ satisfying the following conditions:
\begin{itemize}
\item Since $\textbf{U}\in \mathcal{F}_{\textbf{R}^{\textbf{F}(\mu-\textbf{e}_{k+1},k)}}=\Fmu,$ we get $\dim U_i/U_{i-1}=\mu_i$ for all $1<i\le n.$ In particular, we get $U_k=U_{k-1}$ as $\mu_k=0.$
\item Since $(\textbf{U},\textbf{V}')\in \mathbb{O}_{\textbf{F}(\mu-\textbf{e}_{k+1},k)}$ we get $U_i=V'_i$ for $i\neq k$ and $U_k\subset V'_k,$ $\dim V'_k/U_k=1.$
\item Since $(\textbf{V}',\textbf{W})\in \mathbb{O}_{\textbf{E}(\mu-\textbf{e}_{k+1},k)}$ we get $V'_i=W_i$ for $i\neq k$ and $W_k\subset V'_k,$ $\dim V'_k/W_k=1.$ Moreover, we get $W_k=W_{k-1}=U_{k-1}=U_k.$
\end{itemize}
In other words, using Corollary \ref{corollary: convolution with restriction} with $U=\Fnd\times\Fnd$ we get:
\begin{equation*}
Y=\{U_0\subset\ldots\subset U_k=W_k\subset V'_k\subset U_{k+1}\subset\ldots\subset U_n\}\simeq \mathcal{F}_{(\mu_1,\ldots,\mu_k,1,\mu_{k+1}-1,\mu_{k+2},\ldots,\mu_n)}.
\end{equation*}
Moreover,
\begin{equation*}
\pi_{12}^*(\mathfrak{f}_{\mu-\textbf{e}_{k+1},k}(q))=l^q,
\end{equation*}
and
\begin{equation*}
\pi_{23}^*(\mathfrak{e}_{\mu-\textbf{e}_{k+1},k}(p))=l^p,
\end{equation*}
where $l=[V'_k/U_k]=[V'_k/W_k]\in K^G(Y).$ The projection $\pi_{13}|_Y$ forgets the flag $\textbf{V}'.$ One gets
\begin{equation*}
\pi_{13}(Y)=\{U_0\subset\ldots\subset U_k=W_k\subset U_{k+1}\subset\ldots\subset U_n\}=\mathbb{O}_{Diag(\mu)}\simeq \Fmu,
\end{equation*}
and $\pi_{13}$ is the natural projection between the varieties of partial flags forgetting $V'_k.$ Applying Corollary \ref{corollary: pi_*^-} one gets
\begin{align*}
-\mathfrak{f}_{\mu-\textbf{e}_{k+1},k}&(q)\star\mathfrak{e}_{\mu-\textbf{e}_{k+1},k}(p)\\
&\!\!\!\! =-\pi_{13*}(l^{p+q})=
\left\{
\begin{array}{rl}
(-1)^{\mu_{k+1}}[\Lambda_{\mu_{k+1}}(\mathcal{T}_{k+1})][Sym_{(p+q)-\mu_{k+1}}(\mathcal{T}_{k+1})],& p+q\ge \mu_{k+1},\\
0,& 0<p<\mu_{k+1},\\
-\lbrack Sym_{-(p+q)}(\mathcal{T}_{k+1}^*)],& p+q\le 0.\\
\end{array}
\right.
\end{align*}
\end{proof}

\begin{lemma}\label{lemma: local H as commutator spec 2} 
Let $1\le k<n$ be an integer and $\mu$ be a composition of $d$ of length $n$ such that $\mu_k>0$ and $\mu_{k+1}=0.$ Then the class 
\begin{equation*}
\mathfrak{h}_{\mu,k}(p+q):=\mathfrak{e}_{\mu-\textbf{e}_k,k}(p)\star\mathfrak{f}_{\mu-\mathfrak{e}_k,k}(q)
\end{equation*}
is supported on $\mathbb{O}_{Diag(\mu)}\subset \Fnd\times \Fnd$ and only depends on the sum $p+q.$ Furthermore, 

\begin{equation*}
\mathfrak{h}_{\mu,k}(p+q)=
\left\{
\begin{array}{rl}
(-1)^{\mu_k-1}[\Lambda_{\mu_k}(\mathcal{T}^*_k)][Sym_{-(p+q)-\mu_k}(\mathcal{T}^*_k)],& p+q\le -\mu_k,\\
0,& -\mu_k<p+q<0,\\
\left[Sym_{p+q}(\mathcal{T}_k)\right],& p+q\ge 0,
\end{array}
\right.
\end{equation*}
where $T_k$ is the tautological bundles over $\mathbb{O}_{Diag(\mu)}\simeq \Fmu$ corresponding to the $k$th step of the flag.
\end{lemma}

\begin{proof}
Consider the intersection of preimages 
\begin{equation*}
X:=\pi_{12}^{-1}(\mathbb{O}_{\textbf{E}(\mu-\textbf{e}_k,k)})\cap\pi_{23}^{-1}(\mathbb{O}_{\textbf{F}(\mu-\textbf{e}_k,k)})\subset \Fnd\times \Fnd\times \Fnd.
\end{equation*}
Note that according to Lemma \ref{lemma: transversality}, the intersection is transversal. Similar to Lemma \ref{lemma: local H as commutator}, we get triples of flags $(\textbf{U},\textbf{V},\textbf{W})$ satisfying the following conditions:
\begin{itemize}
\item Since $\textbf{U}\in \mathcal{F}_{\textbf{R}^{\textbf{E}(\mu-\textbf{e}_k,k)}}=\Fmu,$ we get $\dim U_i/U_{i-1}=\mu_i$ for all $1<i\le n.$ In particular, we get $U_{k+1}=U_k$ as $\mu_k=0.$
\item Since $(\textbf{U},\textbf{V})\in \mathbb{O}_{\textbf{E}(\mu-\textbf{e}_k,k)}$ we get $U_i=V_i$ for $i\neq k$ and $V_k\subset U_k,$ $\dim U_k/V_k=1.$
\item Since $(\textbf{V},\textbf{W})\in \mathbb{O}_{\textbf{F}(\mu-\textbf{e}_k,k)}$ we get $V_i=W_i$ for $i\neq k$ and $V_k\subset W_k,$ $\dim W_k/V_k=1.$ Moreover, we get $W_k=W_{k+1}=U_{k+1}=U_k.$
\end{itemize}
In other words, using Corollary \ref{corollary: convolution with restriction} with $U=\Fnd\times\Fnd$ we get:
\begin{equation*}
X=\{U_0\subset\ldots\subset U_{k-1}\subset V_k\subset U_k=W_k\subset U_{k+1}\subset\ldots\subset U_n\}\simeq \mathcal{F}_{(\mu_1,\ldots,\mu_{k-1},\mu_k-1,1,\mu_{k+1},\ldots,\mu_n)}.
\end{equation*}
Moreover,
\begin{equation*}
\pi_{12}^*(\mathfrak{e}_{\mu-\textbf{e}_k,k}(p))=m^p,
\end{equation*}
and
\begin{equation*}
\pi_{23}^*(\mathfrak{f}_{\mu-\textbf{e}_k,k}(q))=m^q,
\end{equation*}
where $m=[U_k/V_k]=[W_k/V_k]\in K^G(X).$ The projection $\pi_{13}|_X$ forgets the flag $\textbf{V}.$ One gets
\begin{equation*}
\pi_{13}(X)=\{U_0\subset\ldots\subset U_k=W_k\subset U_{k+1}\subset\ldots\subset U_n\}=\mathbb{O}_{Diag(\mu)}\simeq \Fmu,
\end{equation*}
and $\pi_{13}$ is the natural projection between the varieties of partial flags forgetting $V_k.$ Applying Corollary \ref{corollary: pi_*^+} one gets
\begin{align*}
\mathfrak{e}_{\mu-\textbf{e}_k,k}(p)&\star\mathfrak{f}_{\mu-\textbf{e}_k,k}(q)\\
&=\pi_{13*}(l^{p+q})=
\left\{
\begin{array}{rl}
[Sym_{p+q}(\mathcal{T}_k],& p+q\ge 0,\\
0,& -\mu_k<p<0,\\
(-1)^{\mu_k-1}[\Lambda_{\mu_k}(\mathcal{T}_k^*)][Sym_{-(p+q)-\mu_k}(\mathcal{T}_k^*)],& p+q\le -\mu_k.
\end{array}
\right.
\end{align*}
\end{proof}

\begin{definition}\label{def: H=0 for muk=muk+1=0}
Let $1\le k<n$ be an integer and $\mu$ be a composition of $d$ of length $n$ such that $\mu_k=\mu_{k+1}=0.$ Then we set
\begin{equation*}
\mathfrak{h}_{\mu,k}:=0\in K^G(\mathbb{O}_{Diag(\mu)})\subset \Snd.
\end{equation*}
\end{definition}

\begin{definition}
Let $\mu$ be a composition of $d$ of length $n,$ and let $r\in\mathbb{Z}$ be an integer. Suppose that $\mu_n>0,$ then define
\begin{equation*}
K^G(\mathbb{O}_{Diag(\mu)})\ni\mathfrak{h}_{\mu,n}(r):=
\left\{
\begin{array}{rl}
(-1)^{\mu_n-1}[\Lambda_{\mu_n}(\mathcal{T}^*_n)][Sym_{-(r)-\mu_n}(\mathcal{T}^*_n)],& r\le -\mu_n,\\
0,& -\mu_n<r<0,\\
\left[Sym_r(\mathcal{T}_n)\right],& r\ge 0.
\end{array}
\right.
\end{equation*}
If $\mu_n=0$ then set
\begin{equation*}
K^G(\mathbb{O}_{Diag(\mu)})\ni\mathfrak{h}_{\mu,n}(r)=0.
\end{equation*}
\end{definition}

Consider the embedding $i:\Fnd\times \Fnd\to \mathcal{F}_{n+1}^d\times \mathcal{F}_{n+1}^d$ extending every flag $\textbf{U}$ in the trivial way $U_{n+1}=U_n.$ It is not hard to see that $i$ is closed and respects the convolution product, therefore one gets that $i_*:\Snd\to \mathbb{S}^{\operatorname{aff}}_0(n+1,d)$ is an embedding of algebras. Then, according to Lemma \ref{lemma: local H as commutator spec 2}, we get
\begin{equation*}
i_*(\mathfrak{h}_{\mu,n}(p+q))=\mathfrak{h}_{\bar{\mu},n}(p+q)=\mathfrak{e}_{\bar{\mu}-\textbf{e}_n,n}(p)\star\mathfrak{f}_{\bar{\mu}-\textbf{e}_n,n}(q),
\end{equation*}
where $\bar{\mu}:=(\mu_1,\ldots,\mu_n,0).$ This allows one to deal with $\mathfrak{h}_{\mu,n}(r)$ in the same way as with $\mathfrak{h}_{\mu,k}(r)$ for $1\le k<n.$

Consider the diagonal $D\subset \Fnd\times \Fnd$ consisting of pairs of identical flags: $D:=\{(\textbf{U},\textbf{U})\}.$ Clearly, $D$ is closed and invariant, therefore we get $K^G(D)\subset K^G(\Fnd\times \Fnd).$

\begin{lemma}\label{lemma: local H conv=tensor}
The subspace $K^G(D)\subset K^G(\Fnd\times \Fnd)$ is closed under the convolution product and, moreover, the convolution product on $K^G(D)$ coincides with the usual tensor product.
\end{lemma}

\begin{proof}
Note that the intersection of preimages $X:=\pi^{-1}_{12}(D)\cap\pi^{-1}_{23}(D)\subset \Fnd\times \Fnd\times \Fnd$ is transversal according to Lemma \ref{lemma: transversality} and consists of triples of identical flags. Therefore, one gets
\begin{equation*}
\pi^{-1}_{12}(D)\cap\pi^{-1}_{23}(D)\simeq D\simeq \Fnd,
\end{equation*}
and the restrictions of the projections $\pi_{12}|_X,$ $\pi_{23}|_X,$ and $\pi_{13}|_X$ are all natural isomorphisms. In particular, using Corollary \ref{corollary: convolution with restriction} with $U=\Fnd\times\Fnd,$ for any $\alpha,\beta\in K^G(D)$ one gets
\begin{equation*}
\alpha\star\beta=\pi_{13*}\left(\pi^*_{12}(\alpha)\times\pi^*_{23}(\beta)\right)=\alpha\times\beta.
\end{equation*}
\end{proof}

One immediately gets the following Corollary:

\begin{corollary}\label{cor: local H commutes}
One gets $\mathfrak{h}_{\mu,k}(p)\star\mathfrak{h}_{\nu,l}(q)=0$ unless $\mu=\nu,$ and
\begin{equation*}
\mathfrak{h}_{\mu,k}(p)\star\mathfrak{h}_{\mu,l}(q)=\mathfrak{h}_{\mu,l}(q)\star\mathfrak{h}_{\mu,k}(p).
\end{equation*}
\end{corollary}

\begin{proof}
Indeed, all the classes $\mathfrak{h}_{\mu,k}(p)$ are supported on the diagonal $D\subset \Fnd\times \Fnd.$ Therefore, the convolution product coincides with the tensor product, which is commutative in $K^G(\Fnd\times \Fnd)$. Also, for the first statement, note that the supports of the classes $\mathfrak{h}_{\mu,k}(p)$ and $\mathfrak{h}_{\nu,l}(q)$ are disjoint for $\mu\neq\nu.$
%The first statement is clear because the intersection of the corresponding preimages is empty:
%$$
%\pi_{12}^{-1}(\mathbb{O}_{Diag(\mu)})\cup \pi_{23}^{-1}(\mathbb{O}_{Diag(\nu)})=\emptyset.
%$$
%The second statement is also straightforward as both intersections of the preimages coincide with the triples $(\textbf{U},\textbf{U},\textbf{U})$ of identical flags, where $\textbf{U}\in \Fmu,$ and the pull-backs of the tautological bundles are isomorphic:
%\begin{equation*}
%\pi_{12}^*(T_k)\simeq\pi_{23}^*(T_k)\ \textrm{and}\ \pi_{12}^*(T_l)\simeq\pi_{23}^*(T_l).
%\end{equation*} 
%(Here $T_k$ and $T_l$ are the tautological vector bundles on $\mathbb{O}_{Diag(\mu)}\simeq \Fmu$ corresponding to the $k$th and the $l$th steps respectively.)
\end{proof}

\begin{lemma}\label{lemma: local H_n vs E and F}
Let $\mu$ be a composition of $d-1$ of length $n.$ One gets the following relations:
\begin{equation*}
\mathfrak{h}_{\mu+\textbf{e}_{n-1},n}(p)\star\mathfrak{e}_{\mu,n-1}(q)=\mathfrak{e}_{\mu,n-1}(q)\star\mathfrak{h}_{\mu+\textbf{e}_n,n}(p)-\mathfrak{e}_{\mu,n-1}(q+1)\star\mathfrak{h}_{\mu+\textbf{e}_n,n}(p-1),
\end{equation*}
for all $p,q\in\mathbb{Z}.$
\begin{equation*}
\mathfrak{h}_{\mu+\textbf{e}_k,n}(p)\star\mathfrak{e}_{\mu,k}(q)=\mathfrak{e}_{\mu,k}(q)\star\mathfrak{h}_{\mu+\textbf{e}_{k+1},n}(p),
\end{equation*}
for all $1\le k<n-1$ and $p,q\in\mathbb{Z}.$

Similarly, for classes $\mathfrak{f}:$
\begin{equation*}
\mathfrak{f}_{\mu,n-1}(p)\star\mathfrak{h}_{\mu+\textbf{e}_{n-1},n}(q)=\mathfrak{h}_{\mu+\textbf{e}_n,n}(q)\star\mathfrak{f}_{\mu,n-1}(p)-\mathfrak{h}_{\mu+\textbf{e}_n,n}(q-1)\star\mathfrak{f}_{\mu,n-1}(p+1),
\end{equation*}
for all $p,q\in\mathbb{Z}.$
\begin{equation*}
\mathfrak{h}_{\mu+\textbf{e}_{k+1},n}(p)\star\mathfrak{f}_{\mu,k}(q)=\mathfrak{f}_{\mu,k}(q)\star\mathfrak{h}_{\mu+\textbf{e}_k,n}(p),
\end{equation*}
for all $1\le k<n-1$ and $p,q\in\mathbb{Z}.$
\end{lemma}

\begin{proof}
All four formulas are proved in a similar manner using the embedding $i_*:\Snd\to \mathbb{S}^{\operatorname{aff}}_0(n+1,d)$ described above, definition of elements $\mathfrak{h}$ and the relations on the elements $\mathfrak{e}$ and $\mathfrak{f}.$ Let us prove the first formula. Applying the embedding $i_*$ and the definition of $\mathfrak{h}_{\bar{\mu}+\textbf{e}_{n-1},n}(p)$ to the left hand side, one gets (note that $\bar{d}_{n+1}=0$):
\begin{align*}
\mathfrak{h}_{\bar{\mu}+\textbf{e}_{n-1},n}&(p)\star\mathfrak{e}_{\bar{\mu},n-1}(q)=\mathfrak{e}_{\bar{\mu}+\textbf{e}_{n-1}-\textbf{e}_n,n}(p)\star\mathfrak{f}_{\bar{\mu}+\textbf{e}_{n-1}-\textbf{e}_n,n}(0)\star\mathfrak{e}_{\bar{\mu},n-1}(q)\\
&=\mathfrak{e}_{\bar{\mu}+\textbf{e}_{n-1}-\textbf{e}_n,n}(p)\star\mathfrak{e}_{\bar{\mu}+\textbf{e}_{n+1}-\textbf{e}_n,n-1}(q)\star\mathfrak{f}_{\bar{\mu},n}(0)\\
&=\left(\mathfrak{e}_{\bar{\mu},n-1}(q)\star\mathfrak{e}_{\bar{\mu},n}(p)-\mathfrak{e}_{\bar{\mu},n-1}(q+1)\star\mathfrak{e}_{\bar{\mu},n}(p-1)\right)\star\mathfrak{f}_{\bar{\mu},n}(0)\\
&=\mathfrak{e}_{\bar{\mu},n-1}(q)\star\left(\mathfrak{e}_{\bar{\mu},n}(p)\star\mathfrak{f}_{\bar{\mu},n}(0)\right)-\mathfrak{e}_{\bar{\mu},n-1}(q+1)\star\left(\mathfrak{e}_{\bar{\mu},n}(p-1)\star\mathfrak{f}_{\bar{\mu},n}(0)\right)\\
&=\mathfrak{e}_{\bar{\mu},n-1}(q)\star\mathfrak{h}_{\bar{\mu}+\textbf{e}_n,n}(p)-\mathfrak{e}_{\bar{\mu},n-1}(q+1)\star\mathfrak{h}_{\bar{\mu}+\textbf{e}_n,n}(p-1)\\
&=i_*\left(\mathfrak{e}_{\mu,n-1}(q)\star\mathfrak{h}_{\mu+\textbf{e}_n,n}(p)-\mathfrak{e}_{\mu,n-1}(q+1)\star\mathfrak{h}_{\mu+\textbf{e}_n,n}(p-1)\right).
\end{align*}
Since $i_*$ is injective, this completes the proof.
\end{proof}

\begin{lemma}\label{lemma: local H and E rel}
Let $\mu$ be a composition of $d-1$ of length $n,$ and let $1\le k<n.$ Then one gets:
\begin{equation*}
\mathfrak{h}_{\mu+\textbf{e}_k,k}(p)\star\mathfrak{e}_{\mu,k}(q)=-\mathfrak{e}_{\mu,k}(q-1)\star\mathfrak{h}_{\mu+\textbf{e}_{k+1},k}(p+1).
\end{equation*}
\end{lemma}

\begin{proof}
Consider the intersection of preimages
\begin{equation*}
X:=\pi_{12}^{-1}(\mathbb{O}_{Diag(\mu+\textbf{e}_k)})\cap\pi_{23}^{-1}(\mathbb{O}_{\textbf{E}(\mu,k)})\subset \Fnd\times \Fnd\times \Fnd.
\end{equation*}
Note that according to Lemma \ref{lemma: transversality}, the intersection is transversal. Furthermore, it consists of triples of flags $(\textbf{U},\textbf{U},\textbf{V})$ satisfying conditions:
\begin{itemize}
\item Since $\textbf{U}\in \mathcal{F}_{\mu+\textbf{e}_k},$ one gets $\dim V_i/V_{i-1}=\mu_i$ for $i\neq k,$ and $\dim U_k/U_{k-1}=\mu_k+1.$
\item Since $(\textbf{U},\textbf{V})\in \mathbb{O}_{\textbf{E}(\mu,k)}$ one gets $V_i=U_i$ for $i\neq k$ and $V_k\subset U_k,$ $\dim U_k/V_k=1.$
\end{itemize}
In other words, one gets $X\simeq \mathcal{F}_{(\mu_1,\ldots,\mu_k,1,\mu_{k+1},\ldots,\mu_n)}.$
The restriction of the projection $\pi_{13}$ to $X$ is an isomorphism to $\mathbb{O}_{\textbf{E}(\mu,k)}\simeq \mathcal{F}_{(\mu_1,\ldots,\mu_k,1,\mu_{k+1},\ldots,\mu_n)}.$ 

Similarly, consider the intersection of preimages
\begin{equation*}
Y:=\pi_{12}^{-1}(\mathbb{O}_{\textbf{E}(\mu,k)})\cap\pi_{23}^{-1}(\mathbb{O}_{ Diag(\mu+\textbf{e}_{k+1})})\subset \Fnd\times \Fnd\times \Fnd.
\end{equation*}
It consists of triples of flags $(\textbf{U},\textbf{V},\textbf{V}),$ where the flags $\textbf{U}$ and $\textbf{V}$ satisfy the same conditions as before. Therefore, one gets $Y\simeq \mathcal{F}_{(\mu_1,\ldots,\mu_k,1,\mu_{k+1},\ldots,\mu_n)},$ and the restriction of the projection $\pi_{13}$ to $Y$ is also an isomorphism to $\mathbb{O}_{\textbf{E}(\mu,k)}\simeq \mathcal{F}_{(\mu_1,\ldots,\mu_k,1,\mu_{k+1},\ldots,\mu_n)}.$ 

Let $\mathcal{T}_k$ and $\mathcal{T}_{k+1}$ be the tautological vector bundles on $\mathbb{O}_{\textbf{E}(\mu,k)}$ corresponding to the steps $\mu_k$ and $\mu_{k+1}$ respectively, and let $\mathcal{L}$ be the line bundle corresponding to the one dimensional step, so that $[\mathcal{L}]=\mathfrak{e}_{\mu,k}(1)\in K^G(\mathbb{O}_{\textbf{E}(\mu,k)})$. Then, using Corollary \ref{corollary: convolution with restriction} with $U=\Fnd\times\Fnd,$ one gets
\begin{enumerate}
\item If $p\le -(\mu_k+1) \Leftrightarrow p+1\le -\mu_k,$ then
\begin{align*}
\mathfrak{h}_{\mu+\textbf{e}_k,k}(p)\star\mathfrak{e}_{\mu,k}(q)&=
(-1)^{\mu_k}[\mathcal{L}]^q[\Lambda_{\mu_k+1}(\mathcal{T}^*_k\oplus \mathcal{L}^*)][Sym_{-p-(\mu_k+1)}(\mathcal{T}^*_{k+1}\oplus \mathcal{T}^*_k\oplus \mathcal{L}^*)]\\
&=-(-1)^{\mu_k-1}[\mathcal{L}]^{q-1}[\Lambda_{\mu_k}(\mathcal{T}^*_k)][Sym_{-(p+1)-\mu_k}(\mathcal{T}^*_{k+1}\oplus \mathcal{T}^*_k\oplus \mathcal{L}^*)]\\
&=-\mathfrak{e}_{\mu,k}(q-1)\star\mathfrak{h}_{\mu+\textbf{e}_{k+1},k}(p+1).
\end{align*}
\item If $-(\mu_k+1)<p<\mu_{k+1}\Leftrightarrow -\mu_k<p+1<\mu_{k+1}+1,$ then
\begin{equation*}
\mathfrak{h}_{\mu+\textbf{e}_k,k}(p)\star\mathfrak{e}_{\mu,k}(q)=0=-\mathfrak{e}_{\mu,k}(q+1)\star\mathfrak{h}_{\mu+\textbf{e}_{k+1},k}(p-1).
\end{equation*}
\item If $p\ge \mu_{k+1}\Leftrightarrow p+1\ge \mu_{k+1}+1,$ then
\begin{align*}
\mathfrak{h}_{\mu+\textbf{e}_k,k}(p)\star\mathfrak{e}_{\mu,k}(q)&=
(-1)^{\mu_{k+1}}[\mathcal{L}]^q[\Lambda_{\mu_{k+1}}(\mathcal{T}_{k+1})][Sym_{p-\mu_{k+1}}(\mathcal{T}_{k+1}\oplus \mathcal{T}_k\oplus \mathcal{L})]\\
&=-(-1)^{\mu_{k+1}+1}[\mathcal{L}]^{q-1}[\Lambda_{\mu_{k+1}+1}(\mathcal{T}_{k+1}\oplus \mathcal{L})][Sym_{(p+1)-(\mu_{k+1}+1)}(\mathcal{T}_{k+1}\oplus \mathcal{T}_k\oplus \mathcal{L})]\\
&=-\mathfrak{e}_{\mu,k}(q-1)\star\mathfrak{h}_{\mu+\textbf{e}_{k+1},k}(p+1).
\end{align*}
\end{enumerate}
\end{proof}

\begin{lemma}
Let $\mu$ be a composition of $d-1$ of length $n,$ and let $1\le k<n.$ Then one gets:
\begin{equation*}
\mathfrak{h}_{\mu+\textbf{e}_{k+1},k}(p)\star\mathfrak{f}_{\mu,k}(q)=-\mathfrak{f}_{\mu,k}(q+1)\star\mathfrak{h}_{\mu+\textbf{e}_k,k}(p-1).
\end{equation*}
\end{lemma}

\begin{proof}
We skip the proof as it is very similar to the proof of Lemma \ref{lemma: local H and E rel}.
\end{proof}

\begin{theorem}\label{theorem: Snd generated by E F H}
The algebra $\Snd=K^G(\Fnd\times \Fnd)$ is generated by the classes $\mathfrak{e}_{\mu,k}(p),\ \mathfrak{f}_{\mu,k}(p),$ and $\mathfrak{h}_{\nu,n}(p),$ where $\mu$ runs through all compositions of $d-1$ of length $n,$ $\nu$ run through all compositions of $d$ of length $n,$ $k$ runs through $\{1,\ldots, n-1\},$ and $p$ runs through $\mathbb{Z}.$ 
\end{theorem}

\begin{proof}
According to Theorem \ref{theorem: Snd generated by diag and almost diag}, it is enough to show that by combining classes $\mathfrak{e}_{\mu,k}(p),\ \mathfrak{f}_{\mu,k}(p),$ and $\mathfrak{h}_{\nu,n}(p)$ one can obtain any class supported on an orbit $\mathbb{O}_M,$ where $M$ is either diagonal or almost diagonal.

Let us start with the case when $M$ is diagonal. Let $\nu=(\nu_1,\ldots,\nu_n)$ be a composition of $d.$ According to equation \ref{equation: K^G(O_M)}, we have an isomorphism
\begin{equation*}
K^G(\mathbb{O}_{Diag(\nu)})\simeq\mathbb{C}[x_1^{\pm 1},\ldots,x_d^{\pm 1}]^{S_{\nu}}\simeq\bigotimes_{k=1}^n \Lambda^{\pm}[\textbf{x}_\textbf{k}],
\end{equation*}
where $\textbf{x}_\textbf{k}:=(x_{(\sum_{i=1}^{k-1} \nu_i)+1},\ldots,x_{\sum_{i=1}^k \nu_i}),$ and $\Lambda^{\pm}[\textbf{x}_\textbf{k}]$ is the ring of symmetric Laurent polynomials of $\textbf{x}_\textbf{k}.$ We will prove by induction that elements $\mathfrak{h}_{\nu,l}(r),$ $r\in\mathbb{Z},$ $1\le l\le n,$ generate $\Lambda^{\pm}[\textbf{x}_{\textbf{k}}]$ for all $1\le k\le n$ starting the induction from $k=n$ and proceeding in the decreasing order. Note that according to Corollary \ref{lemma: local H conv=tensor} the convolution product of classes supported on the diagonal coincides with the usual tensor product (i.e. the ordinary product of the corresponding partially symmetric functions).

Under the above isomorphism one has $\mathfrak{h}_{\nu,n}(r)=h_r(\textbf{x}_\textbf{n})$ for $r\ge 0,$ and 
$$
\mathfrak{h}_{\nu,n}(-\nu_n)=(-1)^{\nu_n-1}e_{\nu_n}(\textbf{x}_\textbf{n}^{-1})=\frac{(-1)^{\nu_n-1}}{x_{(\sum_{i=1}^{k-1} \nu_i)+1}\cdot\ldots\cdot x_{\sum_{i=1}^k \nu_i}}.
$$ 
It follows that $\Lambda^{\pm}(\textbf{x}_\textbf{n})$ is generated by $\mathfrak{h}_{\nu,n}(r),$ $r\in\{-\nu_n,0,1,\ldots\}.$

Suppose now that we proved that $\Lambda^{\pm}[\textbf{x}_\textbf{i}]$ are generated by $\mathfrak{h}_{\nu,l}(r)$ for $i>k.$ Then for $r\ge \nu_{k+1}$ one has 
$$
\mathfrak{h}_{\nu,k}(r)=(-1)^{\nu_{k+1}}e_{\nu_{k+1}}(\textbf{x}_{\textbf{k+1}})h_{r-\nu_{k+1}}(\textbf{x}_{\textbf{k}},\textbf{x}_{\textbf{k+1}}),
$$
%Since $e_{\nu_{k+1}}(\textbf{x}_{\textbf{k+1}}^{-1})\in\Lambda^{pm}[\textbf{x}_{\textbf{k+1}}],$ we conclude that then $h_i(\textbf{x}_{\textbf{k}},\textbf{x}_{\textbf{k+1}})$ can all be generated, and therefore all the positive degree symmetric polynomials in $(\textbf{x}_{\textbf{k}},\textbf{x}_{\textbf{k+1}})$ can too.
where $h_{r-\nu_{k+1}}(\textbf{x}_{\textbf{k}},\textbf{x}_{\textbf{k+1}})$ is the complete homogeneous symmetric polynomial in the union of variable from $\textbf{x}_{\textbf{k}}$ and $\textbf{x}_{\textbf{k+1}}.$ For $r=-\nu_k$ one has
$$
\mathfrak{h}_{\nu,k}(-\nu_k)=(-1)^{\nu_k-1}e_{\nu_k}(\textbf{x}_\textbf{k}^{-1}).
$$
Since $\Lambda^{\pm}[\textbf{x}_\textbf{k+1}]$ is already proven to be generated by elements $\mathfrak{h}_{\nu,l}(r)$ and 
$$
e_{\nu_k+\nu_{k+1}}(\textbf{x}_\textbf{k}^{-1},\textbf{x}_\textbf{k+1}^{-1})=e_{\nu_k}(\textbf{x}_\textbf{k}^{-1})e_{\nu_{k+1}}(\textbf{x}_\textbf{k+1}^{-1}),
$$ 
we conclude that $\Lambda^{\pm}[\textbf{x}_\textbf{k},\textbf{x}_\textbf{k+1}]$ is generated too. Finally, we use the elementary fact that $\Lambda^{\pm}[\textbf{x}_\textbf{k},\textbf{x}_\textbf{k+1}]$ and $\Lambda^{\pm}[\textbf{x}_\textbf{k+1}]$ together generate $\Lambda^{\pm}[\textbf{x}_\textbf{k}]\otimes\Lambda^{\pm}[\textbf{x}_\textbf{k+1}].$

Let now $\mu=(\mu_1,\ldots,\mu_n)$ be a composition of $d-1$ and $1\le k<n$ be an integer. Consider the orbit $\mathbb{O}_{\textbf{E}(\mu,k)}\simeq \mathcal{F}_{(\mu_1,\ldots,\mu_k,1,\mu_{k+1},\ldots,\mu_n)}.$ According to equation \ref{equation: K^G(O_M)}, we have an isomorphism
\begin{equation*}
K^G(\mathbb{O}_{\textbf{E}(\mu,k)})\simeq\mathbb{C}[x_1^{\pm 1},\ldots,x_d^{\pm 1},y^{\pm 1}]^{S_{\mu_1,\ldots,\mu_n,1}}\simeq\left(\bigotimes_{k=1}^n \Lambda^{\pm}[\textbf{x}_\textbf{k}]\right)\otimes \Lambda^{\pm}[y],
\end{equation*}
where for each $k$ the group of variables $\textbf{x}_\textbf{k}$ corresponds to the step of dimension $\mu_k,$ and $y$ corresponds to the off diagonal entry $1$ in $\textbf{E}(\mu,k),$ i.e. $y=\mathfrak{e}_{\mu,k}(1).$ We get that the subspace $\Lambda^{\pm}[y]\subset K^G(\mathbb{O}_{\textbf{E}(\mu,k)})$ is generated by $\mathfrak{e}_{\mu,k}(p),$ $p\in \mathbb{Z}.$ Let 
$$
\nu=\mu+\textbf{e}_k=(\mu_1,\ldots,\mu_k+1,\ldots,\mu_n),
$$ 
and consider the convolution product
\begin{equation*}
\star:K^G(\mathbb{O}_{Diag(\nu)})\otimes \Lambda^{\pm}[y]\to K^G(\mathbb{O}_{\textbf{E}(\mu,k)}).
\end{equation*}
Note that similar to the proof of Lemma \ref{lemma: local H and E rel}, the restriction of the projection $\pi_{13}:\Fnd\times \Fnd\times \Fnd \to \Fnd\times \Fnd$ to the intersection of the preimages $\pi_{12}^{-1}(\mathbb{O}_{Diag(\nu)})\cap\pi_{23}^{-1}(\mathbb{O}_{\textbf{E}(\mu,k)})$ is an isomorphism to $\mathbb{O}_{\textbf{E}(\mu,k)}\subset \Fnd\times \Fnd.$ Therefore, in terms of the partial symmetric functions the convolution product is simply the multiplication map:
\begin{equation*}
\left(\Lambda^{\pm}[\textbf{x}_\textbf{k},y]\otimes\bigotimes_{i\neq k} \Lambda^{\pm}[\textbf{x}_\textbf{i}]\right)\otimes \Lambda^{\pm}[y]\to \left(\bigotimes_{k=1}^n \Lambda^{\pm}[\textbf{x}_\textbf{k}]\right)\otimes \Lambda^{\pm}[y],
\end{equation*}
which is clearly surjective.

The case when $M$ is a lower-triangular almost diagonal matrix is done in a similar manner.
\end{proof}

\section{The map $\phi_d:\Unz\to \Snd$}\label{section: the map phi_d}

\subsection{Definition of the map $\phi_d$}\label{subsection: definition of phi_d}

Recall that we use notations $\mu\vDash d,$ $l(\mu)=n$ for a composition of $d$ of length $n.$ All compositions in this section are of length $n.$

\begin{definition}
Fix positive integers $n$ and $d,$ and for $1\le k<n$ and $p\in\mathbb{Z}$ set
\begin{equation*}
\mathfrak{e}_k(p):=\sum_{\mu \vDash d-1}\mathfrak{e}_{\mu,k}(p),
\end{equation*}
and 
\begin{equation*}
\mathfrak{f}_k(p):=\sum_{\mu \vDash d-1} \mathfrak{f}_{\mu,k}(p).
\end{equation*}
%where the sums are taken over all compositions $\mu$ of $d-1$ of length $n.$
Also, for $1\le k\le n$ set 
\begin{equation*}
\mathfrak{h}_k(p):=\sum_{\nu \vDash d} \mathfrak{h}_{\nu,k}(p).
\end{equation*}
%where the sum is taken over all compositions $\nu$ of $d$ of length $n.$
\end{definition}

\begin{theorem}
The map sending $E_k(p)\to\mathfrak{e}_k(p),$ $F_k(p)\to \mathfrak{f}_k(p),$ and $H_n(p)\to \mathfrak{h}_n(p)$ for all $1\le k<n$ and $p\in\mathbb{Z}$ extends to a well defined homomorphism of algebras $\phi_d:\Unz\to \Snd.$
\end{theorem}

\begin{proof}
All we need to do is to check that the elements $\mathfrak{e}_k(p),$ $\mathfrak{f}_k(p),$ and $\mathfrak{h}_n(p)$ satisfy the relations on the generators $E_k(p),$ $F_k(p),$ and $H_n(p)$ of the algebra $\Unz.$ This is done by a direct application of the relations on the local generators $\mathfrak{e}_{\mu,k}(p),$ $\mathfrak{f}_{\mu,k}(p),$ and $\mathfrak{h}_{\nu,n}(p)$ developed in the previous section. Since all the relations are proved in a similar and rather routine way, we prove relation \ref{rel:3 term E} as an example and leave the rest as an exercise.
%\begin{itemize}
%\item Let us start with relation \ref{rel:2 term E}. We use Lemmas \ref{lemma: local gens prod=0 unless} and \ref{Lemma: local relation E two terms} to get
%\begin{align*}
%\mathfrak{e}_k(p)&\star\mathfrak{e}_k(q)=\sum_{\mu \vDash d-1} \mathfrak{e}_{\mu,k}(p)\star\sum_{\nu \vDash d-1} \mathfrak{e}_{\nu,k}(q)=\sum_{\textbf{C}^{\textbf{E}(\mu,k)}=\textbf{R}^{\textbf{E}(\nu,k)}} \mathfrak{e}_{\mu,k}(p)\star\mathfrak{e}_{\nu,k}(q)\\
%&=\sum_{\mu+\textbf{e}_{k+1}=\nu+\textbf{e}_k} \mathfrak{e}_{\mu,k}(p)\star\mathfrak{e}_{\nu,k}(q)=\sum_{\mu\vDash d,\ \mu_k>0,\ \mu_{k+1}>0} \mathfrak{e}_{\mu-\textbf{e}_{k+1},k}(p)\star\mathfrak{e}_{\mu-\textbf{e}_k,k}(q)\\
%&=-\sum_{\mu \vDash d,\ \mu_k>0,\ \mu_{k+1}>0} \mathfrak{e}_{\mu-\textbf{e}_{k+1},k}(q-1)\star\mathfrak{e}_{\mu-\textbf{e}_k,k}(p+1)=-\mathfrak{e}_k(q-1)\star\mathfrak{e}_k(p+1).
%\end{align*}
%\item 
For Relation \ref{rel:3 term E} we start by using  Lemmas \ref{lemma: local gens prod=0 unless} and \ref{Lemma: local relation E three terms}:
\begin{align*}
\mathfrak{e}_{k+1}(p)\star\mathfrak{e}_k(q)&=\sum_{\mu\vDash d-1} \mathfrak{e}_{\mu,k+1}(p)\star\sum_{\nu\vDash d-1} \mathfrak{e}_{\nu,k}(q)=\sum_{\textbf{C}^{\textbf{E}(\mu,k+1)}=\textbf{R}^{\textbf{E}(\nu,k)}} \mathfrak{e}_{\mu,k+1}(p)\star\mathfrak{e}_{\nu,k}(q)\\
&=\sum_{\mu+\textbf{e}_{k+2}=\nu+\textbf{e}_k} \mathfrak{e}_{\mu,k+1}(p)\star\mathfrak{e}_{\nu,k}(q)=\sum_{\mu\vDash d-2} \mathfrak{e}_{\mu+\textbf{e}_k,k+1}(p)\star\mathfrak{e}_{\mu+\textbf{e}_{k+2},k}(q)\\
&=\sum_{\mu\vDash d-2} \left(\mathfrak{e}_{\mu+\textbf{e}_{k+1},k}(q)\star\mathfrak{e}_{\mu+\textbf{e}_{k+1},k+1}(p)\right.\\
&\qquad\qquad\qquad\qquad\qquad\left.-\mathfrak{e}_{\mu+\textbf{e}_{k+1},k}(q+1)\star\mathfrak{e}_{\mu+\textbf{e}_{k+1},k+1}(p-1)\right).
\end{align*}
Recall that according to Lemma \ref{lemma: local E three term special} one has
\begin{equation*}
\mathfrak{e}_{\mu,k}(q)\star\mathfrak{e}_{\mu,k+1}(p)-\textbf{e}_{\mu,k}(q+1)\star\mathfrak{e}_{\mu,k+1}(p-1)=0,
\end{equation*}
for any composition $\mu\vDash d-1$ such that $\mu_{k+1}=0.$ In particular, one gets
\begin{equation*}
\sum_{\mu\vDash d-1,\mu_{k+1}=0}\left(\mathfrak{e}_{\mu,k}(q)\star\mathfrak{e}_{\mu,k+1}(p)-\mathfrak{e}_{\mu,k}(q+1)\star\mathfrak{e}_{\mu,k+1}(p-1)\right)=0.
\end{equation*}
Adding this, one continues:
\begin{align*}
&\sum_{\mu\vDash d-2}\left( \mathfrak{e}_{\mu+\textbf{e}_{k+1},k}(q)\star\mathfrak{e}_{\mu+\textbf{e}_{k+1},k+1}(p)-\mathfrak{e}_{\mu+\textbf{e}_{k+1},k}(q+1)\star\mathfrak{e}_{\mu+\textbf{e}_{k+1},k+1}(p-1)\right)\\
&=\sum_{\mu\vDash d-1} \mathfrak{e}_{\mu,k}(q)\star\mathfrak{e}_{\mu,k+1}(p)-\sum_{\mu\vDash d-1}\mathfrak{e}_{\mu,k}(q+1)\star\mathfrak{e}_{\mu,k+1}(p-1)\\
&=\sum_{\mu\vDash d-1} \mathfrak{e}_{\mu,k}(q)\star\mathfrak{e}_{\mu,k+1}(p)-\sum_{\mu\vDash d-1}\mathfrak{e}_{\mu,k}(q+1)\star\mathfrak{e}_{\mu,k+1}(p-1)\\
&=\sum_{\textbf{C}^{\textbf{E}(\mu,k)}=\textbf{R}^{\textbf{E}(\nu,k+1)}} \mathfrak{e}_{\mu,k}(q)\star\mathfrak{e}_{\nu,k+1}(p)-\sum_{\textbf{C}^{\textbf{E}(\mu,k)}=\textbf{R}^{\textbf{E}(\nu,k+1)}}\mathfrak{e}_{\mu,k}(q+1)\star\mathfrak{e}_{\nu,k+1}(p-1)\\
&=\mathfrak{e}_k(q)\star\mathfrak{e}_{k+1}(p)-\mathfrak{e}_k(q+1)\star\mathfrak{e}_{k+1}(p-1).
\end{align*}

\end{proof}

\subsection{Surjectivity of $\phi_d$}\label{subsection: surjectivity of phi_d}

\begin{theorem}
The homomorphism $\phi_d:\Unz\to \Snd$ is surjective.
\end{theorem}

\begin{proof}
Fix a composition $\mu\vDash d,$ of length $n.$ According to the proof of Theorem \ref{theorem: Snd generated by E F H}, classes $\mathfrak{h}_{\mu,k}(p),$ $1\le k\le n,$ $p\in\mathbb{Z},$ generate $K^G(\mathbb{O}_{Diag(\mu)})\subset \Snd.$ It then follows that there exists an element $a_\mu\in\phi_d(\Unz)$ supported on the union of orbits corresponding to the diagonal matrices, and such that $a_\mu|_{\mathbb{O}_{Diag(\mu)}}=[\det V].$ (Note that $\det V$ is trivial as an ordinary vector bundle, but carries a non-trivial action of $G$). In fact, $a_\mu$ can be expressed as a polynomial in $\mathfrak{h}_k(p),$  $1\le k\le n,$ $p\in\mathbb{Z}.$ Recall that according to Lemma \ref{lemma: local H conv=tensor} the convolution product of the elements supported on the diagonal coincides with the usual tensor product. 

Consider the element 

\begin{equation*}
b_\mu:=\mathfrak{h}_1(-\mu_1)\star\mathfrak{h}_2(-\mu_2)\star\ldots\star\mathfrak{h}_n(-\mu_n)\in\phi_d(\Unz)\subset \Snd.
\end{equation*}

We claim that $b_{\mu}$ is supported on $\mathbb{O}_{Diag(\mu)}.$ Indeed, it is clear that $b_\mu$ is supported on the orbits corresponding to the diagonal matrices. Let $\nu\vDash d$ be a composition of length $n,$ different from $\mu.$ Then there exists $1\le k\le n$ such that $\nu_k>\mu_k.$ Therefore,

\begin{equation*}
\mathfrak{h}_k(-\mu_k)|_{\mathbb{O}_{Diag(\nu)}}=\mathfrak{h}_{\nu,k}(-\mu_k)=0, 
\end{equation*}
and
\begin{equation*}
b_\mu|_{\mathbb{O}_{Diag(\nu)}}=0.
\end{equation*} 

Furthermore,
\begin{equation*}
b_\mu|_{\mathbb{O}_{Diag(\mu)}}=\prod_{k=1}^n[\det T_k^*]=[\det V^*]=[\det V]^{-1},
\end{equation*} 
where for every $1\le k\le n,$ $T_k$ is the tautological vector bundle on $\mathbb{O}_{Diag(\mu)}\simeq \Fmu$ corresponding to the $k$th step. We conclude that the product $a_\mu\times b_\mu$ is supported on $\mathbb{O}_{Diag(\mu)}$ and that $(a_\mu\times b_\mu)|_{\mathbb{O}_{Diag(\mu)}}=1.$ Furthermore, for any composition $\mu\vDash d$ of length $n$ and any $p\in\mathbb{Z}$ on gets
\begin{equation*}
(a_\mu\times b_\mu)\star\mathfrak{h}_n(p)=\mathfrak{h}_{\mu,n}(p)\in \phi_d(\Unz),\\
\end{equation*}
and for any composition $\mu\vDash d-1$ of length $n,$ any $1\le k< n,$ and any $p\in\mathbb{Z}$ on gets
\begin{align*}
(a_{\mu+\textbf{e}_k}\times b_{\mu+\textbf{e}_k})\star\mathfrak{e}_k(p)&=\mathfrak{e}_{\mu,k}(p)\in \phi_d(\Unz),\\
(a_{\mu+\textbf{e}_{k+1}}\times b_{\mu+\textbf{e}_{k+1}})\star\mathfrak{f}_k(p)&=\mathfrak{f}_{\mu,k}(p)\in \phi_d(\Unz).
\end{align*}
We apply Theorem \ref{theorem: Snd generated by E F H} to conclude the proof.
\end{proof}

\end{document}